\documentclass[onefignum,onetabnum]{siamart190516}

\usepackage{accents}
\usepackage[toc,page]{appendix}
\usepackage{enumerate}
\usepackage{amssymb}

\newsiamthm{lmm}{Lemma}
\newsiamthm{thrm}{Theorem}
\newsiamthm{rmrk}{remark}
\newsiamthm{dfntn}{Definition}
\newsiamthm{system}{System}
\newsiamthm{assumption}{Assumption}
\newsiamthm{crllr}{Corollary}

\newenvironment{matrixbj}{\begin{bmatrix}}{\end{bmatrix}} 
\newcommand{\bv}[1]{{\mathbf #1}} 
\newcommand{\bt}[1]{{\mathbf{#1}}}
\newcommand{\scalco}[1]{{\mathfrak{#1}}} 

\newcommand{\ubar}[1]{\underaccent{\bar}{#1}}
\newcommand{\ubv}[1]{\ubar{\bv{#1}}}
\newcommand{\funSpace}[1]{\mathcal{#1}} 
\newcommand{\Lscal}{\langle}	\newcommand{\Rscal}{\rangle} 

\newcommand{\LegDens}{g}
\newcommand{\LegPDE}{\mathcal{G}}
\newcommand{\HamPDE}{\mathcal{H}}
\newcommand{\MassPDE}{\mathcal{M}}
\newcommand{\LegCo}{{G}}
\newcommand{\HamCo}{{H}}

 \newcommand{\blsPipe}{{{\Omega}}}
 \newcommand{\blsVertex}{{{\nu}}} 
\newcommand{\derst}[1]{{\nabla_{#1}}} 
 \newcommand{\TraceOp}{{\mathcal{T}}}
 \newcommand{\Onepipe}{\omega}
 
 \newcommand{\bcef}{\bv{e}} 
  \newcommand{\bcfl}{\bv{f}} 
  \newcommand{\nbn}{p} 
    
\usepackage{tikz}
\usetikzlibrary{calc,positioning,shapes,arrows,decorations.shapes} 
\tikzstyle{rect}=[
   rectangle, draw=black, font=\small\sffamily\bfseries, inner sep=9pt]

\begin{document}

\title{On port-Hamiltonian approximation of a \\nonlinear flow problem on networks\thanks{
The support of the German Federal Ministry for Economic Affairs and Energy (BMWI) via the project \textit{MathEnergy} and the German Federal Ministry of Education and Research (BMBF) via the project \textit{EiFer} is acknowledged.}}

\author{Bj\"orn Liljegren-Sailer\thanks{Universit\"at Trier, FB IV - Mathematik, Lehrstuhl Modellierung und Numerik, D-54286 Trier, Germany
  (\email{Corresponding author: bjoern.sailer@uni-trier.de}).}
\and Nicole Marheineke\footnotemark[2]}

\maketitle
\date{\today}
\begin{abstract} 
This paper deals with the systematic development of structure-preserving approximations for a class of nonlinear partial differential equations on networks. The class includes, for example, gas pipe network systems described by barotropic Euler equations. Our approach is guided throughout by energy-based modeling concepts (port-Hamiltonian formalism, theory of Legendre transformation), which provide a convenient and general line of reasoning. Under mild assumptions on the approximation, local conservation of mass, an energy bound, and the inheritance of the port-Hamiltonian structure can be shown. Our approach is not limited to conventional space discretization but also covers complexity reduction of the nonlinearities by inexact integration. Thus, it can serve as a basis for structure-preserving model reduction. Combined with an energy-stable time integration, we demonstrate the applicability and good stability properties using the example of the Euler equations on networks.
\end{abstract}

%

\section{Introduction}

Structure-preserving discretization is an active research area in the last decades.
By preserving or mimicking relevant geometric structures such as, e.g., conservation laws, dissipative relations, or symplecticities, unphysical solution behavior and numerical instabilities can be avoided in many cases, cf.\ \cite{book:leveque-finite-volume-methods,art:topics-in-strpresdisc,art:avfPDE, art:morParamHamHesthaven}. 
The model problem we consider in this paper describes nonlinear flows on networks. It covers a hierarchy of models used to describe gas network systems, including particularly the barotropic Euler equations \cite{art:hierarchGasMindt2019, art:domschke-adj-based2015, art:egger-mfem-compressEuler,inproc:bls-hyp2018, inbook:mathEnergy21}, but also p-systems \cite{inbook:ruggeri05-dissiphyp,art:mei-psystems09} and more general symmetrizable hyperbolic systems \cite{art:harten1983,art:mock1980}. 
Further applications are in the contexts \cite{art:weak-pH-koty, inbook:Maschke2001}, e.g., in the modeling of electric transmission lines. 
The network is assumed to be described by a directed graph. Each edge $\Onepipe$ of the graph can be identified with an interval. Given a strictly convex smooth function $h: \mathbb{R}^2 \rightarrow \mathbb{R}$ and a non-negative function $\tilde{r}: \mathbb{R}^2 \rightarrow \mathbb{R}$, the edgewise states $\ubar{\bv{z}}^\Onepipe = [z_1^\Onepipe;z_2^\Onepipe]: [0,T] \times \omega \rightarrow \mathbb{R}^2$ are governed by
\begin{align*} 
	\partial_t z_1^\Onepipe(t,x) &= - \partial_x \nabla_{2} h(\ubar{\bv{z}}^\Onepipe(t,x)), \\
	\partial_t z_2^\Onepipe(t,x) &= - \partial_x \nabla_{1} h(\ubar{\bv{z}}^\Onepipe(t,x)) - \tilde{r}(\ubar{\bv{z}}^\Onepipe(t,x)) \nabla_{2} h(\ubar{\bv{z}}^\Onepipe(t,x)),
\end{align*}
with $\derst{i} h(\ubar{\bv{z}}^\Onepipe(t,x)) = \partial_{z_i} h( [z_1;z_2])_{|[z_1;z_2] = \ubar{\bv{z}}^\Onepipe(t,x)}$ for $i=1,2$. The expressions $\MassPDE(\ubar{\bv{z}}) = \sum_{\Onepipe \in \mathcal{E}} \int_{\Onepipe} z_1^\Onepipe dx$ and $\tilde{\HamPDE}(\ubar{\bv{z}}) = \sum_{\Onepipe \in \mathcal{E}} \int_{\Onepipe}  h(\ubar{\bv{z}}^\Onepipe) dx$, with $\mathcal{E}$ the set of all edges, represent the total mass and the Hamiltonian of the system. Fundamental properties of the hyperbolic model problem are that, under appropriate coupling conditions on the edgewise equations, conservation of mass and dissipation of the Hamiltonian (energy dissipation) hold up to the exchange with the boundary. Moreover, the convective terms can be related to a certain skew symmetric geometric structure. Our aim is to derive a structure-preserving approximation approach, which is applicable for both conventional space discretization and projection-based (Galerkin) model order reduction with additional complexity reduction.
The latter denotes a sparse approximation of nonlinear terms (also known as hyper-reduction). Note that we use the expression model reduction as an umbrella term for model order- and complexity-reduction.

There exists a rich literature on structure-preserving discretization. The mimetic finite difference and finite element methods \cite{book:arnold2006-comp, art:lee-mixededmimFE} are designed to mimic conservation laws and dissipative relations on the discrete level. These approaches typically involve the separate approximation of operators and geometric objects and do therefore not completely fit in the Galerkin framework, or are tailored towards specific discrete structures or problems, see e.g., \cite{art:egger-mfem-compressEuler, art:euler-reigstad, art:lilsailer-dwemor, inproc:bls-hyp2018} in the gas network context. Many entropy-stable methods are based on a flux-centered point of view and the so-called entropy-flux pair. This is the case for finite volume methods \cite{book:leveque-finite-volume-methods, art:Jameson2008} or methods using summation-by-parts approximation \cite{art:FISHER13-entrop, art:Nordstrom14-SumByParts}. These methods are particularly well-suited in the presence of shocks. The summation-by-parts methods have the additional advantage that inexact integration can be taken into account, but they are not purely projection-based and heavily rely on the flux-interpretation, which complicates their adaption for model reduction. A contribution in this direction we are aware of is \cite{art:chan-EntropROM}. For systems in Lagrangian form, structure-preserving Galerkin approximations are studied in \cite{art:hyperreduction-preserving-lagrangian-structure, art:comred-ecsw}. The Hamiltonian formulation is dual to the Lagrangian one \cite{art:topics-in-strpresdisc, art:gieselmann-rel-energy}, and the so-called Legendre transformation links the Hamiltonian and the Lagrangian functions of a system \cite{book:RockWets98,book:Zeidler1984NonlinearFA}. 

For the approximation of our model problem, we consider the port-Hamiltonian framework, a generalization of the Hamiltonian formalism, which is particularly well-suited for network problems.
It has its origins in the analysis of finite-dimensional connected systems \cite{art:pcH-maschke92, inbook:Maschke2001, art:SchM13, art:SchJ14} but has been extended to the infinite-dimensional setting of partial differential equations \cite{art:badlyan2018open, art:GoloTSM04} and has been systematically generalized to constrained dynamical systems \cite{art:beattie2017porthamiltonian, art:SchJ14}. The lumped port-Hamiltonian approximation of infinite-dimensional systems has also gained interest recently. Most contributions in this direction either focus on linear systems \cite{inbook:interpolation-based-port-Hamiltonian-Systems,art:WOLF2010401} or are tailored to very specific discrete structures such as mimetic discretizations \cite{art:FISHER13-entrop, inproc:farle-pH-FE,art:PasumarthyAS12}. 
There exist a few works on structure-preserving model reduction of nonlinear port-Hamiltonian systems, most of them consider the finite-dimensional case, see \cite{art-mor-structure-preserving-nonlinear-pH,art:morParamHamHesthaven,art:moramHamDissHesthaven,art:symplHamMor}. A comprehensive literature is available for standard model reduction. 
cf., e.g., \cite{art:morBauBF14, book:dimred2003, art:comred-ecsw, art:hyperreduction-preserving-lagrangian-structure}. 

In this paper, we propose a conforming Galerkin ansatz, which respects the Hamiltonian and geometric structure of our nonlinear flow problem. On the one hand, we make use of compatibility conditions on the ansatz spaces, which also play a fundamental role in mimetic finite element methods and symplectic model reduction \cite{art:symplHamMor,art:morParamHamHesthaven}. On the other hand, a certain variable transformation related to the Hamiltonian of the system and its systematic analysis using the theory on the partial Legendre transformation is crucial in our approach. To the best of the authors' knowledge, the latter has not yet been used in this extend in the context. The use of structured variable transformations widens the range of formulations, for which structure-preserving Galerkin approximations can be derived with purely variational arguments. Moreover, complexity reduction of nonlinear terms by inexact quadrature can be included without much difficulty, similar to other Galerkin-based approaches \cite{art:discret-dissip-Egger19, art:comred-ecsw}. 
We present a complete modeling work flow consisting of the following steps: Port-Hamiltonian modeling of our flow problem on networks; analytical investigation of a variable transformation induced by the partial Legendre transformation; approximation in space by Galerkin projection and quadrature-based complexity reduction (hyper-reduction);  energy-stable time discretization.

The structure of this paper is as follows: The underlying energy-based modeling concepts, i.e., the port-Hamilto\-nian formalism and the partial Legendre transformation, are shortly introduced in Section~\ref{sec:centr-concepts}. In Section~\ref{sec-p1b-netproblem} we present our model problem together with an appropriate parametrization of the solution and a variational principle. The latter provides the basis for the structure-preserving approximations derived in Section~\ref{sec:approx-framew}. The approximations inherit, among others, port-Hamiltonian structure, which is revealed by their structured coordinate representations in Section~\ref{sec:Part1b:coordrepres}. An energy-stable time discretization is presented in Section~\ref{sec:p1b:time-disc} and the applicability of our approach is numerically demonstrated for the barotropic Euler equations in two test cases in Section~\ref{chap:p1b:isothEuler}.
Our conclusion and outlook give a summary and point towards possible directions for future research. The appendix extends our approach to systems with an additional dissipation term as well as with edge weighting.

\subsection*{Notation}
In this paper matrices, vectors and scalars are indicated by capital boldfaced, small boldfaced and normal letters, respectively, whereby vector-ex{\-}pres{\-}sions always refer to column vectors. Given two vectors $\bv a \in \mathbb{R}^n$ and $\bv b \in \mathbb{R}^m$, we write $[\bv a ; \bv b] \in \mathbb{R}^{n+m}$ for the concatenation to a new (column) vector. For $m=n$ the Euclidean scalar product is denoted by $\bv{a} \cdot \bv{b}$.
We distinguish between function-valued vector spaces (e.g.,~$\funSpace{L}^2$, $\funSpace{V}$) and  real-valued spaces and sets (e.g.,~$\mathbb{R}^N, \mathbb{S}$) in the typesetting.

Given a scalar field $h:\mathbb{R}^n \rightarrow \mathbb{R}$ and a partitioning of the argument into vector-components ${\bv{z}}=[{\bv{z}}_1;{\bv{z}}_2]$ with ${\bv{z}}_i \in \mathbb{R}^{n_i}$, $n_1 +n_2 =n$, we often structure the gradient $\nabla h$ and Hessian $\nabla^2 h$ in sub-blocks, i.e.,
\begin{align*}
		 \nabla h = \begin{matrixbj}
			\nabla_{1} h \\
			\nabla_{2} h
\end{matrixbj},\quad
\nabla^2 h =
		\begin{matrixbj}
			\nabla_{11} h & \nabla_{12} h \\
			\nabla_{21} h & \nabla_{22} h 
\end{matrixbj}, \hspace{0.3cm}
 \nabla_i h:\mathbb{R}^n \rightarrow \mathbb{R}^{n_i}, \quad  \nabla_{ij} h:\mathbb{R}^n \rightarrow \mathbb{R}^{{n_i \times n_j}}.
\end{align*}

The solutions of the partial differential equations in this paper depend on time $t$ and space $x$, e.g., $z(t,x) \in \mathbb{R}$. When it is convenient, we interpret them as functions in time with values in a function space, e.g., $z(t) \in \funSpace{L}^2$. The grouping of solution components in a vector is underscored to distinguish it from an ordinary vector (coordinate representation).


\section{Modeling concepts for structured systems} \label{sec:centr-concepts}

In this section, we briefly present some basic concepts for characterizing properties in a structured system. The concepts are applicable on the continuous level of our model problem as well as on the space-discrete level of its approximations. The core is the port-Hamiltonian framework that encodes geometric structures in an algebraic way. The geometric structures typically reflect fundamental physical properties, hence it is desirable to preserve them throughout all approximation steps. The partial Legendre transformation allows us to systematically investigate certain variable transformations related to the Hamiltonian of the system.

\subsection{Port-Hamiltonian systems} \label{sec:finitdim-pH}
In recent years, the port-Hamiltonian framework has been found in many different areas of applications, see, e.g., \cite{inbook:Maschke2001, art:SchJ14, art:GoloTSM04, art:beattie2017porthamiltonian, art:weak-pH-koty}. Central is the Hamiltonian that is often related to an energy or entropy function in applications. In the finite-dimensional setting, $\HamCo:\mathbb{R}^n \rightarrow \mathbb{R}$ is typically assumed to be convex and continuously differentiable. An important class of port-Hamiltonian systems with state $\bv{z}\in \mathbb{R}^n$ reads as follows, \cite{art:SchJ14,art:beattie2017porthamiltonian}: 

Find $\bv{z} \in \funSpace{C}^1([0,T]; \mathbb{R}^{n})$, $\bcef \in \funSpace{C}([0,T]; \mathbb{R}^{p})$ and $\bcfl \in \funSpace{C}^1([0,T]; \mathbb{R}^{p})$ such that
\begin{align} 
\begin{aligned}\label{eq:pH-first-prototype}
\frac{d}{dt}\bv{z}(t)
&= \left( \bar{\bt{J}}(\bv{z}(t)) - \bar{\bt{R}}(\bv{z}(t)) \right)
 \nabla_{} \HamCo(\bv{z}(t)) 
+ \bt{K}
 \bcef(t) \qquad &   \\
  \bcfl(t)   &= 
  \bt{K}^T  \nabla_{} \HamCo(\bv{z}(t)) , \hspace{2cm} \bv{z}(0)= \bv{z}_0. 
  \end{aligned}
\end{align}
Here, the matrix $\bar{\bt{J}}(\bv{z})$ is anti-symmetric (i.e., $\bar{\bt{J}}(\bv{z}) = -\bar{\bt{J}}(\bv{z})^T$), $\bar{\bt{R}}(\bv{z})$ is symmetric positive semi-definite (i.e., $\bar{\bt{R}}(\bv{z}) = \bar{\bt{R}}(\bv{z})^T$ has non-negative eigenvalues), and \mbox{$\bt{K}\in \mathbb{R}^{n \times p}$.} 
The choice of closing conditions depends on the application. Given an input $\bv{u}:[0,T]\rightarrow \mathbb{R}^p$, the system~\eqref{eq:pH-first-prototype} can, e.g., be closed by the equations $\bcef(t) = \bv{u}(t)$, or $\bcfl(t) = \bv{u}(t)$. In the first case, the system reduces to an ordinary differential equation for $\bv{z}$, whereas the second case yields a differential-algebraic equation of index 2. In the port-Hamiltonian wording, $\bv{z}$ is called the energy variable, $\nabla_{} H(\bv{z})$ the effort variable, and $\bcef$ and $\bcfl$ the boundary effort and boundary flow, respectively. The system structure readily implies, among others, that the Hamiltonian is dissipated over time up to exchange with the boundary, i.e., for $t \geq 0$
\begin{align*} 
	\frac{d}{dt} \HamCo(\bv{z}(t)) \leq \bcef(t) \cdot \bcfl(t).
\end{align*}

\subsection{Partial Legendre transformation} \label{sec:part-Leg-traf}
The Legendre transformation plays an important role in classical mechanics, as it represents a link between the Lagrangian and the Hamiltonian modeling framework \cite{book:Zeidler1984NonlinearFA}. Further, there exists a rich theory on the implication of the Legendre transformation for duality principles, which are, e.g.,  used in nonlinear optimization \cite{rockafellar-1970a}. In this paper we deal with a slight generalization of the standard approach, the so-called partial Legendre transformation \cite{book:RockWets98}.
Accordingly, we consider the following partitioning of a vector, $\bv{z} = [\bv{z}_1; \bv{z}_2]\in \mathbb{R}^n$ with $\bv{z}_i \in \mathbb{R}^{n_i}$, $n = n_1+n_2$.

\begin{dfntn}[Partial Legendre transformation] \label{def:partialLeg-GEN}
	Let $\mathbb{S} = \mathbb{S}_1 \times \mathbb{S}_2\subset \mathbb{R}^n$ be a convex set with $\mathbb{S}_i\subset \mathbb{R}^{n_i}$, $n=n_1+n_2$. The partial Legendre transformation of the function $h:\mathbb{S} \to \mathbb{R}$ with respect to the second (sub-vector) component is defined as $g:\mathbb{S} \to \mathbb{R}$, 
\begin{align*}
	 \LegDens(\bv{a}) = \sup_{\bv{z}_2\in \{\bar{\bv{z}}_2: [\bv{a}_1;\bar{\bv{z}}_2] \in \mathbb{S}\}} \bv{a}_2 \cdot \bv{z}_2  - h([\bv{a}_1 ;\bv{z}_2]), \qquad \bv{a} = [\bv{a}_1; \bv{a}_2].
\end{align*}	
\end{dfntn}
For differentiable functions a coordinate transformation can be associated to the Legendre transformation.

\begin{thrm} \label{theor:partialLeg-GEN}
Let $\mathbb{S}$, $h$ and $g$ be given as in Definition~\ref{def:partialLeg-GEN}. Let, additionally, $\bv{z}_2 \mapsto h([\bar{\bv{z}}_1; \bv{z}_2])$ be strictly convex and continuously differentiable for fixed $\bar{\bv{z}}_1\in \mathbb{S}_1$. Then, $\bv{a}_2 \mapsto g([\bar{\bv{z}}_1; \bv{a}_2])$ is strictly convex, differentiable and $\hat{\bv{z}}: D(\hat{\bv{z}})\rightarrow \mathbb{S}$ with 
\begin{align*}
		\hat{\bv{z}}(\bv{a}) = 
			[\bv{a}_1;
			\nabla_2 g(\bv{a})]
		, \qquad \text{for } \bv{a} = [\bv{a}_1;\bv{a}_2]
		 \text{ and }  D(\hat{\bv{z}}) = \{ \bv a \in \mathbb{R}^n:\, \hat{\bv{z}}(\bv a) \in \mathbb{S} \}
\end{align*}
is a homeomorphism. Moreover, it holds $g(\bv{a}) = \bv a_2 \cdot \nabla_2 g(\bv a) - h([\bv a_1; \nabla_2 g(\bv a)])$.
 \end{thrm}

\begin{proof}
The function $\phi:\bv{a}_2 \mapsto g([\bar{\bv{z}}_1; \bv{a}_2])$ can be considered as the Legendre transformation of $\bv{z}_2 \mapsto h([\bar{\bv{z}}_1; \bv{z}_2])$ for any fixed $\bar{\bv{z}}_1$. According to \cite[Theorem 11.13]{book:RockWets98}, thus $\phi$ inherits strict convexity and continuous differentiability and the gradients of the two mappings are inverse functions of each other, i.e., $\nabla_{2} h(\bv{z})_{|\bv{z}=\hat{\bv{z}}([\bar{\bv{z}}_1;\bv{a}_2])} = \bv{a}_2$ for $[\bar{\bv{z}}_1;{\bv{a}}_2] \in D(\hat{\bv{z}})$ holds. For bijectivity of $\hat{\bv{z}}$, it remains to show injectivity, i.e., $\hat{\bv{z}}(\breve{\bv{a}}) \neq \hat{\bv{z}}(\tilde{\bv{a}})$ holds for any $\breve{\bv{a}}, \tilde{\bv{a}}\in D(\hat{\bv{z}})$ with $\breve{\bv{a}} \neq \tilde{\bv{a}}$. If the first (sub-vector) components differ ($\breve{\bv{a}}_1 \neq \tilde{\bv{a}}_1$), this holds trivially. Let us therefore assume $\breve{\bv{a}}_1 = \tilde{\bv{a}}_1$ and $\breve{\bv{a}}_2 \neq \tilde{\bv{a}}_2$. From the strict convexity of $\phi$,
	$ \left( \nabla_{2} g([\breve{\bv{a}}_1; \breve{\bv{a}}_2]) - \nabla_{2} g([\breve{\bv{a}}_1; \tilde{\bv{a}}_2]) \right) \cdot (\breve{\bv{a}}_2 - \tilde{\bv{a}}_2 ) > 0$ follows,
	hence $\nabla_{2} g(\breve{\bv{a}}) \neq \nabla_{2} g(\tilde{\bv{a}})$ and, consequently, also $\hat{\bv{z}}(\breve{\bv{a}}) \neq \hat{\bv{z}}(\tilde{\bv{a}})$. 
Finally, the fact that  $\phi$ is a Legendre transformation of a parametrized function also yields the equality $g(\bv{a}) = \bv a_2 \cdot \nabla_2 g(\bv a) - h([\bv a_1; \nabla_2 g(\bv a)])$ under the posed smoothness assumptions, cf.~\cite{rockafellar-1970a}.
\end{proof}
By considering $\bv{a}_2 \mapsto g([\bar{\bv{z}}_1; \bv{a}_2])$ as a function parametrized in $\bar{\bv{z}}_1$, the function $h$ can be characterized as the partial Legendre transformation of $g$ by standard results on Legendre transformations, cf.~\cite{book:RockWets98,rockafellar-1970a}.
\begin{lmm}\label{lem:legOfleg}
Under the assumptions of Theorem~\ref{theor:partialLeg-GEN}, the relations $h(\hat{\bv z}(\bv a)) = \nabla_2 g({\bv{a}}) \cdot \bv{a}_2  - \LegDens({\bv{a}})$ and $\nabla_{1} \LegDens(\bv{a}) = - \nabla_{1} h(\bv{z})_{|\bv{z} = \hat{\bv{z}}(\bv{a})}$ are valid for $\bv a = [\bv{a}_1;\bv{a}_2]$.
\end{lmm}


\section{Model problem} \label{sec-p1b-netproblem}

A class of nonlinear flows on networks is described by our model problem that covers, e.g., the barotropic Euler equations or electromagnetic waves.  We write it in a specific form using a variable transformation that is related to a partial Legendre transformation. This is advantageous since parts of the coupling conditions become linear, which significantly simplifies the handling of network aspects. Moreover, it is key in the analysis of a variational principle. We show that the variational principle encodes an energy bound, local mass conservation and inherits a gradient structure strongly connected to the Hamiltonian of our model problem. 

\subsection{Network description} \label{subsec-netdescrip}
 
\begin{figure}[t] 
\begin{center}
\begin{tabular}{l | l}

\begin{minipage}{.24\textwidth}
\hspace*{-2.4em}
{%
\begin{tikzpicture}[scale=1.7]
\node[circle,draw,inner sep=2pt] (v1) at (-1.87,0) {$\nu_1$};
\node[circle,draw,inner sep=2pt] (v2) at (-0.87,0) {$\nu_2$};
\node[circle,draw,inner sep=2pt] (v3) at (0,0.5) {$\nu_3$};
\node[circle,draw,inner sep=2pt] (v4) at (0,-0.5) {$\nu_4$};
\draw[->,thick,line width=1.5pt] (v1) -- node[above] {$\Onepipe_1$} ++(v2);
\draw[->,thick,line width=1.5pt] (v2) -- node[above,sloped] {$\Onepipe_2$} ++(v3);
\draw[->,thick,line width=1.5pt] (v2) -- node[above,sloped] {$\Onepipe_3$} ++(v4);
\end{tikzpicture}
} 
\end{minipage}
&
\begin{minipage}[c]{.62\textwidth}
\includegraphics[width= 1\textwidth]{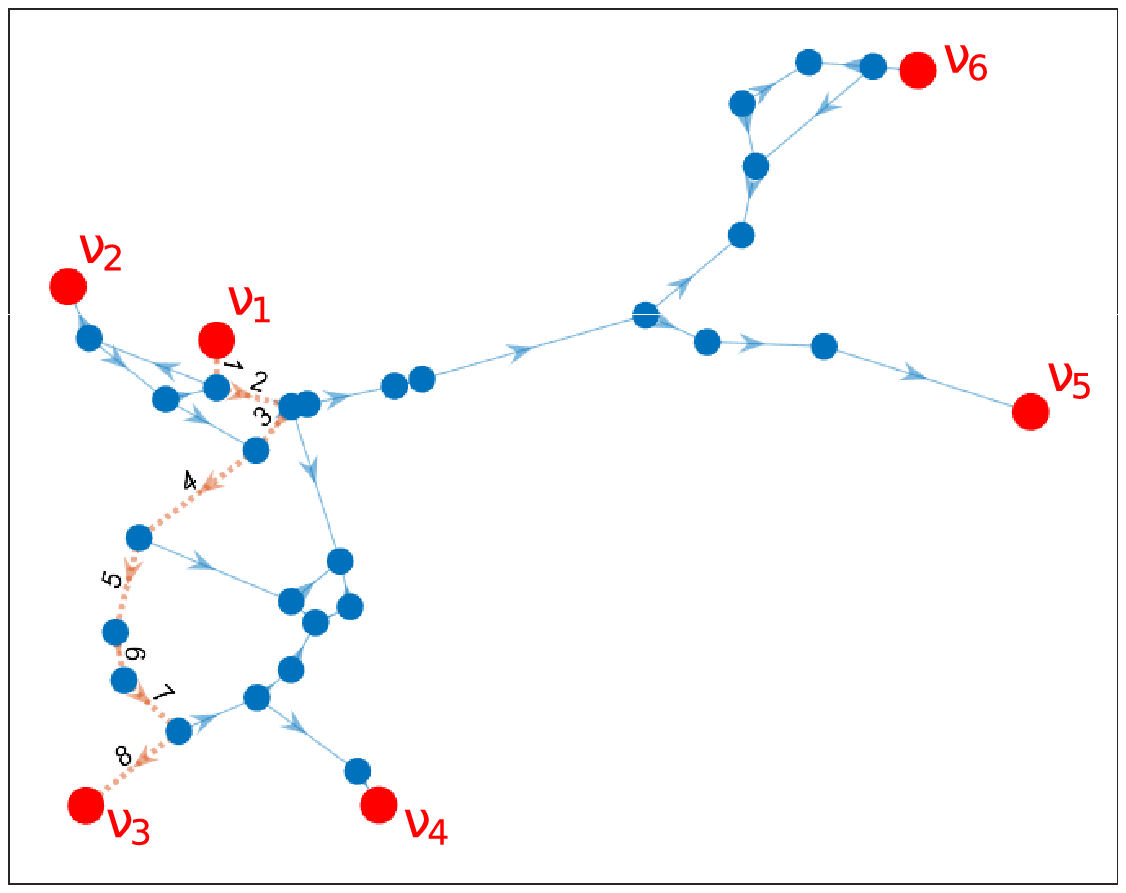}
\end{minipage}
\end{tabular}
\end{center}
\caption{
\textit{Left:} Illustration of graph notation for $\mathcal{G}=(\mathcal{N},\mathcal{E})$ with nodes $\mathcal{N}=\{\nu_1,\nu_2,\nu_3,\nu_4\}$ and edges $\mathcal{E}=\{\Onepipe_1,\Onepipe_2,\Onepipe_3\}$
defined by $\Onepipe_1=(\nu_1,\nu_2)$, $\Onepipe_2=(\nu_2,\nu_3)$ and $\Onepipe_3=(\nu_2,\nu_4)$. 
Thus, $\mathcal{N}_0~=\{\nu_2\}$, $\mathcal{N}_\partial~=~\{\nu_1,\nu_3,\nu_4\}$ and $\mathcal{E}(\nu_2)=\{\Onepipe_1,\Onepipe_2,\Onepipe_3\}$.
Further, $n^{\Onepipe_1}[\nu_1]=n^{\Onepipe_2}[\nu_2]=n^{\Onepipe_3}[\nu_2]=1$ and 
$n^{\Onepipe_1}[\nu_2]=n^{\Onepipe_2}[\nu_3]=n^{\Onepipe_3}[\nu_4]=-1$.\newline
\textit{Right:} Topology for gas pipeline network used in the numerical studies (Section~\ref{subsec:realize-num}). The larger red nodes $\nu_i$, $i=1,...,6$, correspond to boundary nodes. The spatial domain of the pipes $\omega_j$, $j=1,...,8$, i.e., the path from $\nu_1$ to $\nu_3$, is highlighted by a dotted brown line.
\label{fig:net-topol-ex}}
\end{figure}

Concerning the network aspects, we rely on the framework of \cite{art:lilsailer-dwemor,art:kugler-dwe-net}. Let
a network be described by a directed graph $\mathcal{G} = (\mathcal{N}, {\mathcal{E}} )$ with set of nodes ${\mathcal{N}} = \{\blsVertex_1, \ldots, \blsVertex_l\}$ and edges $\mathcal{E} = \{ \Onepipe_1, \ldots, \Onepipe_k \} \subset {\mathcal{N}} \times {\mathcal{N}} $.
To every edge $\Onepipe$, we associate a length $\ell^\Onepipe$. The set of all edges adjacent to the node $\blsVertex$ is denoted by $\mathcal{E}(\blsVertex) =  \{\Onepipe \in \mathcal{E}: \, \Onepipe= ( \blsVertex, \bar{\blsVertex}), \text{ or } \Onepipe=( \bar{\blsVertex}, \blsVertex) \}$. For $\Onepipe \in \mathcal{E}(\nu)$, the incidence mapping $n^\Onepipe[\blsVertex]$ is defined by
\begin{align*}
	 n^\Onepipe[\blsVertex] =
	 \begin{cases}
	 	\,\,\,\, 1 & \text{ for } \Onepipe = ( \blsVertex, \bar{\blsVertex}) \text{ for some }  \bar{\blsVertex} \in \mathcal{N} \\
	 	 -1  & \text{ for } \Onepipe = ( \bar{\blsVertex}, \blsVertex ) \text{ for some }  \bar{\blsVertex} \in \mathcal{N} .
	 \end{cases}
\end{align*}
The nodes are grouped into interior nodes ${\mathcal{N}}_0 \subset {\mathcal{N}}$ and boundary nodes \mbox{${\mathcal{N}}_\partial = {\mathcal{N}} \text{\textbackslash} {\mathcal{N}}_0$}, for an illustration of the network notation see Fig.~\ref{fig:net-topol-ex}.

On the network, function spaces are constructed by compositions of standard Sobolev spaces for every edge. The spatial domain is given as the union of edges $\Omega_{} = \{x: \, x\in \Onepipe, \text{ for } \Onepipe \in \mathcal{E}\}$. Note that every edge $\Onepipe$ can be identified with an interval $(0,l^\Onepipe)$ which is tacitly employed in the upcoming expressions. The space of square-integrable functions on $\mathcal{E}$ reads $\funSpace{L}^2(\mathcal{E}) =  \left\{b: \Omega_{}\rightarrow \mathbb{R} \text{ with } b_{|\Onepipe} \in \funSpace{L}^2(\Onepipe) \text{  for all } \Onepipe\in \mathcal{E} \right\}$, where the subscript $._{|\Onepipe}$ indicates the restriction of a function to the edge $\Onepipe$. The respective scalar product and norm read $\Lscal b, \tilde{b} \Rscal = \sum_{\Onepipe \in \mathcal{E}} \int_\Onepipe b \, \tilde{b} \, dx$ and $|| b || = \sqrt{\Lscal b, b \Rscal}$ for $b$, $\tilde{b} \in \funSpace{L}^2(\mathcal{E})$. 
The weak (broken) derivative operator for functions on the network is defined by $(\partial_x b)_{|\Onepipe} =  \partial_x b_{|\Onepipe}$ for $\Onepipe \in \mathcal{E}$. The space of functions with square-integrable weak broken derivative is given as $\funSpace{H}_{pw}^1(\mathcal{E}) = \left\{b\in \funSpace{L}^2(\mathcal{E}) :  \partial_x b \in \funSpace{L}^2(\mathcal{E}) \right\}.$ Accordingly, $\funSpace{C}_{pw}^k(\mathcal{E}) = \left\{b: \Omega_{} \rightarrow \mathbb{R} \text{ with } b_{|\Onepipe} \in \funSpace{C}^k(\Onepipe) \text{  for all } \Onepipe\in \mathcal{E}  \right\}$ denotes the space of piecewise smooth functions for $k\geq 0$.

The domains for boundary and coupling conditions are the sets of nodes $\mathcal{N}_\partial$ and $\mathcal{N}_0$. For $b \in \funSpace{H}_{pw}^1(\mathcal{E})$ we indicate node evaluations with squared brackets, i.e., $	b_{|\Onepipe}[\blsVertex] \in \mathbb{R}$ for $\blsVertex \in \mathcal{N}$. They are well-defined by means of the trace theorem, \cite{book:braess07}. Note that $b_{|\Onepipe}[\blsVertex]$ and $b_{|\tilde{\Onepipe}}[\blsVertex]$ may in general differ for $\Onepipe \neq \tilde{\Onepipe} \in \mathcal{E}$.
Following \cite{art:kugler-dwe-net}, a Sobolev space incorporating a certain coupling condition at inner nodes is defined as
\begin{align*}
	\funSpace{H}_{div}^1(\mathcal{E}) &= \{b\in \funSpace{H}^1_{pw}(\mathcal{E}) : \sum_{\Onepipe \in \mathcal{E}(\blsVertex) } n^{\Onepipe}[\blsVertex] b_{|\Onepipe}[\blsVertex] = 0, \text{ for }  \blsVertex \in {\mathcal{N}}_0 \}.
\end{align*}
The boundary nodes $\mathcal{N}_\partial = \{ \blsVertex_1, \ldots, \blsVertex_p \}$ are assumed to be connected to exactly one edge each. Thus, a boundary operator $\TraceOp^{ }: \funSpace{H}_{pw}^1(\mathcal{E}) \rightarrow \mathbb{R}^{p}$ can be defined by $\left[ \TraceOp^{ } b \right]_i =
	 n^\Onepipe[\blsVertex_i] b_{|\Onepipe}[\blsVertex_i] $  for  $\Onepipe \in \mathcal{E}(\blsVertex_i)$, $ i =1,\ldots, p $ and  $b \in H^1_{pw}(\mathcal{E})$.

\subsection{Strong form and variable transformation} \label{subsec:modelprob-full}

As model problem we consider a class of prototypical nonlinear partial differential equations for a flow on a network: The state $\ubar{\bv{z}} = [z_1;z_2]: [0,T] \times \blsPipe \rightarrow \mathbb{R}^2 $ is governed by
\begin{subequations} \label{bls-eq:abstr}
\begin{align} \label{bls-eq:abstr-a}
	\partial_t	\ubar{\bv{z}}(t,x) =
	\begin{bmatrix}
		 & -\partial_x \\
		 -\partial_x & -\tilde{r}(\ubar{\bv{z}}(t,x)	)
	\end{bmatrix}
	 \derst{} h(\ubar{\bv{z}}(t,x)), \hspace{0.8cm} x \in \blsPipe, \quad t \in (0,T]
\end{align}
with $\tilde{r}:\mathbb{R}^2\rightarrow \mathbb{R}$ such that $\tilde{r}(\ubar{\bv{z}}(t,x))\geq 0$. The solution components are interconnected by the coupling conditions
\begin{align} \label{bls-eq:abstr-coup}
	\sum_{\Onepipe \in \mathcal{E}(\blsVertex) } n^{\Onepipe}[\blsVertex] \nabla_{2} h(\ubar{\bv{z}}_{|\Onepipe}(t,\blsVertex))= 0  , \qquad
	\nabla_{1} h(\ubar{\bv{z}}_{|\Onepipe}(t,\blsVertex))  = \nabla_{1} h(\ubar{\bv{z}}_{|\tilde{\Onepipe}}(t,\blsVertex)) 
\end{align}
for $\Onepipe,\tilde{\Onepipe} \in \mathcal{E}(\blsVertex)$,
at $\nu \in \mathcal{N}_0$. The system is closed by initial and boundary conditions of the form
\begin{align}\label{bls-eq:abstr-bc}
\ubv{z}(0,x) = \ubv{z}_0(x) \quad  \text{for }  x \in \blsPipe, \hspace{1cm} \tilde{t}^\blsVertex(\ubv{z}(t,\nu),u_{\nu}(t)) = 0 
\end{align}
for $t \in [0,T]$, $\blsVertex\in \mathcal{N}_\partial$.
We particularly assume one boundary data $u_{\nu}: [0,T] \rightarrow \mathbb{R}$ per boundary node  $\blsVertex\in \mathcal{N}_\partial$ to be given, where $\tilde{t}^\blsVertex$ denotes an appropriately chosen real-valued function.
We refer to $h$ as the Hamiltonian density and to
\begin{align}
	\MassPDE(\ubar{\bv{z}}) = \Lscal z_1,1 \Rscal  = \sum_{\Onepipe \in \mathcal{E}} \int_{\Onepipe} z_1 dx , \hspace{1cm} \tilde{\HamPDE}(\ubar{\bv{z}})&= \Lscal h(\ubar{\bv{z}}), 1 \Rscal = \sum_{\Onepipe \in \mathcal{E}} \int_{\Onepipe} h(\ubar{\bv{z}}) dx 	
\end{align}
\end{subequations}
as the (total) mass and the Hamiltonian, respectively. If \eqref{bls-eq:abstr} has a strong solution, the following two relations can be shown
\begin{subequations}\label{eq:relation}
\begin{align}\label{eq:relation1}
		\frac{d}{dt}\MassPDE(\ubar{\bv{z}}) &=  \sum_{\blsVertex \in \mathcal{N}_\partial, \, \Onepipe \in \mathcal{E}(\blsVertex) } n^{\Onepipe}[\blsVertex] \derst{2} h(\ubar{\bv{z}}_{|\Onepipe}[\blsVertex]), \\\label{eq:relation2}
	\frac{d}{dt} \tilde{\HamPDE}(\ubar{\bv{z}}) & \leq \sum_{\blsVertex \in \mathcal{N}_\partial, \, \Onepipe \in \mathcal{E}(\blsVertex) } n^{\Onepipe}[\blsVertex] \derst{1} h(\ubar{\bv{z}}_{|\Onepipe}[\blsVertex]) \derst{2} h(\ubar{\bv{z}}_{|\Onepipe}[\blsVertex]) .
\end{align}
\end{subequations}
In our application \eqref{eq:relation1} has the interpretation of mass conservation, and \eqref{eq:relation2} is referred to as energy dissipation and $\ubar{\bv{z}}$ as energy variable from now on.  Note that the structural properties also depend on the use of appropriate coupling conditions, in our case \eqref{bls-eq:abstr-coup} relates to conservation of the mass and the Hamiltonian at inner nodes, cf. \cite{art:euler-reigstad,art:hierarchGasMindt2019}. 

Throughout our analysis, we assume that \eqref{bls-eq:abstr} has a smooth and unique solution. This assumption typically holds in the context of gas networks, which are governed by friction-dominated regimes \cite{art:egger2020stability,art:gugat17noblowup, inbook:ruggeri05-dissiphyp}, i.e., large $\tilde{r}(\cdot)$. Moreover, in many applications there is an additional dissipation term, see Appendix~\ref{subsec:dissip}. 
Our approximation scheme is constructed such that it mimics \eqref{eq:relation} on a discrete level, which enhances its stability in comparison to standard methods. The type of boundary conditions has no special influence on the proposed approximation ansatz, but many choices can be included very naturally. 

\begin{assumption} \label{assum:h-sconvex}
The domain $\mathbb{S} \subset \mathbb{R}^2$ of the Hamiltonian density $h:\mathbb{S}\rightarrow \mathbb{R}$ is an open convex set. Moreover, $h$ is twice continuously differentiable with symmetric positive definite Hessian $\nabla^2 h(\bv{z})$ for all $\bv{z} \in \mathbb{S}$.
\end{assumption}

\begin{lmm}\label{corr:h-sconvex}
Let $g$ be the partial Legendre transformation of the Hamiltonian density $h$ with respect to the second component. Under Assumption~\ref{assum:h-sconvex}, the variable transformation
	\begin{align*}
		\hat{\bv{z}}: D(\hat{\bv{z}})\rightarrow \mathbb{S}, \qquad \hat{\bv{z}}: \, {\bv{a}} \mapsto 
		[a_1; \derst{2}g({\bv{a}})] \quad \text{ with } \bv{a} = [a_1; a_2]
	\end{align*}
and its inverse are continuously differentiable.
\end{lmm}
\begin{proof}
	Theorem~\ref{theor:partialLeg-GEN} is applicable, which implies that
	\begin{align*}
		\hat{\bv{a}}:\, {\bv{z}} \mapsto 
			[z_1;
			\derst{2}h({\bv{z}})],  \quad \text{ for } \bv{z} = 
			[z_1; z_2] \in \mathbb{S}
	\end{align*}	
is the inverse of $\hat{\bv{z}}$. By Assumption~\ref{assum:h-sconvex} the derivative $\frac{d}{d\bv{z}} \hat{\bv{a}}(\bv{z})$ is continuous, and a direct calculation shows that its eigenvalues $\lambda_1 = 1$ and $\lambda_2 =  \derst{22}h({{\bv{z}}})$ are strictly positive due to the positive definiteness of $\nabla^2 h ({{\bv{z}}})$. Therefore, the derivative of $\hat{\bv{z}}$ is also a continuous function and reads
	\begin{align*}
		\frac{d}{d\bv{a}}\hat{\bv{z}}(\bv{a})  = \left(\frac{d}{d\bv{z}}\hat{\bv{a}}(\bv{z})^{-1}\right)_{|\bv{z}=\hat{\bv{z}}(\bv{a}) } = 	\begin{matrixbj}
		1 & 0 \\
		-\frac{ \nabla_{12} h(\hat{\bv{z}}(\bv{a}))} {\nabla_{22} h(\hat{\bv{z}}(\bv{a}))} & \frac{1}{\nabla_{22} h(\hat{\bv{z}}(\bv{a}))}
	\end{matrixbj},
\end{align*}		
which follows from applying the implicit function theorem on $\hat{\bv{a}}$ locally for each $\bv{z} \in \mathbb{S}$.
\end{proof}
Because of Assumption~\ref{assum:h-sconvex}, our model problem \eqref{bls-eq:abstr} is of strictly hyperbolic type. The underlying system matrix has one positive and one negative eigenvalue, such that one boundary condition per boundary node is required for a well-posed setup. We refer to \cite{phd:liljegren,book:leveque-finite-volume-methods} for details.
Change of variables have played a crucial role in the theoretical and numerical analysis of hyperbolic systems, see e.g., \cite{art:mock1980,art:harten1983,art:CharlotHughes,art:Jameson2008}. The variable-transformed formulation considered in this work is related to the partial Legendre transformation of the Hamiltonian density.

\begin{crllr}[Partial Legendre-transformed strong form] \label{corol:partial-legtrafo-sys}
Let $\ubar{\bv{a}} =[a_1;a_2] \in  \funSpace{C}^1([0,T];\funSpace{C}_{pw}^1(\mathcal{E}) \times \funSpace{C}_{pw}^1(\mathcal{E}))$ satisfy
\begin{subequations}\label{eq:System_a}
\begin{align} 
	\partial_t	
	\begin{matrixbj}
		a_1(t) \\
		\derst{2} \LegDens(\ubar{\bv{a}}(t))
	\end{matrixbj}	
	 =
	\begin{bmatrix}
		 & -\partial_x \\
		 -\partial_x & -r(\ubar{\bv{a}}(t))
	\end{bmatrix}
		\begin{matrixbj}
		-\derst{1} \LegDens(\ubar{\bv{a}}(t)) \\
		a_2(t)
	\end{matrixbj}	, \qquad  \,x \in \blsPipe, \quad t \in [0,T]  \label{eq:S_a} \\
		\sum_{\Onepipe \in \mathcal{E}(\blsVertex) } n^{\Onepipe}[\blsVertex]a_{2|\Onepipe}(t)[\blsVertex]= 0  , \quad 
	\nabla_{1} g(\ubar{\bv{a}}_{|\Onepipe}(t)[\blsVertex] = \nabla_{1} g(\ubar{\bv{a}}_{|\tilde{\Onepipe}}(t)[\blsVertex], \quad \Onepipe,\tilde{\Onepipe} \in \mathcal{E}(\blsVertex)
\end{align}
(supplemented with closing conditions according to \eqref{bls-eq:abstr-bc}), where $r(\ubar{\bv{a}}(t))= \tilde{r}(\hat{\bv{z}}(\ubar{\bv{a}}(t)))$ is non-negative and $g$ the partial Legendre transformation of $h$. Then $\hat{\bv{z}}(\ubar{\bv{a}})$ fulfills \eqref{bls-eq:abstr-a}-\eqref{bls-eq:abstr-coup}, and the Hamiltonian $\HamPDE$ of the system is given by
\begin{align}\label{eq:H_a}
	\HamPDE(\ubar{\bv{a}}(t))  = \Lscal \nabla_2 g(\ubar{\bv{a}}(t)) a_2(t)  - \LegDens(\ubar{\bv{a}}(t)), 1\Rscal .
\end{align}
\end{subequations}
\end{crllr}
Note that here and in the following, the solution is interpreted as a function in time with values in a function space. Moreover, by slight abuse of notation, the point-evaluations $\ubv{a}_{|\Onepipe}(t,\blsVertex) \in \mathbb{R}^2$ are written as $\ubv{a}_{|\Onepipe}(t)[\blsVertex]$.

\begin{rmrk}
The Hamiltonian $\HamPDE$ as function of $\ubv a$ is characterized in terms of the Legendre transformation $g$ in \eqref{eq:H_a}, cf.\ Lemma~\ref{lem:legOfleg}. The more direct characterization in the energy variable with respect to the Hamiltonian density $h$ reads $\HamPDE(\ubar{\bv{a}}) = \tilde{\HamPDE}(\hat{\bv{z}}(\ubar{\bv{a}}))= \Lscal h(\hat{\bv{z}}(\ubar{\bv{a}}) ,1 \Rscal$. Note that Assumption~\ref{assum:h-sconvex} implies structural properties, e.g., strict convexity, on $\tilde{\HamPDE}$ rather than on $\HamPDE$.
\end{rmrk}

\begin{rmrk}
Due to the variable transformation the first coupling condition in \eqref{bls-eq:abstr-coup} becomes linear in the variable $a_2$, although the model problem is of a general nonlinear form. This simplifies the handling of the network aspects, since we can use linear ansatz spaces in our approach, similar to \cite{art:egger-mfem-compressEuler}. Note that the incorporation, and beforehand, the identification of consistent coupling conditions for our model problem is a non-trivial issue, cf.\ \cite{art:euler-reigstad, art:hierarchGasMindt2019}.
\end{rmrk}

\subsection{Variational principle} \label{subsec:var-principle}
\begin{thrm} \label{theor:var-principle-twopH-a1z2}
A strong solution $\ubar{\bv{a}}\in \funSpace{C}^1([0,T];\funSpace{C}_{pw}^1(\mathcal{E}) \times \funSpace{C}_{pw}^1(\mathcal{E}))$ of the transformed model problem fulfills for all $b_1 \in \funSpace{L}^2(\mathcal{E})$, $b_2 \in \funSpace{H}_{div}^1(\mathcal{E})$ the variational principle
\begin{align*} 
 \Lscal \partial_t a_1(t), b_1 \Rscal &= - \Lscal \partial_x a_2(t), b_1 \Rscal \\
  \Lscal \partial_t \derst{2}{g}(\ubar{\bv{a}}(t)), b_2 \Rscal &= -\Lscal \derst{1}{g}(\ubar{\bv{a}}(t)), \partial_x b_2 \Rscal + \bcef(t) \cdot \TraceOp b_2 - \Lscal r(\ubar{\bv{a}}(t)) a_2(t), b_2 \Rscal  \qquad & \\
  \bcfl(t)   &=  \TraceOp a_2(t),
\end{align*}
where the functions $\bcef$, $\bcfl: [0,T] \rightarrow \mathbb{R}^{p}$  with one entry for each boundary node $\blsVertex_i \in \mathcal{N}_\partial$, $p=| \mathcal{N}_\partial |$, satisfy
\begin{align*}
	e_i=   -\derst{1}{g}(\ubar{\bv{a}}_{|\Onepipe}[\blsVertex_i]), \quad \text{} \quad 
	f_i =  n^{\Onepipe}[\blsVertex_i]	 a_{2|\Onepipe}[\blsVertex_i] ,\quad \text{for } \Onepipe \in \mathcal{E}(\blsVertex_i), \, i =1,\ldots, p.
\end{align*}
\end{thrm}
\begin{proof}
By testing the second equation of \eqref{eq:S_a} with $b_2 \in \funSpace{H}^1_{div}(\mathcal{E})$, integrating it over one edge $\Onepipe = (\blsVertex, \tilde{\blsVertex}) \in \mathcal{E}$ and using integration by parts, we obtain
\begin{align*}
		 \Lscal \partial_t \nabla_2 g(\ubv{a}(t), b_2 \Rscal_\Onepipe &= -\Lscal \derst{1}{g}(\ubv{a}(t)), \partial_x b_2 \Rscal_\Onepipe - \Lscal r(\ubv{a}(t)) a_2(t), b_2 \Rscal_\Onepipe  \\
		 & \quad + \left[ - \derst{1}{g}(\ubv{a}_{|\Onepipe}(t)[x]) (-b_{2|\Onepipe}[x]) \right]_{x={\blsVertex}}^{\tilde{\blsVertex}}.
\end{align*}
Here, the subscript $._\Onepipe$ indicates the restriction onto the edge $\Onepipe$. Repeating the calculation for all edges and then summing up all equations, the interface terms at the inner nodes drop out, as $b_2 \in  \funSpace{H}^1_{div}(\mathcal{E})$. This gives the second equation of the variational principle. The other equations follow similarly.
\end{proof}
We use the variational principle to show local mass conservation and the energy dissipation equality. 
\begin{thrm}[Structural properties]\label{theor:var-principle-energy-diss}
Let $\ubar{\bv{a}}\in \funSpace{C}^1([0,T];\funSpace{C}_{pw}^1(\mathcal{E}) \times \funSpace{C}_{pw}^1(\mathcal{E}))$ fulfill the variational principle of Theorem~\ref{theor:var-principle-twopH-a1z2} for some $\bcef \in  \funSpace{C}([0,T];\mathbb{R}^{\nbn})$ and $\bcfl \in  \funSpace{C}^1([0,T];\mathbb{R}^{\nbn})$. Then it holds for $[w_a,w_b] \subset \Onepipe$, $\Onepipe \in \mathcal{E}$ and the point evaluations $a_{2|\Onepipe}[w_a]$, $a_{2|\Onepipe}[w_b]$ of $a_2$ that
\begin{align*}
	\frac{d}{dt} \int_{[w_a,w_b]} a_1(t) dx &= a_{2|\Onepipe}(t)[w_a] - a_{2|\Onepipe}(t)[w_b] , \qquad  \frac{d}{dt} \MassPDE(\ubar{\bv{a}}(t)) =  \sum_{i:\,{\blsVertex_i \in \mathcal{N}_\partial}} f_i(t),  \\ 
	\frac{d}{dt} \HamPDE(\ubv{a}(t))  &= \bcef(t) \cdot \bcfl(t) - \left\Lscal r(\ubar{\bv{a}}(t)) a_2(t), a_2(t) \right\Rscal \leq \bcef(t) \cdot \bcfl(t) . 
\end{align*}
\end{thrm}
\begin{proof}
The first relation, the local mass conservation, follows from testing the first equation of the variational principle (Theorem~\ref{theor:var-principle-twopH-a1z2}) with $b_1= \chi _{[w_a,w_b]}$, the indicator function of the domain $[w_a,w_b]$. That is,
\begin{align*}
	\frac{d}{dt} \int_{[w_a,w_b]} a_1 dx &=  \Lscal \partial_t a_1 , \chi _{[w_a,w_b]} \Rscal = - \Lscal \partial_x a_2 , \chi _{[w_a,w_b]} \Rscal = a_{2|\Onepipe}[w_a] - a_{2|\Onepipe}[w_b].
\end{align*}
The second equation, the global mass conservation, results from summing up the individual masses on the edges. 

As a preliminary step for the energy dissipation equality, we prove the existence of $\xi(\ubar{\bv{a}})$ for solutions $\ubar{\bv{a}}$ of the variational principle with $\xi(\ubar{\bv{a}}(t)) \in \funSpace{L}^2(\mathcal{E})$ such that
\begin{align} \label{bls:eq-proof1}
\Lscal -\derst{1}g(\ubar{\bv{a}}), b_1 \Rscal = \Lscal \xi(\ubar{\bv{a}}), b_1 \Rscal \hspace{1cm} \text{for  } b_1 \in \funSpace{L}^2(\mathcal{E}).
\end{align}
As $\ubar{\bv{a}}$ fulfills the variational principle, $\Lscal -\derst{1}g(\ubar{\bv{a}}), \partial_x b_2 \Rscal$ has to be bounded for all $b_2 \in \funSpace{H}_{div}^1(\mathcal{E})$ and $t\in [0,T]$ due to the second equation of Theorem~\ref{theor:var-principle-twopH-a1z2}. The compatibility $\funSpace{L}^2(\mathcal{E}) \subset \left \{\xi: \text{It exists } \zeta \in \funSpace{H}_{div}^1(\mathcal{E}) \text{ with } \partial_x \zeta = \xi \right \}$ follows  from the definition of the broken derivative $\partial_x$. This shows, in turn, , which
that $b_1 \mapsto \Lscal \derst{1}h(\hat{\bv{z}}(\ubar{\bv{a}})), b_1 \Rscal$ is an element of the dual space of $\funSpace{L}^2(\mathcal{E})$. Thus, the existence of $\xi(\ubar{\bv{a}})$ fulfilling \eqref{bls:eq-proof1}  follows from the Riesz representation theorem, as $\funSpace{L}^2(\mathcal{E})$ with $\Lscal \cdot,\cdot \Rscal$ is a Hilbert space.

Let us now turn to the proof of the energy dissipation. A formal application of the chain rule leads to
\begin{align*}
 	\frac{d}{dt} \HamPDE(\ubar{\bv{a}}) &= \frac{d}{dt} \Lscal \nabla_2 g(\ubar{\bv{a}}) a_2  - \LegDens(\ubar{\bv{a}}), 1\Rscal  =
 	\int_\Omega \begin{bmatrix}
 	 -\derst{1}g(\ubar{\bv{a}}) \\
 	 a_2
\end{bmatrix} \cdot 
 	\begin{bmatrix}
 	\partial_t a_1\\
 	 \partial_t \derst{2}g(\ubar{\bv{a}}) 
\end{bmatrix} dx \\
&=
 	 \Lscal \partial_t a_1, -\derst{1}g(\ubar{\bv{a}}) \Rscal  + \Lscal \partial_t \derst{2}g(\ubar{\bv{a}}), a_2 \Rscal.
\end{align*}
By \eqref{bls:eq-proof1} we have $\Lscal \partial_t a_1, -\derst{1}g(\ubar{\bv{a}}) \Rscal= \Lscal \partial_t a_1, \xi(\ubar{\bv{a}}) \Rscal$ and $\Lscal \xi(\ubv{a}), \partial_x b_2 \Rscal_ = \Lscal-\derst{1}g(\ubv{a}) , \partial_x b_2 \Rscal$. Using the variational principle with $\ubar{\bv{b}}=[\xi(\ubar{\bv{a}});a_2]$ yields
\begin{align*}
	&\Lscal \partial_t a_1, -\derst{1}g(\ubar{\bv{a}}) \Rscal  + \Lscal \partial_t \derst{2}g(\ubar{\bv{a}}), a_2 \Rscal \\&=  -\Lscal \partial_x a_2, \xi(\ubar{\bv{a}}) \Rscal - \Lscal \derst{1}g(\ubar{\bv{a}}), \partial_x a_2 \Rscal + \bcef \cdot \TraceOp a_2 
	- \Lscal r(\ubar{\bv{a}}) a_2, a_2 \Rscal 
		 = \bcef \cdot \bcfl - \Lscal r(\ubar{\bv{a}}) , a_2^2 \Rscal \leq \bcef \cdot \bcfl
\end{align*}
due to the non-negativity of $r$, which finishes the proof.
\end{proof}

\begin{rmrk}\label{rem:infin-gradH}
The proof of the energy dissipation relies on the fact that for any function $\ubv{a}$ fulfilling the variational principle (Theorem~\ref{theor:var-principle-twopH-a1z2}), $[- \nabla_1 g(\ubv{a});a_2]$ can be identified with the gradient of the Hamiltonian $\tilde{\HamPDE}(\ubv{z})_{|\ubv{z} = \hat{\bv{z}}(\bv a)}= \Lscal h(\ubv{z})_{|\ubv{z} = \hat{\bv{z}}(\bv a)} ,1 \Rscal$ with respect to the energy variable $\ubv{z}$. In the infinite-dimensional case, gradients of functionals, i.e., identifications of their derivative with an appropriately smooth function, do not exist for all choices of function spaces, cf.\ \cite{book:Zeidler1984NonlinearFA, book:BrezzisFuncAna}. The gradient structure is a non-trivial property that our special variational principle has. It plays a fundamental role for the structure-preserving properties of our proposed approximation approach.
\end{rmrk}

For each boundary node $\blsVertex_i \in \mathcal{N}_\partial$, we assume a boundary data $u_{\nu_i}: [0,T] \rightarrow \mathbb{R}$ to be given and pose a boundary condition of the form $k_i(e_i,f_i,u_{\blsVertex_i})=0$. This equation only depends on the boundary effort $e_i$ and the boundary flow $f_i$ (Theorem~\ref{theor:var-principle-twopH-a1z2}). For a more concise notation, we collect all boundary conditions in a vector-valued function $\bv{k}:\mathbb{R}^{p}\times\mathbb{R}^{p} \times \mathbb{R}^{p} \to \mathbb{R}^{p}$, which yields 
\begin{align*}
\bv{k}(\bcef,\bcfl,\bv{u})=\bv{0}, \qquad \text{ with }\quad \bv u = [u_{\nu_1}; u_{\nu_2}; \ldots; u_{\nu_p}].
\end{align*}  
Note that the type of boundary condition affects the necessary regularity assumptions on the boundary data. To obtain for example a continuously differentiable solution $\ubv a$, it can be seen that a boundary condition of the type $f_i= u_{{\nu_i}}$ requires at least continuously differentiable data $u_{\nu_i}$, whereas for a boundary condition of the type $e_i= u_{\nu_i}$ continuity can be sufficient.

In addition, the initial condition of the model problem are stated as  $\ubv{a}(0)=\ubv{a}_0$.


\section{Spatial approximation approach} \label{sec:approx-framew}

The introduced variational principle is the basis for our structure-preserving spatial approximation approach which covers Galerkin projection and additional complexity reduction of the nonlinear terms. Local mass conservation and an energy bound for the approximations as well as a gradient structure related to the Hamiltonian of the system are kept under mild assumptions. 

\subsection{Galerkin approximation} \label{sec-p1b-galapp}
The Galerkin approximation consists of a projection onto a suitable finite-dimensional ansatz space. It is applicable to a finite element discretization, but also in a more general setting, such as projection-based model reduction. Similar to the previous works \cite{art:kugler-dwe-net,art:lilsailer-dwemor,inproc:bls-hyp2018} for linear and semilinear variants of our model problem, we impose few compatibility conditions as the only restriction.

\begin{assumption}[Compatibility of spaces] \label{assum:compatV1V2}
Let $\funSpace{V}= \funSpace{V}_1 \times \funSpace{V}_2 \subset \funSpace{L}^2(\mathcal{E})\times \funSpace{H}_{div}^1(\mathcal{E})$ be a finite-dimensional subspace that fulfills the compatibility conditions
\begin{enumerate}[({A}1)]
\item $\funSpace{V}_{1} =   \partial_x \funSpace{V}_{2}, \qquad  \text{ with }  \partial_x \funSpace{V}_{2}= \left \{\xi: \text{It exists } \zeta \in \funSpace{V}_{2} \text{ with } \partial_x \zeta = \xi \right \} $.
\item $\{b_2 \in \funSpace{H}^{1}_{div}(\mathcal{E}): \, \partial_x b_2 = 0 \} \subset \funSpace{V}_2$.
\end{enumerate}
\end{assumption}
\begin{rmrk}
In this paper only condition (A1) is explicitly used. Nonetheless, (A2) should be included as the analysis of the related steady state problem reveals, see \cite{art:kugler-dwe-net, inproc:bls-hyp2018} for the linear and semilinear case.
\end{rmrk}
The proposed Galerkin approximation for our model problem reads as follows.
\begin{system}\label{sys:twopH-gallapp}
Assumption~\ref{assum:compatV1V2} is supposed to hold. Given initial and boundary data, $\ubv a_0 \in \funSpace{V}$ and $u_\nu : [0,T]\rightarrow \mathbb{R}$ for $\nu \in \mathcal{N}_\partial$, find $\ubar{\bv{a}}\in \funSpace{C}^1([0,T];\funSpace{V}_1 \times \funSpace{V}_2) $, \mbox{$\bcfl \in  \funSpace{C}^1([0,T];\mathbb{R}^{\nbn})$} and $\bcef \in  \funSpace{C}([0,T];\mathbb{R}^{\nbn})$ that solve
\begin{subequations} \label{eq:system_all}
\begin{align} 
 \Lscal \partial_t a_1(t), b_1 \Rscal &= - \Lscal \partial_x a_2(t), b_1 \Rscal  \label{eq:system_1} \\ 
  \Lscal \partial_t \nabla_2 g(\ubv{a}(t)), b_2 \Rscal &= -\Lscal \nabla_1 g(\ubv{a}(t)), \partial_x b_2 \Rscal + \bcef(t) \cdot \TraceOp b_2 - \Lscal r(\ubar{\bv{a}}(t)) a_2(t), b_2 \Rscal  \label{eq:system_2} \\
  \bcfl(t)   &=  \TraceOp a_2(t)
\end{align}
\end{subequations}
for all $b_1 \in \funSpace{V}_1$, $b_2 \in \funSpace{V}_2$ and fulfill $\ubv a(0)=\ubv a_0$ and $\bv{k}(\bcef,\bcfl,\bv{u})= \bv{0}$.
\end{system}
The compatibility conditions on our ansatz spaces (Assumption~\ref{assum:compatV1V2}) guarantee several structural properties for our Galerkin approximation. 

\begin{lmm}\label{lem:space-comp-Gradb}
Let $\ubv a$ be a solution of System~\ref{sys:twopH-gallapp}, then the derivative of the functional $\LegPDE:\funSpace{V} \to \mathbb{R}$, $\ubv{a} \mapsto \Lscal g(\ubv{a}),1\Rscal$ associated to the partial Legendre transformation can be identified with a function in $\funSpace{V}$ in the sense that there exists $[\xi(\ubv{a}(t));\tilde{\xi}(\ubv{a}(t))]\in \funSpace{V}$ for $t\in [0,T]$, such that for $b_1 \in \funSpace{V}_1, b_2 \in \funSpace{V}_2$ it holds
\begin{align*}
	\Lscal \xi(\ubv{a}(t)), b_1 \Rscal = \Lscal -\nabla_1 g(\ubv{a}(t)), b_1 \Rscal, \qquad 	\Lscal \tilde{\xi}(\ubv{a}(t)), b_2 \Rscal = \Lscal \partial_t \nabla_2 g(\ubv{a}(t)), b_2 \Rscal.
\end{align*}
\end{lmm}
\begin{proof}
If $\ubv a$ fulfills \eqref{eq:system_all}, the expression $\Lscal -\nabla_1 g(\ubv{a}(t), \tilde{b} \Rscal$ is well-defined for any $\tilde{b}\in \partial_x \funSpace{V}_2$. Thus, it is well-defined for any $\tilde{b} \in \funSpace{V}_1$ due to the compatibility condition $\funSpace{V}_1 \subset \partial_x\funSpace{V}_2$. By that, $b_1 \mapsto \Lscal -\derst{1}g(\ubar{\bv{a}}), b_1 \Rscal$ is an element of the dual space of $\funSpace{V}_1$, which can be identified with a $\xi(\ubv{a})\in \funSpace{V}_1$ by the Riesz representation theorem. The claim on the existence of $\tilde{\xi}(\ubv{a})$ follows similarly by using that $\funSpace{V}_2$ with $\Lscal \cdot, \cdot \Rscal$ is a Hilbert space.
\end{proof}
Thus, our approximation inherits the gradient structure of the model problem. Here, in the finite-dimensional setting, the result is slightly stronger than in the infinite-dimensional setting, where $\funSpace{H}_{div}^1(\mathcal{E})$ with $\Lscal \cdot, \cdot \Rscal$ is not a Hilbert space, cf.\ Remark~\ref{rem:infin-gradH}.

\begin{thrm}\label{theor:energy-bound-galapp}
Let $\ubar{\bv{a}}$ be a solution of System~\ref{sys:twopH-gallapp}. Then the following two properties, related to local mass conservation and energy dissipation, are fulfilled for $t\in[0,T]$,
 \begin{align*}
 \partial_t a_1(t) &= -\partial_x a_2(t),\\
 \frac{d}{dt} \HamPDE(\ubar{\bv{a}}(t)) &= \bcef(t) \cdot \bcfl(t) - \left \Lscal r(\ubar{\bv{a}}(t)) a_2(t) , a_2(t) \right \Rscal  \leq \bcef(t) \cdot \bcfl(t).
 \end{align*}
\end{thrm}
\begin{proof}
The first relation implies that mass conservation holds in a pointwise sense, i.e., locally. To show it, note that by construction of the ansatz spaces $\partial_t a_1({{t}})$, $\partial_x a_2({{t}}) \in \funSpace{V}_1$ for ${{t}} \in [0,T]$ holds. Thus, the variational principle with $b_1 = \partial_t a_1(t)+\partial_x a_2(t)$ gives
\begin{align*}
	 \Lscal \partial_t a_1(t) , \partial_t a_1(t) + \partial_x a_2({{t}}) \Rscal  = \Lscal - \partial_x a_2({{t}}) , \partial_t a_1({{t}}) + \partial_x a_2({{t}}) \Rscal .
\end{align*}
Subtracting the right-hand side reveals that the $\funSpace{L}^2$-norm of $\partial_t a_1(t) + \partial_x a_2(t)$ is zero, i.e., $\partial_t a_1(t)= -\partial_x a_2(t) $.

The proof of the energy dissipation equality is the same as that of Theorem~\ref{theor:var-principle-energy-diss}, except for replacing the spaces $\funSpace{L}^2(\mathcal{E})$ and $\funSpace{H}_{div}^1(\mathcal{E})$ by $\funSpace{V}_1$ and $\funSpace{V}_2$ and using Lemma~\ref{lem:space-comp-Gradb}.
\end{proof}
There are a few possibilities for constructing a suitable finite-dimensional approximation space $\funSpace{V}$ in the sense of Assumption~\ref{assum:compatV1V2}, cf.\ \cite{inproc:bls-hyp2018,art:lilsailer-dwemor}. We generate one by means of mixed finite elements. We divide the edges $\Onepipe \in \mathcal{E}$ into sub-parts $T_{\Onepipe,k}$, $k=1,\ldots,J_{\Onepipe}$ and consider the polynomial spaces
\begin{align*}
\funSpace{Q}_q(T_{\Onepipe,k};\mathbb{R}) &= \left\{ \phi: T_{\Onepipe,k} \rightarrow \mathbb{R}: \, \phi(x) = \sum_{j=0}^{q} x^j \xi_j, \text{ for } x\in T_{\Onepipe,k}, \text{ with } \xi_j \in \mathbb{R}  \right\}.
\end{align*}
The sub-parts $T_{\Onepipe,k}$ induce a partitioning $T_{\mathcal{E}}$ of the full network, on which we define the piecewise polynomial spaces
\begin{align*}
	\funSpace{Q}_q(T_{\mathcal{E}}) &= \left\{ \phi: \Omega_{} \rightarrow \mathbb{R}: \, \phi_{|T_{\Onepipe,k}} \in \funSpace{Q}_q(T_{\Onepipe,k}), \text{ for } k=1,\ldots,J_{\Onepipe}, \, \Onepipe \in \mathcal{E}\right\} , \\
	\funSpace{P}_q(T_{\mathcal{E}}) &= \funSpace{Q}_q(T_{\mathcal{E}}) \cap \funSpace{H}^1_{div}(\mathcal{E}).
\end{align*}
Particularly, $\funSpace{P}_q(T_{\mathcal{E}})$ inherits the coupling conditions encoded in $\funSpace{H}^1_{div}(\mathcal{E})$, which predefines one degree of freedom per inner node $\nu \in \mathcal{N}_0$.
The choice $\funSpace{V}_1 = \funSpace{Q}_{q}(T_{\mathcal{E}})$ and $\funSpace{V}_2 = \funSpace{P}_{q+1}(T_{\mathcal{E}})$ for any $q\geq 0$ yields a pair of compatible spaces.

\subsection{Complexity reduction} \label{sec-p1b-comred}
System~\ref{sys:twopH-gallapp} features nonlinearities in the balance equation \eqref{eq:system_2} due to the partial Legendre transformation~$g$ and the model function~$r$. Depending on the model problem and the choice of the Galerkin ansatz space $\funSpace{V}$, an additional complexity reduction can become necessary for an efficient numerical realization. We propose a quadrature-type approximation that preserves the structural properties. In particular, we allow for inexact integration, which is, e.g., crucial in the design of online-efficient reduced order models, cf.\ \cite{art:comred-ecsw,art:efficient-integration-cubature}, but also may be convenient for high-order finite elements.

Given appropriate quadrature points $x_i \in \Omega$ and quadrature weights $w_i$ for $i\in I$ with $I$ describing an index set, we introduce the complexity-reduced approximations
$\Lscal \cdot, \cdot\Rscal_c$ and $|| \cdot ||_c$ of the $\funSpace{L}^2$-scalar product $\Lscal \cdot, \cdot \Rscal$ and $\funSpace{L}^2$-norm $||\cdot ||$  as
\begin{align*}
\Lscal b, \bar{b}\Rscal_c = \sum_{i\in I} w_i b[x_i] \bar{b}[x_i], \qquad ||b||_c = \sqrt{\Lscal b,b \Rscal_c},
\end{align*}
where  $b[x_i]$ refers to point evaluations of $b$.
The following assumption ensures that the quadrature rule is well-defined and stable and that $||\cdot ||_c$ is a norm on $\funSpace{V}_1$ and $\funSpace{V}_2$.

\begin{assumption} \label{assum:quadrat-FIRST}
For $i\in I$, the evaluations $b[x_i]$ are well-defined for  $b\in \funSpace{V}_1 \cup \funSpace{V}_2$, and the quadrature weights are positive, $w_i>0$. Moreover, there exists a constant $\tilde{C}>0$ such that
\begin{align*}
	\frac{1}{\tilde{C}} ||b||_c \leq ||b|| \leq \tilde{C} ||b||_c , \hspace{1.5cm} \text{for } b \in \funSpace{V}_1 \cup \funSpace{V}_2.
\end{align*}
\end{assumption}
Our complexity reduction of System~\ref{sys:twopH-gallapp} modifies only the nonlinear balance equation~\eqref{eq:system_2}, it is replaced by
\begin{align*} 
  \Lscal \partial_t \nabla_2 g(\ubv{a}(t)), b_2 \Rscal_{c} &= \Lscal - \nabla_1 g(\ubv{a}(t)), \partial_x b_2 \Rscal_{c} + \bcef(t) \cdot \TraceOp b_2 - \Lscal r(\ubar{\bv{a}}(t)) a_2(t), b_2 \Rscal_c 
\end{align*}
for $b_2 \in \funSpace{V}_2$. Further, we consider the adjusted functional $\LegPDE_c: \ubv{a} \mapsto \Lscal g(\ubv{a}),1\Rscal_c$ associated to the partial Legendre transformation and the complexity-reduced Hamiltonian $\HamPDE_c(\ubar{\bv{a}})  = \Lscal \nabla_2 g(\ubar{\bv{a}}) a_2  - \LegDens(\ubar{\bv{a}}), 1\Rscal_c$. The gradient structure is kept, since the derivative of $\LegPDE_c$ can be identified with $[\xi(\ubv{a}(t));\tilde{\xi}(\ubv{a}(t))]\in \funSpace{V}$ for $t\in [0,T]$ in the sense that
\begin{align*}
	\Lscal \xi(\ubv{a}(t)), b_1 \Rscal = \Lscal -\nabla_1 g(\ubv{a}(t)), b_1 \Rscal_c, \qquad 	\Lscal \tilde{\xi}(\ubv{a}(t)), b_2 \Rscal = \Lscal \partial_t \nabla_2 g(\ubv{a}(t)), b_2 \Rscal_c
\end{align*}
holds for $b_1 \in \funSpace{V}_1$, $b_2 \in \funSpace{V}_2$.
Note that $b_1 \mapsto \Lscal -\derst{1}g(\ubar{\bv{a}}), b_1 \Rscal_c$ and $b_2 \mapsto \Lscal \partial_t\derst{2}g(\ubar{\bv{a}}), b_2 \Rscal_c$ are bounded linear functionals, i.e., elements of the dual space of $\funSpace{V}_1$ and $\funSpace{V}_2$, respectively. This holds for any choice of inner product, in particular also for the $\funSpace{L}^2$-inner product used above. The energy dissipation equality becomes
\begin{align*}
	 \frac{d}{dt} \HamPDE_c(\ubv{a}(t)) = \bcef(t) \cdot \bcfl(t) - \left \Lscal r(\ubar{\bv{a}}(t))a_2(t) , a_2(t) \right \Rscal_c  \leq \bcef(t) \cdot \bcfl(t).
\end{align*}
The proof of the structural properties follows straightforward the arguments of Lem{\-}ma~\ref{lem:space-comp-Gradb} and Theorem~\ref{theor:energy-bound-galapp} using Assumption~\ref{assum:quadrat-FIRST} on the equivalence of the norms. The local mass conservation is not affected by the complexity reduction at all.

\begin{rmrk}
Methods from literature on model reduction for port-Hamiltonian systems typically employ interpolation-type complexity reduction by DEIM. This yields either approximations, which are only port-Hamiltonian up to an approximation error \cite{art:morParamHamHesthaven, art:moramHamDissHesthaven}, or the Hamiltonian structure is enforced by an additional symmetrization step \cite{art-mor-structure-preserving-nonlinear-pH}, which sacrifices fidelity. For comparisons with our structure-preserving approach we refer to our follow-up paper \cite{art:bls21-snapshotbased}.
\end{rmrk}


\section{Structured representations} \label{sec:Part1b:coordrepres}

In this section we make the transition from the function space setting to the algebraic setting and derive structured coordinate representations for the approximations. The coordinate representations can be transformed into the standard port-Hamiltonian form \eqref{eq:pH-first-prototype}. They are natural discrete counterparts to the formulations of the continuous model \eqref{eq:System_a}, and the underlying Hamiltonian does not degenerate.

Let $\{b^1_i,\ldots, b^{n_i}_i\}$ be the basis of the approximation space $ \funSpace{V}_i$, $\mathrm{dim}( \funSpace{V}_i)=n_i$, for $i=1,2$ with $n = n_1+n_2$.
The bijective mapping between the coordinate representation $\bv{a}= [\bv{a}_1; \bv{a}_2] \in \mathbb{R}^n$, $\bv{a}_i= [\scalco{a}_i^1; \ldots ; \scalco{a}_i^{n_i}] \in \mathbb{R}^{n_i}$ and the function $\ubar{\bv{a}} =[a_1;a_2] \in \funSpace{V}$ is given by
\begin{align}
	\Psi : \mathbb{R}^n \rightarrow \funSpace{V}, \qquad  \Psi(\bv{a}) =
	\begin{bmatrix}
	 \sum_{j=1}^{n_1} b_1^j \scalco{a}_1^j	\\
 	 \sum_{j=1}^{n_2} b^j_2 \scalco{a}_2^j	
	\end{bmatrix} = \begin{bmatrix}
	 a_1 \\
 	 a_2	
	\end{bmatrix}=\ubar{\bv{a}}. \label{eq:coord-imbedding}
\end{align} 
\subsection{Representation of Galerkin approximation} \label{subsec:struct-gal}

For the functional $\LegPDE: \funSpace{V} \rightarrow \mathbb{R}$, $\LegPDE(\ubar{\bv{a}}) = \Lscal \LegDens(\ubar{\bv{a}}), 1 \Rscal$ associated to the partial Legendre transformation $g$ of the Hamiltonian density $h$, we introduce the coordinate representation $\LegCo: \mathbb{R}^n \rightarrow \mathbb{R}$ and a weighted gradient $\nabla^\bt{M} \LegCo: \mathbb{R}^n \rightarrow \mathbb{R}^n$ of block-structure  $\nabla^\bt{M} \LegCo = [\nabla_1^\bt{M} \LegCo;\nabla_2^\bt{M} \LegCo ] $ by
\begin{align*}
	 \LegCo( \bv{a}) &= \LegPDE \left(\Psi(\bv{a}  \right)),\\
		\nabla^{\bt{M}}_i \LegCo(\bv{a}) &=  \bt{M}_i^{-1} \nabla_{i} \LegCo(\bv{a}) = \bt{M}_i^{-1}  \left[ {\partial_{\scalco{a}_i^1}} \LegCo(\bv{a}) ; \ldots ;  {\partial_{\scalco{a}_i^{n_i}}} \LegCo(\bv{a}) \right]\in 
		\mathbb{R}^{n_i},\\
		&\quad \,\, \bt M =  \begin{bmatrix}
	\bt M_1 & \\
			& \bt M_2
	\end{bmatrix}, \quad \bt M_i \in \mathbb{R}^{n_i \times n_i}, \, i=1,2,
\end{align*}
where ${\partial_{\scalco{a}_i^j}} \LegCo(\bv{a}) = \left \Lscal  \derst{i} \LegDens(\ubar{\bv{a}})_{|\ubar{\bv{a}}= \Psi(\bv{a})}, b_i^j \right \Rscal$ holds by the chain rule. The matrix $\bt M$ is assumed to be symmetric positive definite. The weighted gradient and its sub-blocks relate to gradients with respect to the inner products induced by $\bt M$ and $\bt M_i$, respectively, cf. \cite{book:matrixManifoldsAbsil08,book:diffgeo-mcInerney}. By choosing $\bt M$ as the mass matrix we can express the gradient structure observed in Lemma~\ref{lem:space-comp-Gradb} in the algebraic setting.
\begin{crllr} \label{corol:twopH-gall-coordRepres}
Let the following system matrices be given
\begin{align*}
	\bt{M}_i &= \left[ \Lscal b_i^k,b_i^j \Rscal \right]_{j,k=1,\ldots, n_i}, \quad i = 1,2, \hspace{1.5cm}
\bt{M} =
 \begin{matrixbj}
	\bt{M}_1 & \\
		& \bt{M}_2
\end{matrixbj}, \\
\bt{J} &= \left[ \Lscal - \partial_x b_2^k,b_1^j \Rscal  \right]_{j=1,\ldots, n_1,\, k = 1,\ldots, n_2} , \hspace{0.9cm} 
 \bt{R}(\bv{a}) = \left[ \Lscal r(\Psi(\bv{a})) b_2^k,b_2^j \Rscal \right]_{j,k=1,\ldots, n_2},\\
\bt{K}_2 &= \left[\TraceOp b_2^1 | \ldots  | \TraceOp b_2^{n_2} \right]^T, 
\hspace{3.43cm}
\bt{K}= 
\begin{matrixbj}
	 \bt{0}_{n_1,\nbn} \\
	\bt{K}_2
\end{matrixbj},
\end{align*}
then System~\ref{sys:twopH-gallapp} can be equivalently stated as: 

Given initial and boundary data, $\bv a_0 \in \mathbb{R}^{n}$ and $u_\nu: [0,T]\rightarrow \mathbb{R}$ for $\nu \in \mathcal{N}_\partial$, find $\bv{a}\in \funSpace{C}^1([0,T];\mathbb{R}^n) $, $\bcfl \in  \funSpace{C}^1([0,T];\mathbb{R}^{\nbn})$ and $\bcef \in  \funSpace{C}([0,T];\mathbb{R}^{\nbn})$ that solve
\begin{align} 
\begin{aligned} \label{eq:S_a_discrete}
\bt{M}\, 
\frac{d}{dt}\begin{bmatrix}
 \bv{a}_1(t) \\
	\nabla^{\bt{M}}_2 \LegCo(\bv{a}(t))
\end{bmatrix}
&= \begin{bmatrix}
	 & \bt{J} \\
	 -\bt{J}^T & -\bt{R}(\bv{a}(t))
\end{bmatrix} 
\begin{bmatrix}
	-\nabla^{\bt{M}}_1 \LegCo(\bv{a}(t)) \\
	\bv{a}_2(t)
\end{bmatrix}
+ \bt{K}
 \bcef(t)   \\ 
  \bcfl(t)   &= 
  \bt{K}^T \bv{a}(t)
\end{aligned}  
\end{align}
and fulfill $\bv{a}(0) = \bv{a}_0$ and $\bv{k}(\bv{e}(t),\bv{f}(t),\bv{u}(t)) = \bv{0}$ for $t\in[0,T]$. Moreover, the Hamiltonian of the system reads $\HamCo(\bv a) =  \nabla_2 \LegCo(\bv{a}) \cdot \bv a_2 - \LegCo(\bv{a})$.
\end{crllr}
Using Lemma~\ref{lem:space-comp-Gradb} and the upper definitions, the validity of Corollary~\ref{corol:twopH-gall-coordRepres} can be shown.
Given the weighted $\funSpace{L}^2$-inner product $\Lscal \cdot , \cdot \Rscal_M$, for which the basis vectors of~$\funSpace{V}_1$, $\{b^1_1,\ldots, b^{n_1}_1\}$, are orthonormal, it follows for $ b_2^k \in \funSpace{V}_2$ by Assumption~\ref{assum:compatV1V2} that 
\begin{align*}
	\Lscal \nabla_1 g(\ubv{a}), \partial_x b_2^k \Rscal = \left \Lscal \sum_{j=1}^{n_1} \Lscal \nabla_1 g(\ubv{a}), b_1^j \Rscal_M \, b_1^j , \partial_x b_2^k \right \Rscal= 
	\sum_{j=1}^{n_1}  \Lscal  b_1^j , \partial_x b_2^k \Rscal \,{\Lscal \nabla_1 g(\ubv{a}), b_1^j \Rscal_M} ,
\end{align*}
where $\Lscal \nabla_1 g(\ubv{a}), b_1^j \Rscal_M=[\nabla^{\bt{M}}_1 \LegCo(\bv{a})]_j$ with mass matrix $\bt M$. The expression relates to the term $\bt{J}^T  \nabla^{\bt{M}}_1 \LegCo(\bv{a}(t))$ in \eqref{eq:S_a_discrete}. Analogously, the structured algebraic representation of the integral expression with $\nabla_2 g(\ubv{a})$ can be derived. Note that the block matrices 
\begin{align*}
\begin{bmatrix}
 \bt 0 &  \bt{J} \\
  -\bt{J}^T &    \bt 0
\end{bmatrix}
 \quad \text{and} \quad
 \begin{bmatrix}
 \bt 0 &  \bt 0 \\
  \bt 0 &    \bt R(\bv a)
\end{bmatrix},
\end{align*}
which can be found in \eqref{eq:S_a_discrete}, are skew-symmetric and symmetric positive semi-definite, respectively,
due to the non-negativity of $r$.

\begin{thrm} \label{th:comred-LegTrafo-inj-2}
Let $\mathbb{A}_n \subset \mathbb{R}^n$ be an open convex set on the domain of $\LegCo$. Then, the mapping
$\hat{\mathfrak{ {z} }}: \mathbb{A}_n\rightarrow W$ with 
\begin{align*}
\hat{\mathfrak{ {z} }}({\bv{a}}) = 	
			[\bv{a}_1;
			\nabla^{\bt{M}}_2  \LegCo(\bv{a})]
		, \qquad \text{for } \bv{a} = 
			[\bv{a}_1;
			\bv{a}_2]\in \mathbb{R}^{n_1+n_2} \text{ and } W=\hat{\mathfrak{ {z} }}(\mathbb{A}_n)\subset \mathbb{R}^n
\end{align*}
is bijective. Let $\bv{a}$ be a solution of \eqref{eq:S_a_discrete}, then $\hat{\mathfrak{ {z} }}(\bv{a})$ solves a system in standard port-Hamiltonian form, cf.~\eqref{eq:pH-first-prototype}, with the Hamiltonian given as the partial Legendre transformation of $G$.
\end{thrm}

\begin{proof}
For any fixed $\bar{\bv{a}}_1$ the function $\Phi: \bv{a}_2 \mapsto \LegCo([\bar{\bv{a}}_1;\bv{a}_2])$
is strictly convex on every convex open set on its domain, since for $\breve{\bv{a}}_2 \neq \tilde{\bv{a}}_2$ (and $\Psi$ as in \eqref{eq:coord-imbedding})
{\small \begin{align*} 
	\left( \nabla_{} \Phi(\breve{\bv{a}}_2) - \nabla_{} \Phi(\tilde{\bv{a}}_2) \right) &\cdot \left( \breve{\bv{a}}_2- \tilde{\bv{a}}_2 \right)  
	\\ 
	&=\sum_{j=1}^{n_2}	
	\left( \left \Lscal  \derst{2} \LegDens(\ubar{\bv{a}})_{|\ubar{\bv{a}}= \Psi([\bar{\bv{a}}_1;\breve{\bv{a}}_2])}-
	                     \derst{2} \LegDens(\ubar{\bv{a}})_{|\ubar{\bv{a}}= \Psi([\bar{\bv{a}}_1;\tilde{\bv{a}}_2])}, b_2^j \right \Rscal   ( \breve{\scalco{a}}_2^j- \tilde{\scalco{a}}_2^j) \right)
	          \\ 
	 & =\left \Lscal \derst{2} \LegDens([\bar{a}_1;\breve{a}_2])-\derst{2} \LegDens([\bar{a}_1;\tilde{a}_2]),\breve{a}_2-\tilde{a}_2
	  \right \Rscal >0,
\end{align*}} 
due to the strict convexity of $g$ with respect to its second component as consequence of Assumption~\ref{assum:h-sconvex}. The bijectivity of $\hat{\mathfrak{ {z} }}$ follows then from Theorem~\ref{theor:partialLeg-GEN}. The port-Hamiltonian system is obtained from  \eqref{eq:S_a_discrete} via this coordinate transformation.
\end{proof}

\subsection{Representation of complexity reduced system}
An equivalence transformation of our approximation to a port-Hamiltonian system in standard form and a counterpart to Theorem~\ref{th:comred-LegTrafo-inj-2} also exist for our complexity reduced system. Because of the strong similarity to the case without complexity reduction, we only highlight the differences. The structured coordinate representation belonging to the complexity reduction equals  \eqref{eq:S_a_discrete} except of the expressions  $\LegCo$ and $\bt{R}$ that are replaced by
\begin{align*}
	 \LegCo_c(\bv{a}) = \LegPDE_c \left(\Psi(\bv{a}  \right))=\Lscal \LegDens(\Psi(\bv{a})),1 \Rscal_c, \hspace{1cm}
	 \bt{R}_c(\bv{a}) = \left[ \Lscal r(\Psi(\bv{a})) b_2^k,b_2^j \Rscal_c \right]_{j,k=1,\ldots, n_2}. 
\end{align*}
The respective Hamiltonian is $\HamCo_c(\bv a) =   \nabla_2 \LegCo_c(\bv{a}) \cdot \bv a_2 - \LegCo_c(\bv{a})$. The coordinate transformation into standard port-Hamiltonian form is given by
$\hat{\mathfrak{ {z} }}_c(\bv a) = [\bv{a}_1; \nabla^{\bt{M}}_2  \LegCo_c(\bv{a})]$ for $\bv{a} = 
[\bv{a}_1; \bv{a}_2]$. Its bijectivity can be concluded analogously as in the proof of Theorem~\ref{th:comred-LegTrafo-inj-2}. Note that for the strict convexity of $\Phi: \bv{a}_2 \mapsto \LegCo_c([\bar{\bv{a}}_1;\bv{a}_2])$ for fixed~$\bar{\bv{a}}_1$, Assumption~\ref{assum:quadrat-FIRST} is needed.
We can argue that for $\breve{\bv{a}}_2 \neq \tilde{\bv{a}}_2$
\begin{align*}
	\left( \nabla_{} \Phi(\breve{\bv{a}}_2) - \nabla_{} \Phi(\tilde{\bv{a}}_2) \right) \cdot \left( \breve{\bv{a}}_2- \tilde{\bv{a}}_2 \right) = 
	\left \Lscal \derst{2} \LegDens([\bar{a}_1;\breve{a}_2])-\derst{2} \LegDens([\bar{a}_1;\tilde{a}_2]),\breve{a}_2-\tilde{a}_2
	  \right \Rscal_c\\
	=\sum_{i\in I}  w_i \, \left(\derst{2} \LegDens([\bar{a}_1;\breve{a}_2])[x_i]- \derst{2} \LegDens([\bar{a}_1;\tilde{a}_2])[x_i] \right) \, \left( \breve{a}_2[x_i] - \tilde{a}_2[x_i] \right) >0
\end{align*}
is valid, because the quadrature weights are positive ($w_i> 0$) and there exists at least one quadrature point $x_j$, $j\in I$ with $(\breve{a}_2[x_j]-\tilde{a}_2[x_j]) \neq 0$. The latter follows from
\begin{align*}
	\sum_{i\in I}  w_i \, (\breve{a}_2[x_i] -\tilde{ a}_2[x_i])^2 = || \breve{a}_2 - \tilde{a}_2 ||_c^2\geq \frac{1}{\tilde{C}^2} ||\breve{a}_2 - \tilde{a}_2 ||^2 > 0,
\end{align*}
using $\breve{a}_2-\tilde{a}_2 \in \funSpace{V}_2$ and the equivalence of $||\cdot ||_c$ and $|| \cdot ||$ on $\funSpace{V}_2$. The strict convexity of $g$ yields 
\begin{align*}
		\left(\derst{2} \LegDens([\bar{a}_1;\breve{a}_2 ])[x_j] - \derst{2} \LegDens([\bar{a}_1;\tilde{a}_2])[x_j] \right) \,\, \left( \breve{a}_2[x_j] - \tilde{a}_2[x_j] \right)   > 0, 
\end{align*}
(i.e., strict positivity at $x_j$) and ensures non-negativity at all other quadrature points $x_i$, $i\in I$.


\section{Time discretization} \label{sec:p1b:time-disc}
The underlying Hamiltonian structure of our space approximation can be exploited for analysis and control purposes but also for the derivation of energy-stable or energy-preserving time discretization schemes. Probably the simplest energy-stable scheme is provided by the implicit Euler-type method discussed in the following. 

We consider an equidistant time grid $t_k=k\Delta_t \in [0,T]$, $k=0,...,K$, $K=T/\Delta_t$ with grid size $\Delta_t$ and indicate the temporal approximations by a respective super-index, e.g., $\bv{a}^k\approx \bv{a}(t_k)$.  
\begin{system}[Implicit Euler-type scheme] \label{sys:twopH-gallapp-ImplicitEuler}
Given initial data $\bv{a}_0 \in \mathbb{R}^n$ and boundary data $u_{\nu_i}^k \in \mathbb{R}$ for $\nu_i \in \mathcal{N}_\partial$, find $\bv{a}^k = [\bv a_1^k;\bv a_2^k] \in \mathbb{R}^n$ and $\bcfl^k = [f_{\nu_1}^k;\ldots;f_{\nu_p}^k]$, $\bcef^k= [e_{\nu_1}^k;\ldots;e_{\nu_p}^k] \in  \mathbb{R}^{\nbn}$ for $k > 0$ by solving
\begin{align*} 
\frac{1}{\Delta_t} \bt{M}\, 
\left(
\begin{bmatrix} 
 \bv{a}_1^k \\
	\nabla^{\bt{M}}_2 \LegCo(\bv{a}^k) 
\end{bmatrix}
- \begin{bmatrix} \bv{a}_1^{k-1} \\
	\nabla^{\bt{M}}_2 \LegCo(\bv{a}^{k-1})
\end{bmatrix}
\right)
&= \begin{bmatrix}
	 & \bt{J} \\
	 -\bt{J}^T & -\bt{R}(\bv{a}^k)
\end{bmatrix} 
\begin{bmatrix}
	-\nabla^{\bt{M}}_1 \LegCo(\bv{a}^k) \\
	\bv{a}_2^k
\end{bmatrix}
+ \bt{K}
 \bcef^k    \\
  \bcfl^k   &= 
  \bt{K}^T \bv{a}^k 
\end{align*}
with closing conditions $\bv a^0 = \bv a_0$ and $\bv{k}( \bcef^{k}, \bcfl^{k},\bv u^{k}) = \bv{0}$.
\end{system}
System~\ref{sys:twopH-gallapp-ImplicitEuler} can be interpreted as an implicit Euler discretization in the energy variable~$\bv{z}$ and thus is of first-order convergence in time.
A time-discrete counterpart to the energy bound from Theorem~\ref{theor:energy-bound-galapp} can be derived. To do so, we make use of the next auxiliary result, which follows from the strict convexity of the Hamiltonian with respect to the energy variable.
\begin{lmm}\label{lem:HamConvexCond}
For the Hamiltonian $\HamCo$ and the functional $\LegCo$ given as in Corollary~\ref{corol:twopH-gall-coordRepres}, it holds for $\bv a \neq \bar{\bv a }$ that
\begin{align*}
	 \HamCo({\bv a})  -\HamCo(\bar{\bv a}) <
	\begin{bmatrix}
	-\nabla_1 \LegCo(\bv{a}) \\
	\bv a_2 
	\end{bmatrix} \cdot \left(
	\begin{bmatrix}
		\bv a_1 \\
		\nabla_2 \LegCo(\bv{a})
	\end{bmatrix} 
	-
		\begin{bmatrix}
		\bar{\bv{a}}_1 \\
		\nabla_2 \LegCo(\bar{\bv{a}})
	\end{bmatrix} 	
	\right) ,\hspace{0.6cm} \bv a = 
	\begin{bmatrix}
		\bv a_1 \\
		\bv a_2
\end{bmatrix},  \,\bar{\bv{a}}=
\begin{bmatrix}
		\bar{\bv{a}}_1 \\
		\bar{\bv{a}}_2
\end{bmatrix}.
\end{align*}
\end{lmm}
\begin{proof}
By Theorem~\ref{th:comred-LegTrafo-inj-2}, the coordinate transformation $\hat{\mathfrak{z}}: \mathbb{A}_n\rightarrow \mathbb{R}^n$ to the energy variable is well-defined. The functional $\tilde{\HamCo}: \mathbb{R}^n \rightarrow \mathbb{R}$, $\tilde{\HamCo}(\hat{\mathfrak{ {z} }}(\bv a)) = \HamCo(\bv a)$ thus represents the Hamiltonian as a function of the energy variable. As it is strictly convex by construction (due to Assumption~\ref{assum:h-sconvex}), it holds 
\begin{align*}
	\tilde{\HamCo}(\bar{\bv z}) - \tilde{\HamCo}(\bv z) >
	\nabla \tilde{\HamCo}(\bv{z})\cdot
	\left(  \bar{\bv z} - \bv z \right) ,\hspace{1cm} \text{for } \bv z, \bar{\bv{z}} \in \hat{\mathfrak{ z }}(\mathbb{A}_n), \, \bv z \neq \bar{\bv{z}} .
\end{align*}
The claimed inequality is equivalent to the latter one with ${\bv{z}}= \hat{\mathfrak{z}}({\bv{a}})$, $\bar{\bv{z}}= \hat{\mathfrak{z}}(\bar{\bv{a}})$, as can be seen by further employing that $\nabla \tilde{\HamCo}(\bv z)_{|\bv z = \hat{\mathfrak{z}}({\bv{a}})} = [	-\nabla_1 \LegCo(\bv{a}) ;\bv a_2 ]$, cf.\ Lemma~\ref{lem:legOfleg}.
\end{proof}

\begin{thrm}[Energy-dissipation inequality]\label{theor:energy-bound-ImplicitEuler}
For any solution of System~\ref{sys:twopH-gallapp-ImplicitEuler}, it holds for $0\leq \ell < k \leq K$ and $\bv{a}^\ell \neq \bv{a}^k$,
\begin{align*}
	 \HamCo(\bv{a}^k) -  \HamCo(\bv{a}^\ell) < \Delta_t \left( \sum_{j=\ell+1}^k \bcef^j \cdot \bcfl^j - (\bv{a}_2^j)^T \bt{R}(\bv{a}^j) \bv{a}_2^j   \right) \leq   \Delta_t \sum_{j=\ell+1}^k \bcef^j \cdot \bcfl^j.
\end{align*}
\end{thrm}
\begin{proof}
	Clearly, it is sufficient to consider the case $\ell = k-1$, as $\ell<k-1$ can be directly concluded from it. By Lemma~\ref{lem:HamConvexCond} it holds
\begin{align*}
	\HamCo(\bv a^k) - \HamCo(\bv a^{k-1}) <
	\begin{bmatrix}
	-\nabla_1 \LegCo(\bv{a}^k) \\
	\bv a_2^k 
	\end{bmatrix} \cdot \left(
	\begin{bmatrix}
		\bv a_1^k \\
		\nabla_2 \LegCo(\bv{a}^k)
	\end{bmatrix} -	
		\begin{bmatrix}
		\bv{a}_1^{k-1} \\
		\nabla_2 \LegCo(\bv{a}^{k-1})
	\end{bmatrix}
	 \right).
\end{align*}
Testing System~\ref{sys:twopH-gallapp-ImplicitEuler} with $\bv b = [-\nabla_1 \LegCo(\bv{a}^k); \bv{a}_2^k]$ and using 
$ \nabla^\bt{M}_i \LegCo(\bv{a}^k) = \bt{M}_i^{-1}\nabla_i \LegCo(\bv{a}^k)$ for $i=1,2$, and $\bt K^T [\tilde{\bv a}_1; \bv a_2^k] = \bt K_2^T  \bv a_2^k = \bv f^k$ (which holds independently of $\tilde{\bv a}_1$), it follows
\begin{align*}
	\frac{1}{\Delta_t} 
	\begin{bmatrix}
	-\nabla_1 \LegCo(\bv{a}^k) \\
	\bv a_2^k 
	\end{bmatrix} \cdot&
\left(
		\begin{bmatrix}
		\bv a_1^k \\
		\nabla_2 \LegCo(\bv{a}^k)
	\end{bmatrix} -
		\begin{bmatrix}
		\bv{a}_1^{k-1} \\
		\nabla_2 \LegCo(\bv{a}^{k-1})
	\end{bmatrix}
	 \right) 
	  \\
	 	&=\begin{bmatrix}
	-\nabla^{\bt{M}}_1 \LegCo(\bv{a}^k) \\
	\bv a_2^k 
	\end{bmatrix} \cdot \left(
	 \begin{bmatrix}
	 & \bt{J} \\
	 -\bt{J}^T & -\bt{R}(\bv{a}^k)
\end{bmatrix} 
\begin{bmatrix}
	-\nabla^{\bt{M}}_1 \LegCo(\bv{a}^k) \\
	\bv{a}_2^k
\end{bmatrix}
+ \bt{K}
 \bcef^k \right) \\
 &= \bcef^k \cdot \bcfl^k - (\bv{a}_2^k)^T \bt{R}(\bv{a}^k) \bv{a}_2^k  \leq \bcef^k \cdot \bcfl^k,
\end{align*}
due to the positive semi-definiteness of $\bt{R}(\bv{a})$. Setting together the two inequalities yields the assertion.
\end{proof}

The proposed time discretization is energy-dissipative. In contrast to the energy dissipation equality in the time-continuous case (Theorem~\ref{theor:energy-bound-galapp}) we face here a dissipation inequality. The difference is due to the numerical dissipation of the implicit Euler-type method, which is $D(t_k) =   \HamCo(\bv{a}^0)- \HamCo(\bv{a}^k)+\Delta_t ( \sum_{j=1}^k \bcef^j \cdot \bcfl^j - (\bv{a}_2^j)^T \bt{R}(\bv{a}^j) \bv{a}_2^j ) $. Certainly, more sophisticated energy-stable or energy-preserving schemes can also be derived under the use of the underlying Hamiltonian structure, cf., \cite{phd:liljegren, inprod-MehrmannM19, book:hairerGeomInt, art:topics-in-strpresdisc}, but this is beyond the scope of this paper.


\section{Application} \label{chap:p1b:isothEuler}

The one-dimensional barotropic Euler equations and various simplifications of them are important representatives of our model problem. They are used, e.g., for gas network simulations, cf., \cite{phd:liljegren,inproc:bls-hyp2018,art:hierarchGasMindt2019}. In this section we demonstrate the feasibility of our approximation approach at the examples of an undamped dam-break test case with a shock and of a larger network in a friction-dominated regime relevant for pipelines in a gas transport network. For the last example, weighted edges must be considered which requires corresponding minor adjustments to the presented results, as we will briefly comment.

\subsection{Barotropic Euler equations}
We consider the barotropic Euler equations with an additional friction term as an example for our model problem. Flow density and velocity $\rho,v: [0,T] \times  \Omega \rightarrow \mathbb{R}$ are governed by
\begin{align}  \label{eq:Part1b-ISO1}
\partial_t \rho + \partial_x \left( \rho v \right) = 0 , \qquad 
\partial_t v + \partial_x \frac{v^2}{2} + \frac{1}{\rho} \partial_x p(\rho) = - \frac{\lambda}{2 D }  |v| v, \hspace{0.4cm} x\in \Omega, \, t\in [0,T],
\end{align}
where the pressure $p$ is only a function of density. Furthermore, $\lambda$ describes the friction factor that might be state-dependent and $D$ the constant diameter of the pipe. In networks the circular cross-sectional pipe area $A= \pi/4\, D^2$  usually acts as a weighting term for the respective edge $\Onepipe \in \mathcal{E}$, see Appendix~\ref{app:weighted-edges}. Following \cite{book:NovotnyCompressible, art:winters-euler-entropy}, we introduce $P$ as a pressure potential that is characterized by the relation $P''(\rho) =  p'(\rho)/\rho$. Examples for compatible pairs of $p$ and $P$ are given in the test cases. Using the entropy as Hamiltonian $\HamPDE$, system~\eqref{eq:Part1b-ISO1} can be formally rewritten in our standard form~\eqref{bls-eq:abstr} for the energy variable $ \ubar{\bv{z}} =[\rho;v]$, i.e.,
\begin{align*} 
	\partial_t
	\begin{bmatrix}
		\rho \\
		v
	\end{bmatrix}
	 &=
	\begin{bmatrix}
		 & -\partial_x \\
		 -\partial_x & -\tilde{r}([\rho;v])
	\end{bmatrix}
	 \derst{} h([\rho;v])
	 \end{align*}
with Hamiltonian density $h([\rho;v])= \rho {v^2}/{2} + P(\rho)$ and non-negative friction-associated function $\tilde{r}([\rho;v])= {\lambda |v|}/{(2 D \rho) }$.
The coupling conditions for $\blsVertex \in {\mathcal{N}}_0$ take the form
\begin{align*} 
	\sum_{\Onepipe \in \mathcal{E}(\blsVertex) } n^{\Onepipe}[\blsVertex]  (\rho v)_{|\Onepipe}[\blsVertex] &= 0, \qquad
	P'(\rho_{|\Onepipe})[\blsVertex] + \frac{v_{|\Onepipe}[\blsVertex]^2}{2}  = P'(\rho_{|\tilde\Onepipe})[\blsVertex] + \frac{v_{|\tilde\Onepipe}[\blsVertex]^2}{2}
\end{align*}
for $\Onepipe,\tilde{\Onepipe} \in \mathcal{E}(\blsVertex)$. The coupling conditions are entropy-preserving. Note that it is not trivial to find the entropy-preserving conditions without energy-based modeling or other sophisticated analytic tools, cf., \cite{art:euler-reigstad,art:egger-mfem-compressEuler,art:hierarchGasMindt2019}. The system is closed by appropriate initial and boundary conditions that are specified in the test cases.

In our approach, the system is parametrized in the variables $z_1=\rho$ and $\nabla_2 h(\bv z)=\rho v$, i.e., $\ubar{\bv{a}} = [\rho;m]$ with $m = \rho v$ denoting the mass flux.  In the hyperbolic theory these variables are a natural choice for the barotropic Euler equations that is typically used when a conservative formulation is discretized, cf., \cite{book:leveque-finite-volume-methods,art:winters-euler-entropy,art:domschke-adj-based2015,art:euler-reigstad}. Although we do not approach the problem from the hyperbolic point of view, we use the variables as we rely on the port-Hamiltonian framework. Our proposed space discretization reads as follows.

\begin{system} \label{sys:twopH-cr-iso1}
Find $\ubar{\bv{a}}=[\rho;m] \in \funSpace{C}^1([0,T];\funSpace{V}_1 \times \funSpace{V}_2)$, $\bcfl \in \funSpace{C}^1([0,T];\mathbb{R}^{\nbn})$ and $\bcef  \in \funSpace{C}([0,T];\mathbb{R}^{\nbn})$ such that 
\begin{align*} 
 \Lscal \partial_t \rho, b_1 \Rscal &= - \Lscal \partial_x m, b_1 \Rscal   \\
  \left\Lscal \partial_t \frac{m}{\rho}, b_2 \right\Rscal &= \left\Lscal P'(\rho)+ \frac{m^2}{2\rho^2}, \partial_x b_2 \right\Rscal + \bcef \cdot \TraceOp b_2 - \left \Lscal r([\rho;m]) m , b_2 \right \Rscal   \\
  \bcfl   &=  \TraceOp m
\end{align*}
with $ r([\rho;m]) = \lambda |m|/(2D \rho^2)$ hold for all $b_1 \in \funSpace{V}_1$ and $b_2 \in \funSpace{V}_2$. The system is closed with initial and boundary conditions, $\ubv a(0)=\ubv a_0$ and $\bv{k}(\bv{e}(t),\bv{f}(t),\bv{u}(t)) = 0$ for $t\in[0,T]$, for given 
$\ubv a_0 \in \funSpace{V}$ and $u_\nu : [0,T]\rightarrow \mathbb{R}$ for $\nu \in \mathcal{N}_\partial$.
In the case of complexity reduction, the inner products $\Lscal \cdot, \cdot \Rscal$ in the second equation are replaced by $\Lscal \cdot, \cdot \Rscal_c$.
\end{system}

\begin{rmrk}\label{rem:fem-egger}
Note that the preservation of port-Hamiltonian structure is a sufficient but not a necessary condition for a discrete energy bound (as in Theorem~\ref{theor:energy-bound-galapp}) to hold. A similar Galerkin approximation for the barotropic Euler equations that fulfills an energy bound but does not preserve this structure can be found in \cite{art:egger-mfem-compressEuler}. The underlying variational principle differs in the second equation, which reads in \cite{art:egger-mfem-compressEuler}
\begin{align*} 
\left\Lscal \partial_t \frac{m}{\rho} + { \frac{m \partial_t \rho }{2 \rho^2}}, b_2 \right\Rscal &= \left\Lscal P'(\rho)+ \frac{m^2}{2\rho^2}, \partial_x b_2 \right\Rscal + \bcef \cdot \TraceOp b_2 - \left \Lscal r([\rho;m]) m { + \frac{m}{2 \rho^2}\partial_x m} , b_2 \right \Rscal.
\end{align*}
In the reference no complexity reduction is considered and the derivation of an energy-stable time discretization seems to be more involved due to the loss of port-Hamiltonian structure. Note that port-Hamiltonian structure is also desirable for other reasons, e.g., for relative error estimates \cite{art:gieselmann-rel-energy}, or when using control theory \cite{art:pcH-maschke92,art:SchJ14,art:weak-pH-koty}.
\end{rmrk}

\subsection{Numerical studies} \label{subsec:realize-num}

Our numerical studies are performed with mixed finite elements in space and the proposed implicit Euler-type time discretization (System~\ref{sys:twopH-gallapp-ImplicitEuler}) with a uniform step size. The resulting nonlinear algebraic systems are solved by the Newton's method. When using finite elements of the lowest order, i.e., $\funSpace{V}_1 = \funSpace{Q}_{0}(T_{\mathcal{E}})$, $\funSpace{V}_2 = \funSpace{P}_{1}(T_{\mathcal{E}})$, the nonlinear integrals can be evaluated analytically, as they consist of polynomial expressions for the example. Higher-order finite elements are supplemented with complexity reduction by Gaussian quadrature of sufficiently high degree for Assumption~\ref{assum:quadrat-FIRST} to hold.
All results have been generated using \texttt{MATLAB} Version 9.1.0 (R2016b) on an Intel Core i5-7500 CPU with 16.0GB RAM. For better reproducibility, the code and the benchmark data are provided in \cite{code:bls22-phapprox}.

\begin{figure}[tb]
\centering
\begin{tabular}{rlll}
\begin{minipage}{0.022\textwidth}
{
{\small\rotatebox{90}{density $\rho$}\\ 

}}
\end{minipage}
&
{\hspace{-0.8cm}
\begin{minipage}{0.4\textwidth}
\center
\hspace{0.1cm} \includegraphics[height = 0.64\textwidth, width = 0.77\textwidth]{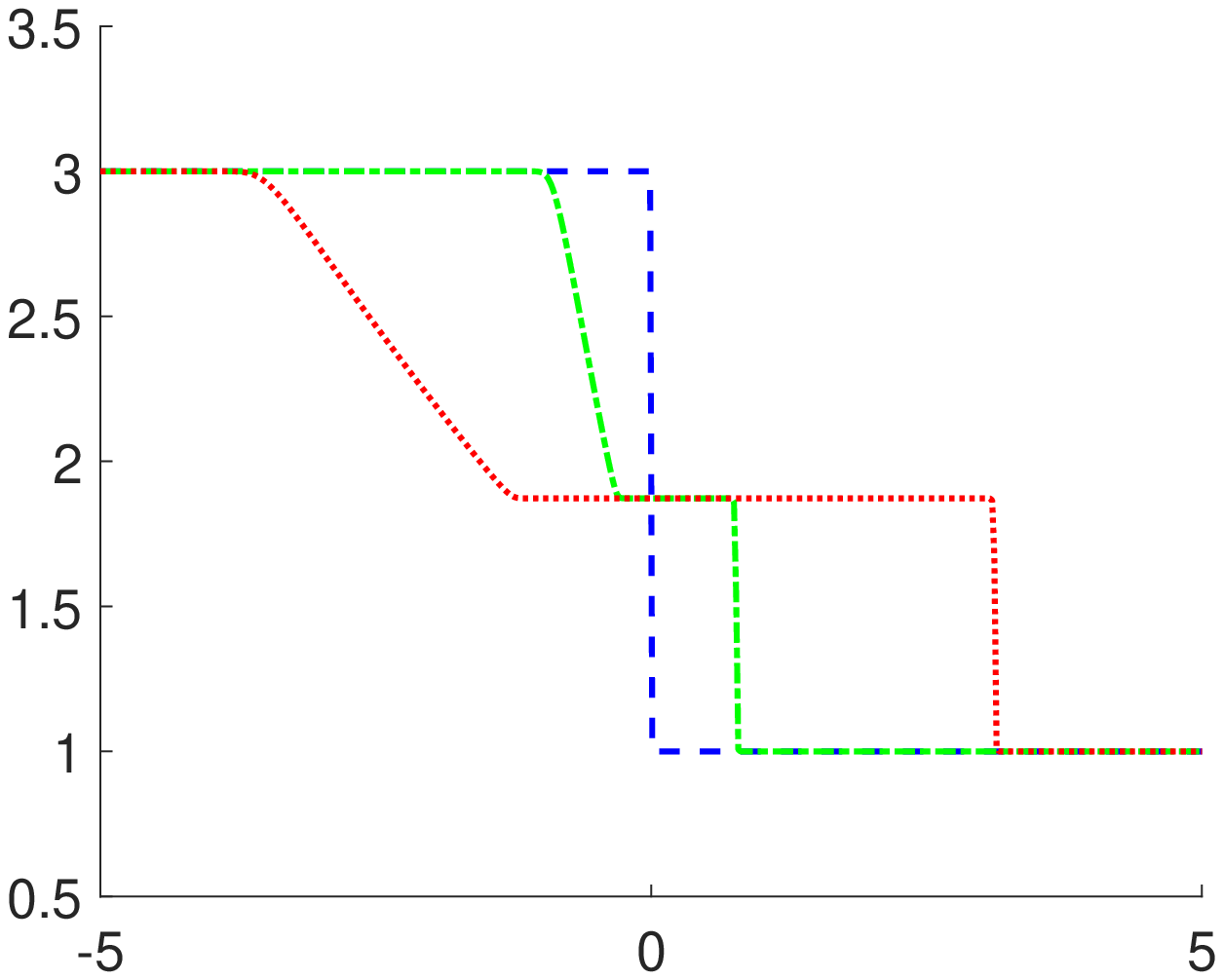} \\
 \hspace{0.5cm} {\small space $x$}
\end{minipage}
}
&
\begin{minipage}{0.022\textwidth}
{\hspace{0.5cm}
{\normalsize\rotatebox{90}{mass flux  $m$ \vspace{0.4cm} } \\ 
}}
\end{minipage}
&
\hspace{-0.7cm}
\begin{minipage}{0.4\textwidth}
\center
\hspace{0.1cm} \includegraphics[height = 0.64\textwidth, width =0.82\textwidth]{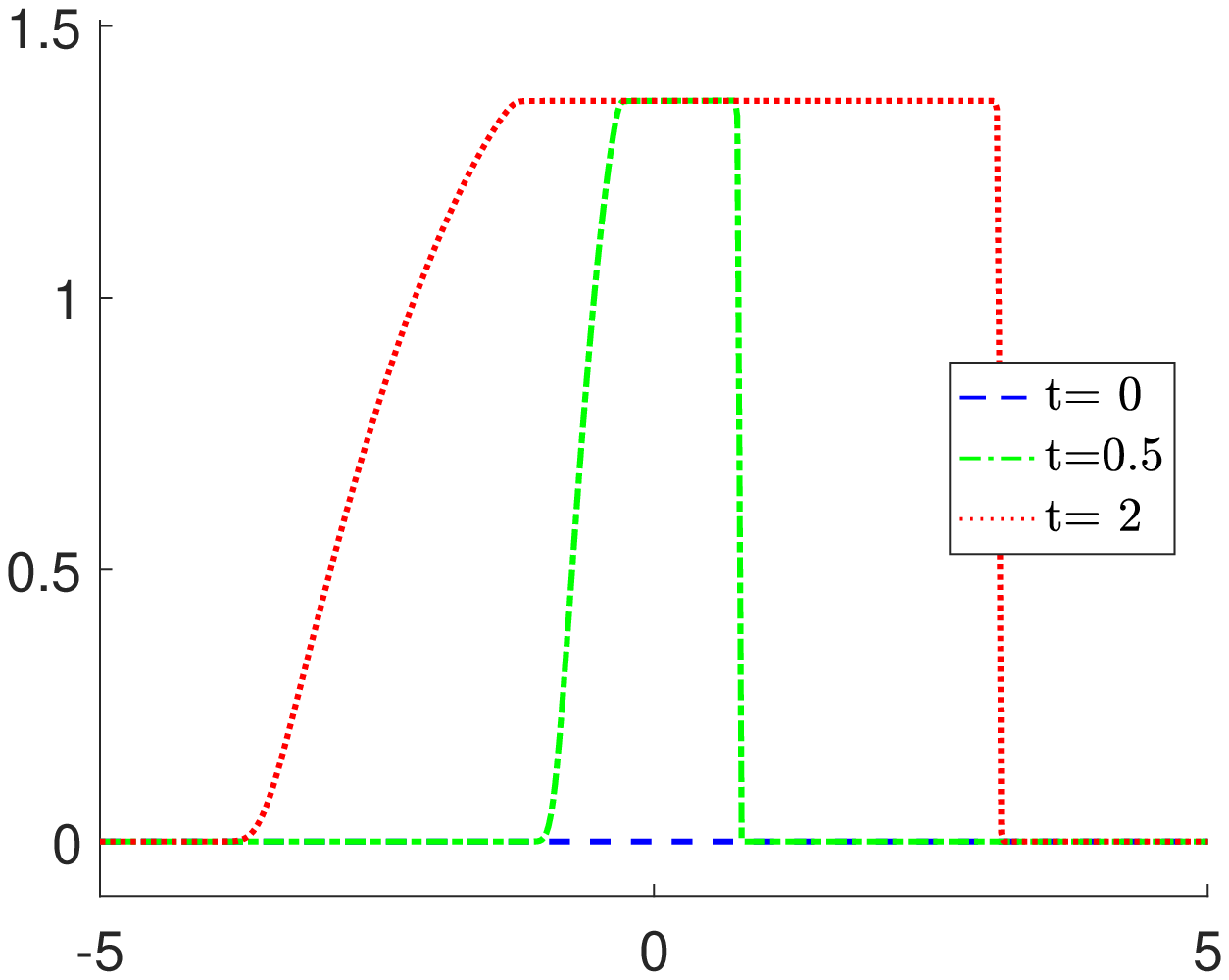} \\
 \hspace{0.5cm} {\small space $x$}
\end{minipage}
\end{tabular}
\caption{Dam-break test case. Spatial representation of solution behavior for three snapshots in time,  generated with  $\funSpace{V}_1 = \funSpace{Q}_{0}(T_{\mathcal{E}})$, $\funSpace{V}_2 = \funSpace{P}_{1}(T_{\mathcal{E}})$, $\Delta_x=0.001$  and $\Delta_t=0.005$. \label{fig:dam-break}}
\end{figure}
\begin{figure}[tb]
\centering
\begin{tabular}{rlll}
\begin{minipage}{0.022\textwidth}
{
{\small\rotatebox{90}{relative error}\\ 

}}
\end{minipage}
&
{\hspace{-1.2cm}
\begin{minipage}{0.44\textwidth}
\center
\hspace{0.1cm} \includegraphics[height = 0.58\textwidth, width = 0.73\textwidth]{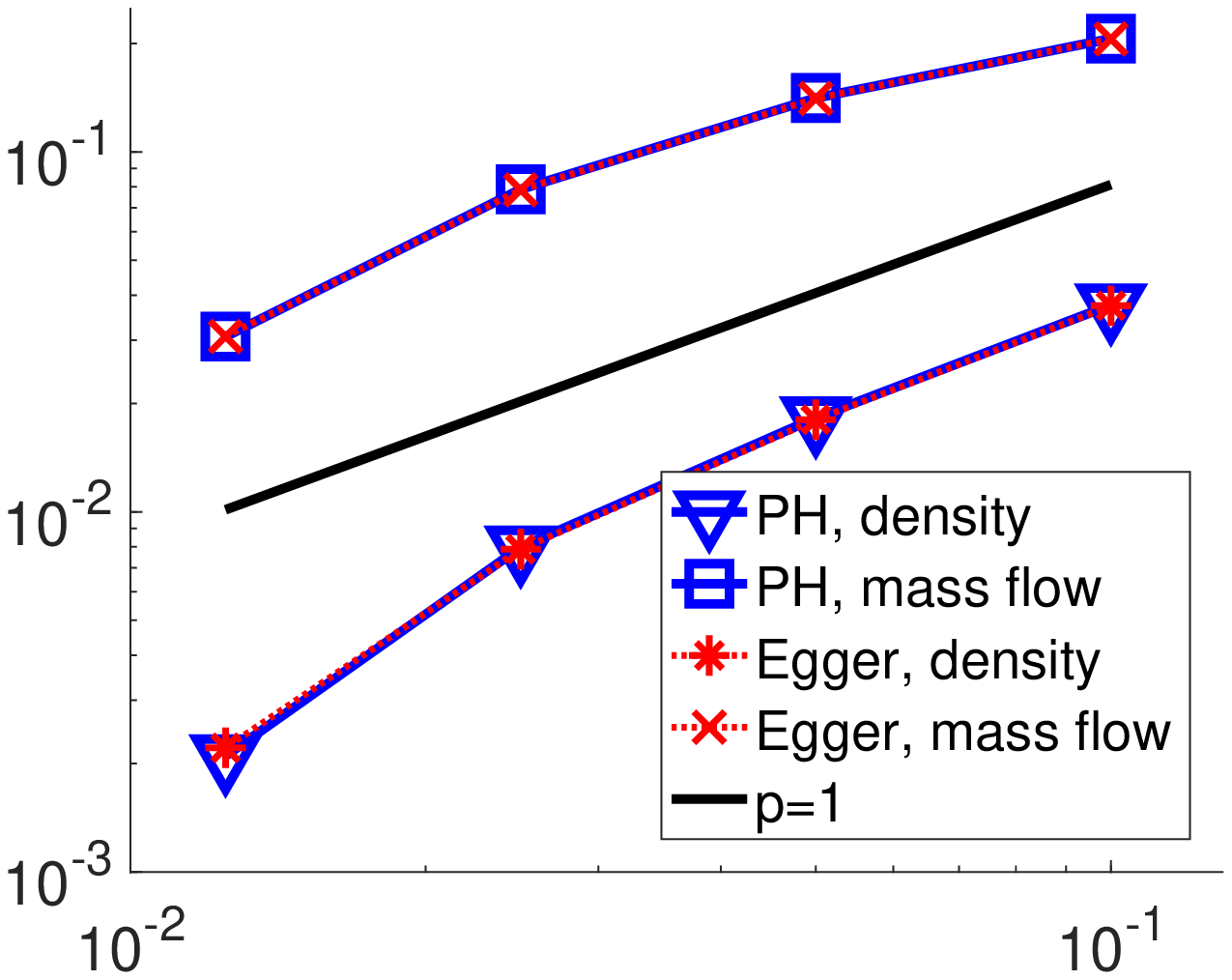} \\
 \hspace{0.5cm} {\small step size $\Delta_x$}
\end{minipage}
}
&
\begin{minipage}{0.032\textwidth}
{\hspace{0.5cm}
{\normalsize\rotatebox{90}{energy loss\vspace{0.1cm} } \\ 
}}
\end{minipage}
&
\hspace{-1.42cm}
\begin{minipage}{0.44\textwidth}
\center
\hspace{0.1cm} \includegraphics[height = 0.58\textwidth, width = 0.73\textwidth]{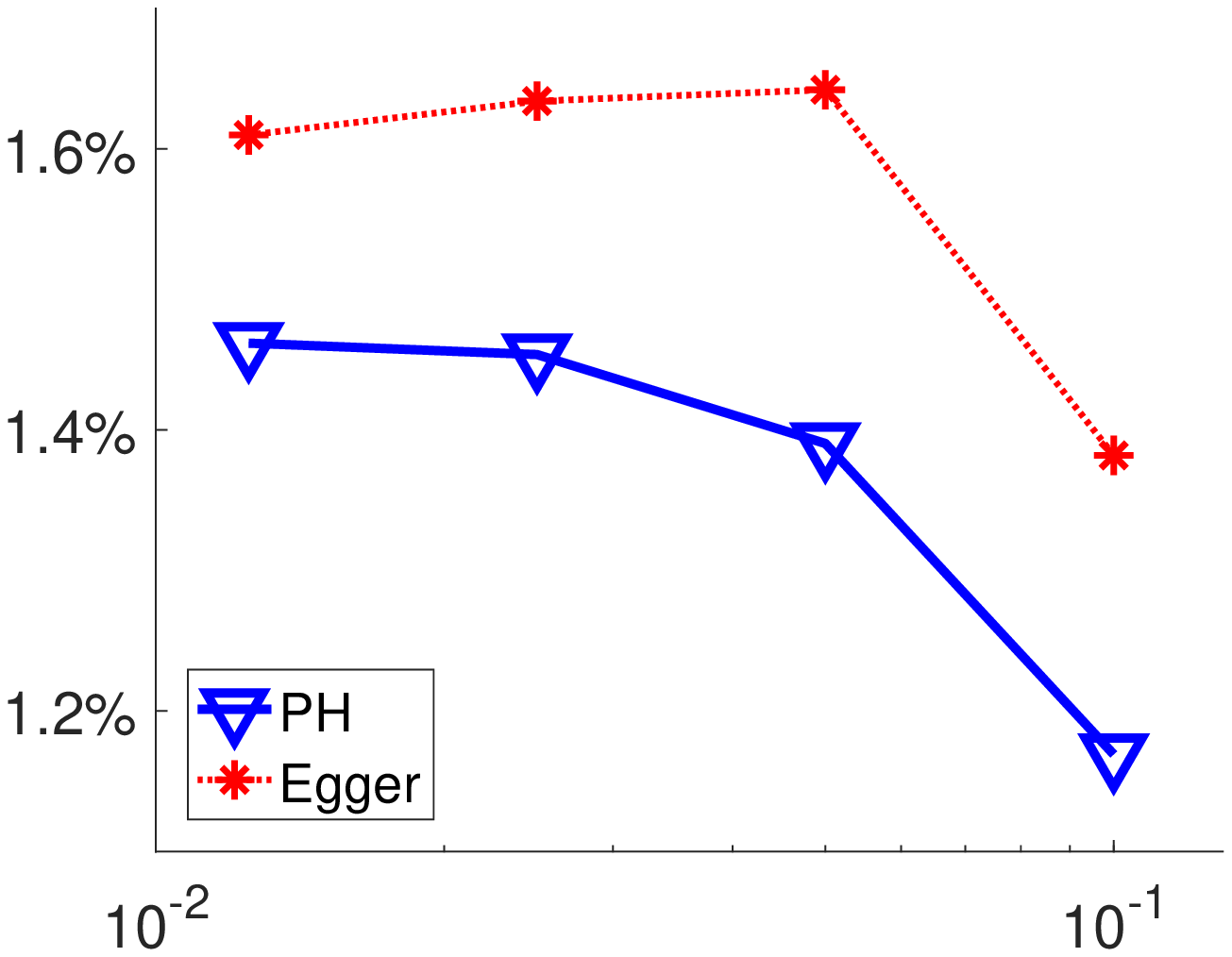} \\
 \hspace{0.5cm} {\small step size $\Delta_x$}
\end{minipage}
\end{tabular}

\caption{Dam-break test case. Comparison of our proposed method (PH) and the method from~\cite{art:egger-mfem-compressEuler} (Egger), cf., Remark~\ref{rem:fem-egger}. \textit{Left:} $\funSpace{L}^2$-errors for spatial domain $[-5,0)$, reference solution obtained by PH with $\Delta_x =  0.001$. \textit{Right:} Energy loss in Hamiltonian, given by $1-\HamPDE(\underline{\bv{a}}(T))/\HamPDE(\underline{\bv{a}}(0))$. Both for end time $T=2$ and $\Delta_t= 0.0005$. \label{fig:dam-err-energy}}
\end{figure}

\begin{figure}[tb]
\centering
\begin{tabular}{rlll}
\begin{minipage}{0.022\textwidth}
{
{\small\rotatebox{90}{density $\rho$}\\ 

}}
\end{minipage}
&
{\hspace{-0.8cm}
\begin{minipage}{0.4\textwidth}
\center
\hspace{0.1cm} \includegraphics[height = 0.67\textwidth, width = 0.95\textwidth]{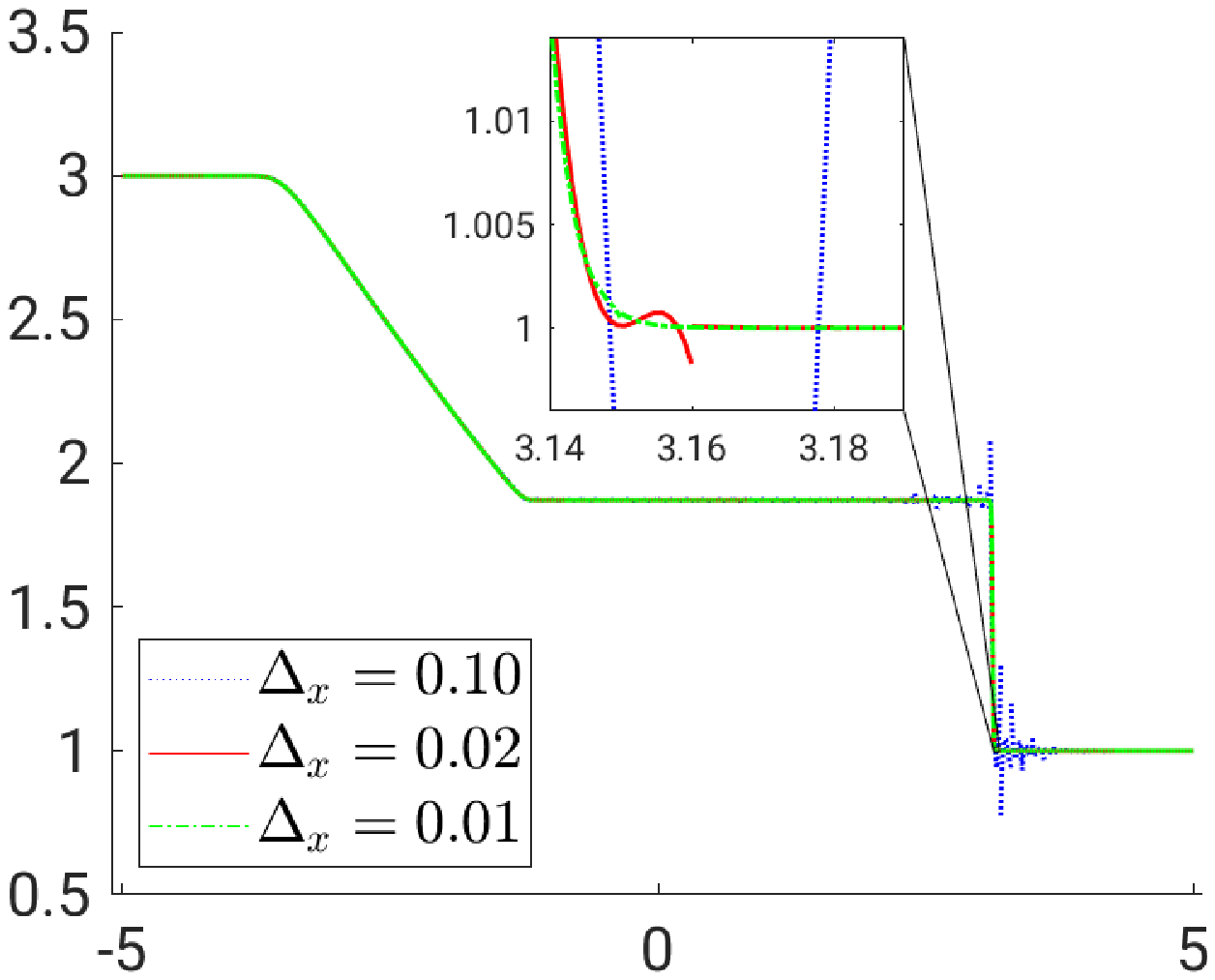} \\
 \hspace{0.5cm} {\small space $x$}
\end{minipage}
}
&
\begin{minipage}{0.032\textwidth}
{\hspace{0.5cm}
{\normalsize\rotatebox{90}{mass flux  $m$ \vspace{0.4cm} } \\ 
}}
\end{minipage}
&
\hspace{-0.7cm}
\begin{minipage}{0.4\textwidth}
\center
\hspace{0.1cm} \includegraphics[height = 0.67\textwidth, width = 1.00\textwidth]{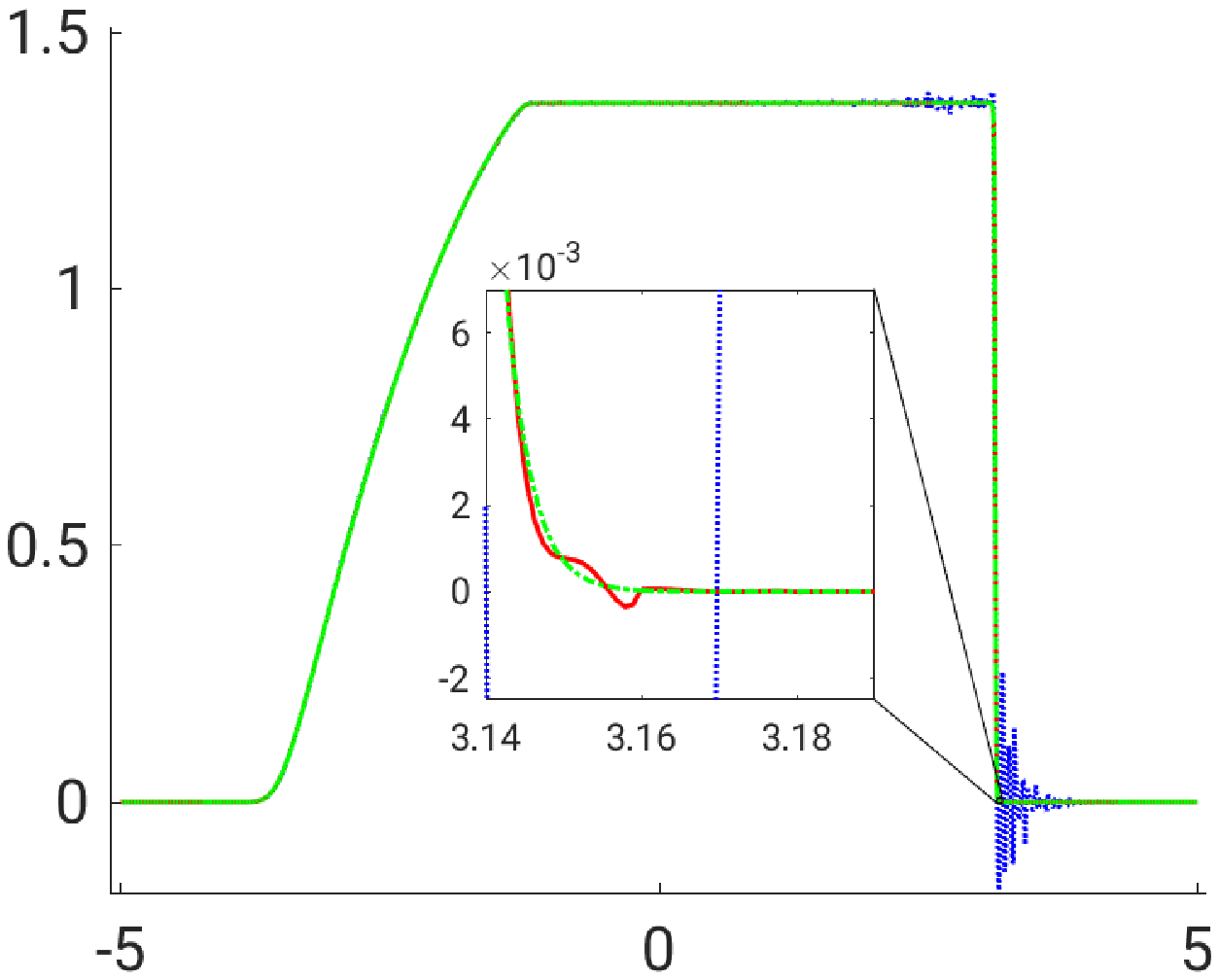} \\
 \hspace{0.5cm} {\small space $x$}
\end{minipage}
\end{tabular}

\caption{Dam-break test case. Spatial representation of solution at $T=2$ obtained by higher-order finite elements for varying element sizes $\Delta_x$ with zoom into region around shock ($\Delta_t= 0.005$). \label{fig:dam-breakHO}}
\end{figure}

\subsubsection{Dam-break test case}
The undamped benchmark example on one edge is taken from \cite{book:leveque-finite-volume-methods,art:egger-mfem-compressEuler} and referred to as dam-break problem. We use it here as a numerical stability test for our approach and therefore purposely omit to add a shock-capturing mechanism. For the latter we refer to, e.g., \cite{art:Nordstrom14-SumByParts,book:leveque-finite-volume-methods} and note that (numerical) dissipation can be included in our analysis without difficulty, see Appendix~\ref{subsec:dissip}. The model parameters are chosen as $\lambda=0$, $A=1$ and $p(\rho) = 0.5 \rho^2$, which corresponds to the isentropic Euler equations with pressure potential $P(\rho) = p(\rho)$. The spatial domain is $\Omega=[-5,5]$, and the initial and boundary conditions read
\begin{align*}
	\rho(0,x) = 2- \text{sgn}(x), \quad m(0,x) = 0 \quad  x \in \Omega, \hspace{0.3cm} \text{ and } \hspace{0.3cm}  m(t,-5) = m(t,5) =0, \quad  t \geq 0.
\end{align*}
As end time we take $T=2$. The initial state has a discontinuity in the density $\rho$ at $x=0$, which propagates from left to right as a shock wave, whereas a rarefaction wave moves in the opposite direction, see Fig.~\ref{fig:dam-break} for an illustration.

For this setup, we compare our method with the finite element discretization proposed in \cite{art:egger-mfem-compressEuler} (cf., Remark~\ref{rem:fem-egger}), using the spaces $\funSpace{V}_1 = \funSpace{Q}_{0}(T_{\mathcal{E}})$ and $\funSpace{V}_2 = \funSpace{P}_{1}(T_{\mathcal{E}})$, see Fig.~\ref{fig:dam-err-energy}.  Both space discretizations show a first-order convergence outside the region of the shock ($x \in [-5, 0)$) and overall very similar approximation quality. But our method is less dissipative, e.g., for $\Delta_x=0.05$ an energy loss of less than $1.4\%$ is observed, i.e., $\HamPDE(\ubv{a}(2))\approx 0.986 \, \HamPDE(\ubv{a}(0))$, in contrast to the other one with a loss of more than $1.6\%$. 

Further, we use the dam-break problem to showcase the entropy stability for higher-order discretization combined with inexact integration. Figure~\ref{fig:dam-breakHO} exemplarily illustrates the results for the end time $T=2$ and the choice $\funSpace{V}_1 = \funSpace{Q}_{3}(T_{\mathcal{E}})$, $\funSpace{V}_2 = \funSpace{P}_{4}(T_{\mathcal{E}})$ supplemented with a five-point Gaussian quadrature of the nonlinear integrals. The size of the elements $\Delta_x$ is varied, otherwise the same parameters as above are used. As to be expected due to the Gibbs phenomenon \cite{art:Nordstrom14-SumByParts, book:arnold2006-comp}, oscillations occur near the shock when the solution is under-resolved, which is the case for $\Delta_x = 0.1$ in this example. Nonetheless, the discrete solution stays bounded because of the entropy stability our method inherits. For the smaller grid sizes $\Delta_x = 0.02$ and $\Delta_x = 0.01$, the numerical dissipation coming from the time discretization seems to be sufficient to filter out the high frequency oscillations almost completely, as can be seen in the enlarged image around the shock in Fig.~\ref{fig:dam-breakHO}.

\subsubsection{Gas pipeline network}

\begin{figure}[tb]
\,{\small \hspace{0.2cm} \underline{input profile $u$} \hspace{1.5cm} \underline{pipe-end 2} \hspace{1.2cm}	\underline{pipe-end 4} \hspace{1.2cm} \underline{pipe-end 7}} \\
\begin{tabular}{l | rl}
\begin{minipage}{0.18\textwidth}
\includegraphics[width=0.96\textwidth, height = 0.86\textwidth]{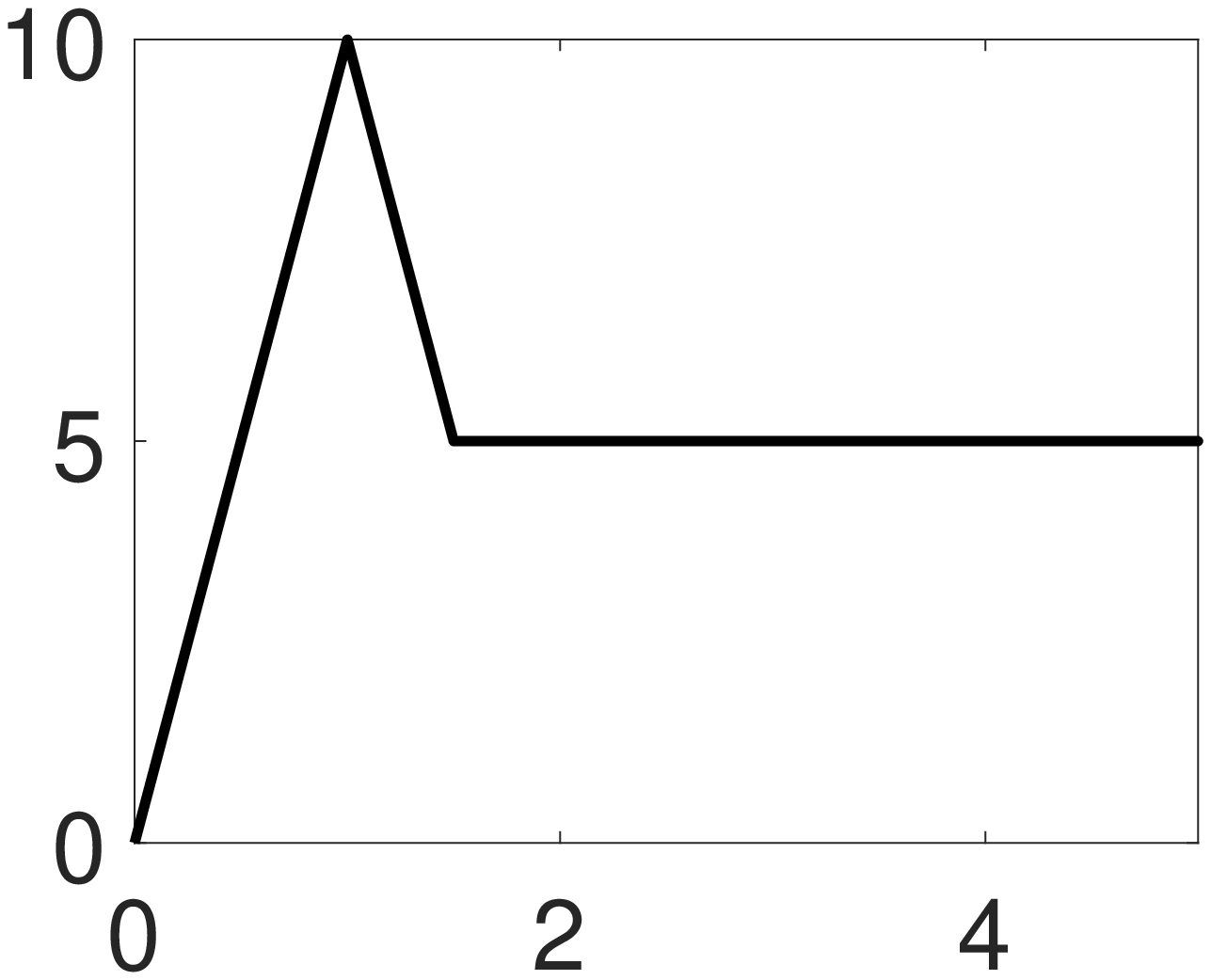}\\
 {\hspace*{0.35cm} \small time $t$ [h]}  
\vspace{2.6cm}
\end{minipage}
&
\begin{minipage}{0.022\textwidth}
{\hspace{-0.0cm}
{\small\rotatebox{90}{density $\rho$}\\ \vspace{0.2cm}\\
\rotatebox{90}{ mass flow $A m$} 
}}
\end{minipage}
&
{\hspace{-0.4cm}
\begin{minipage}{0.71\textwidth}
\includegraphics[width=0.28\textwidth]{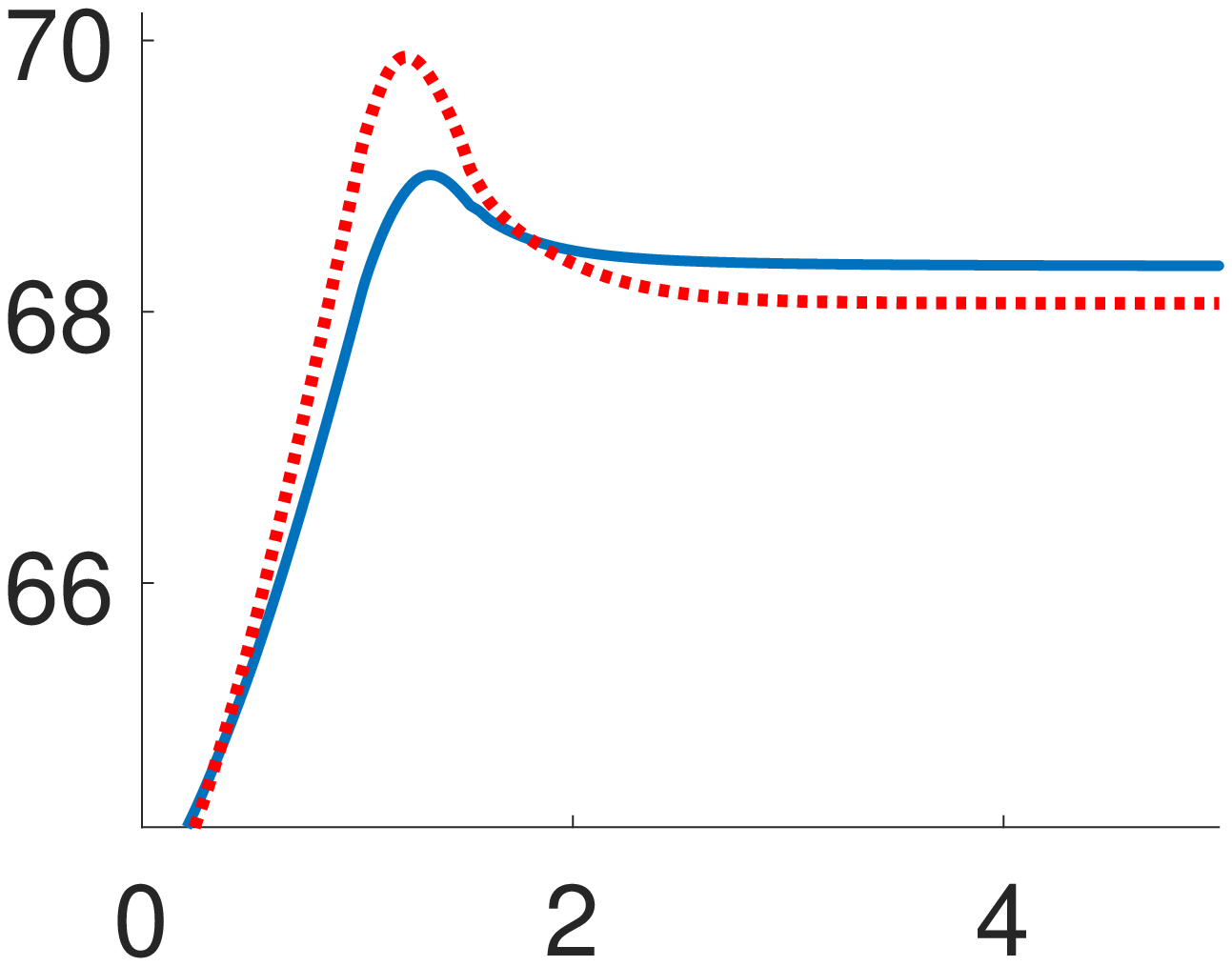} \hspace{0.1cm} 
\includegraphics[width=0.28\textwidth]{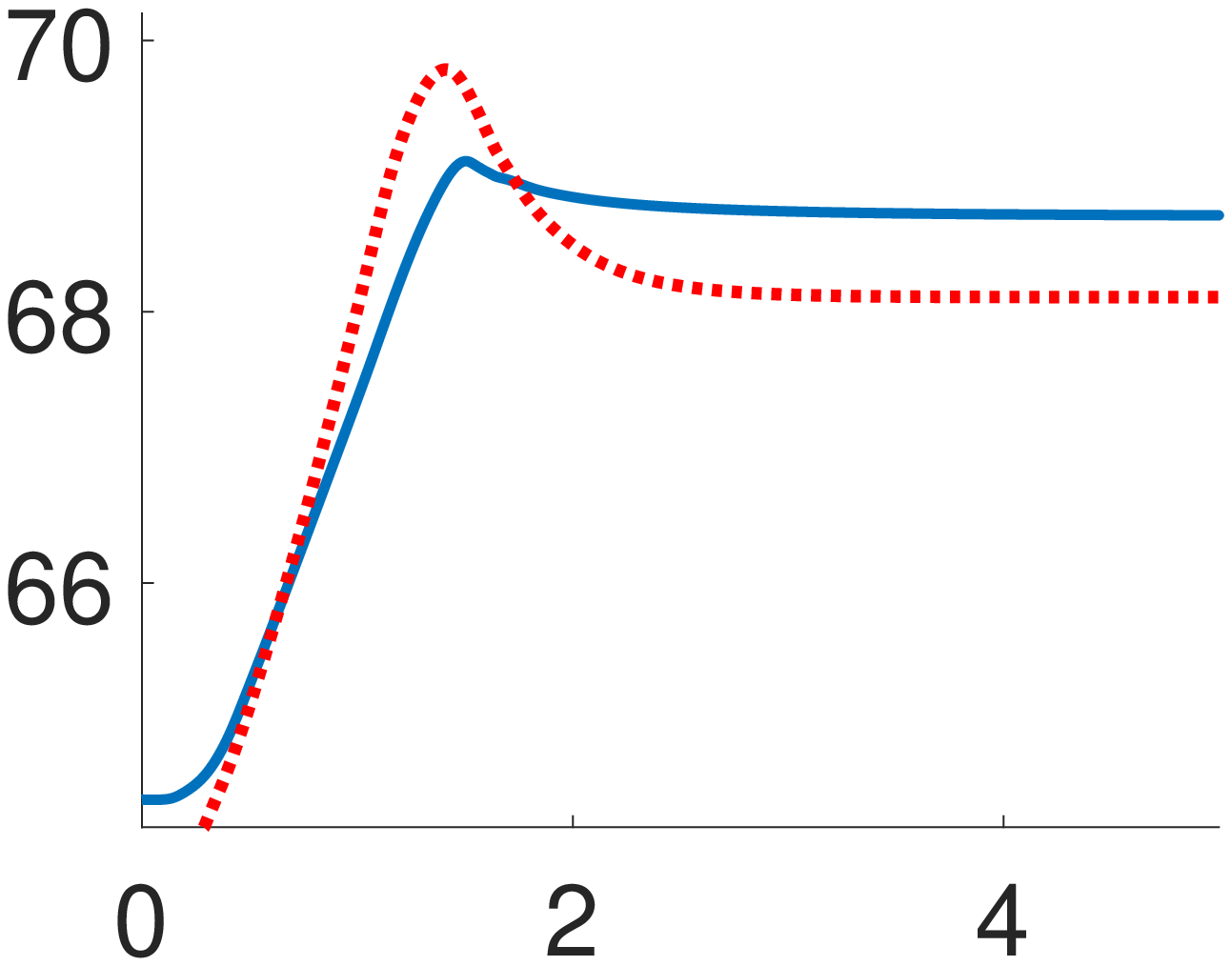}  \hspace{0.1cm} 
\includegraphics[width=0.28\textwidth]{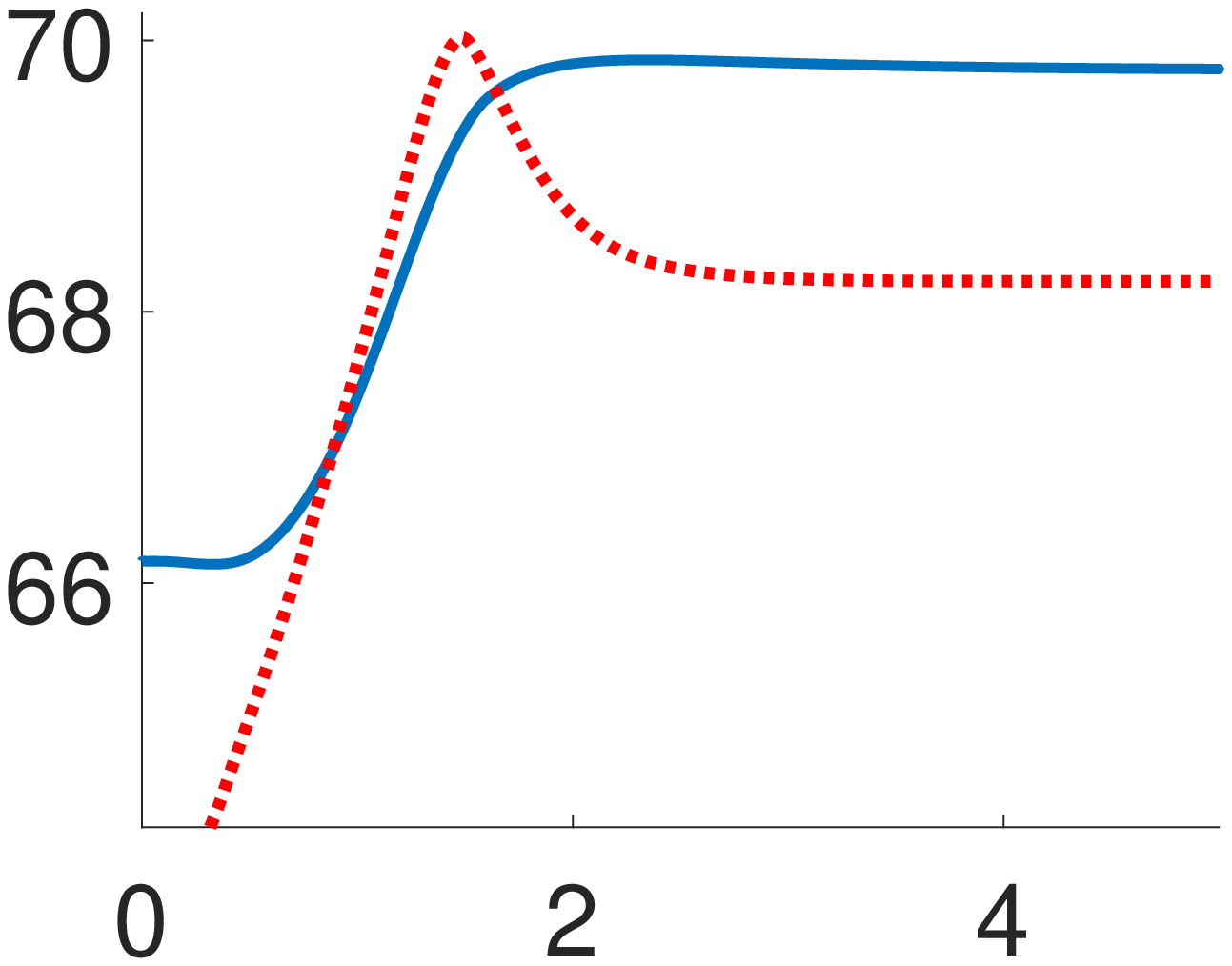}\\
\vspace{0.1cm} \\
\includegraphics[width=0.28\textwidth]{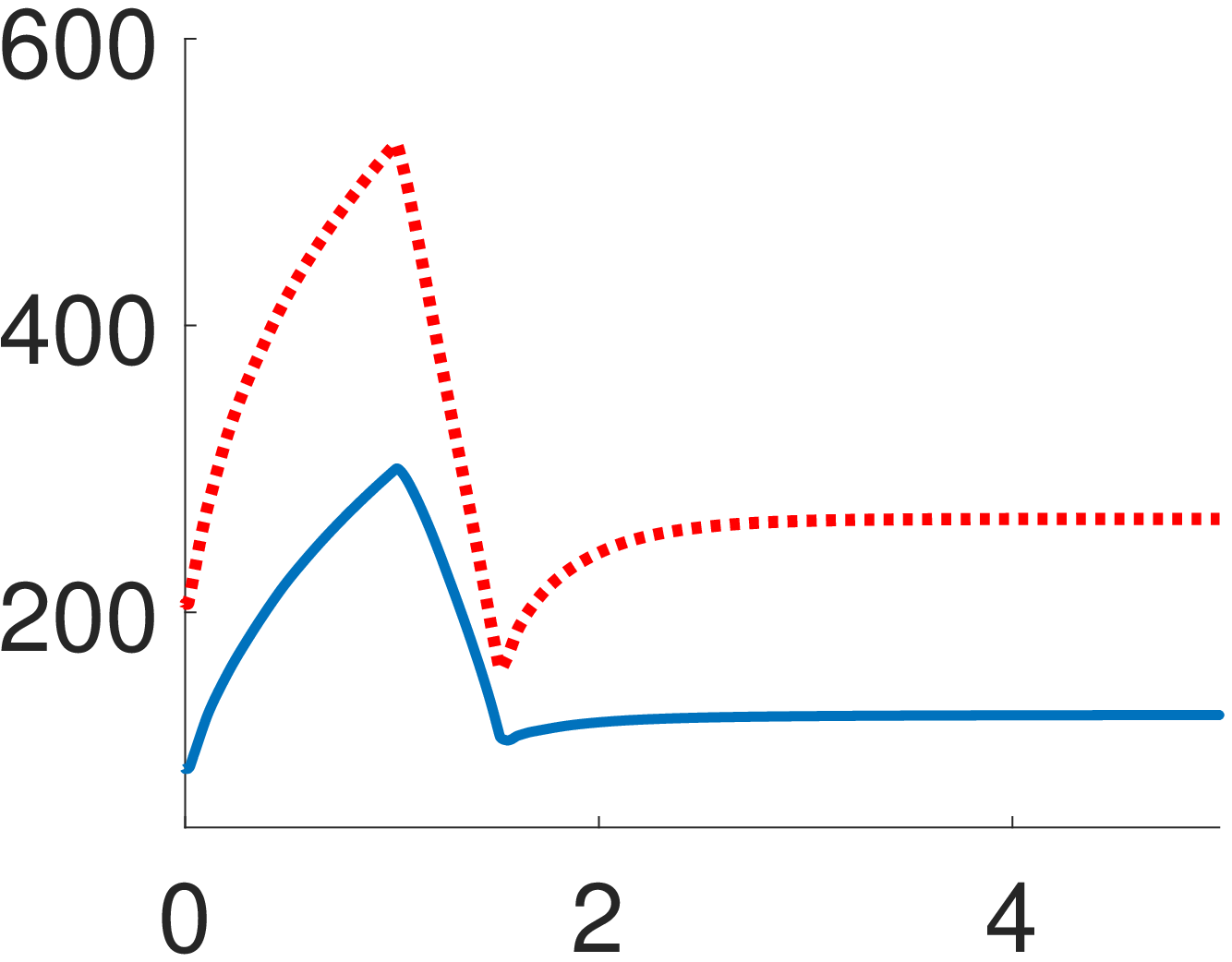} \hspace{0.1cm} 
\includegraphics[width=0.28\textwidth]{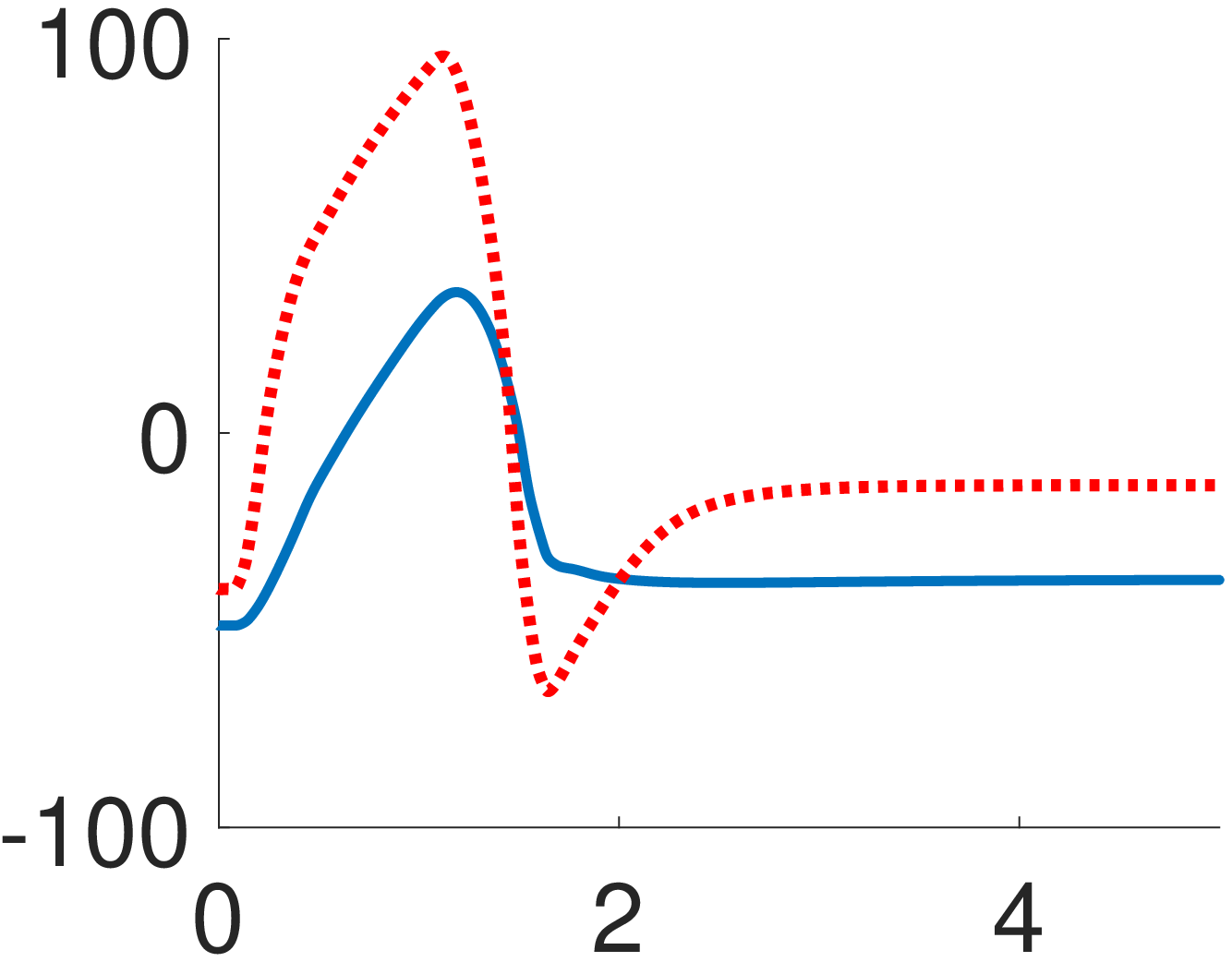} \hspace{0.1cm} 
\includegraphics[width=0.28\textwidth]{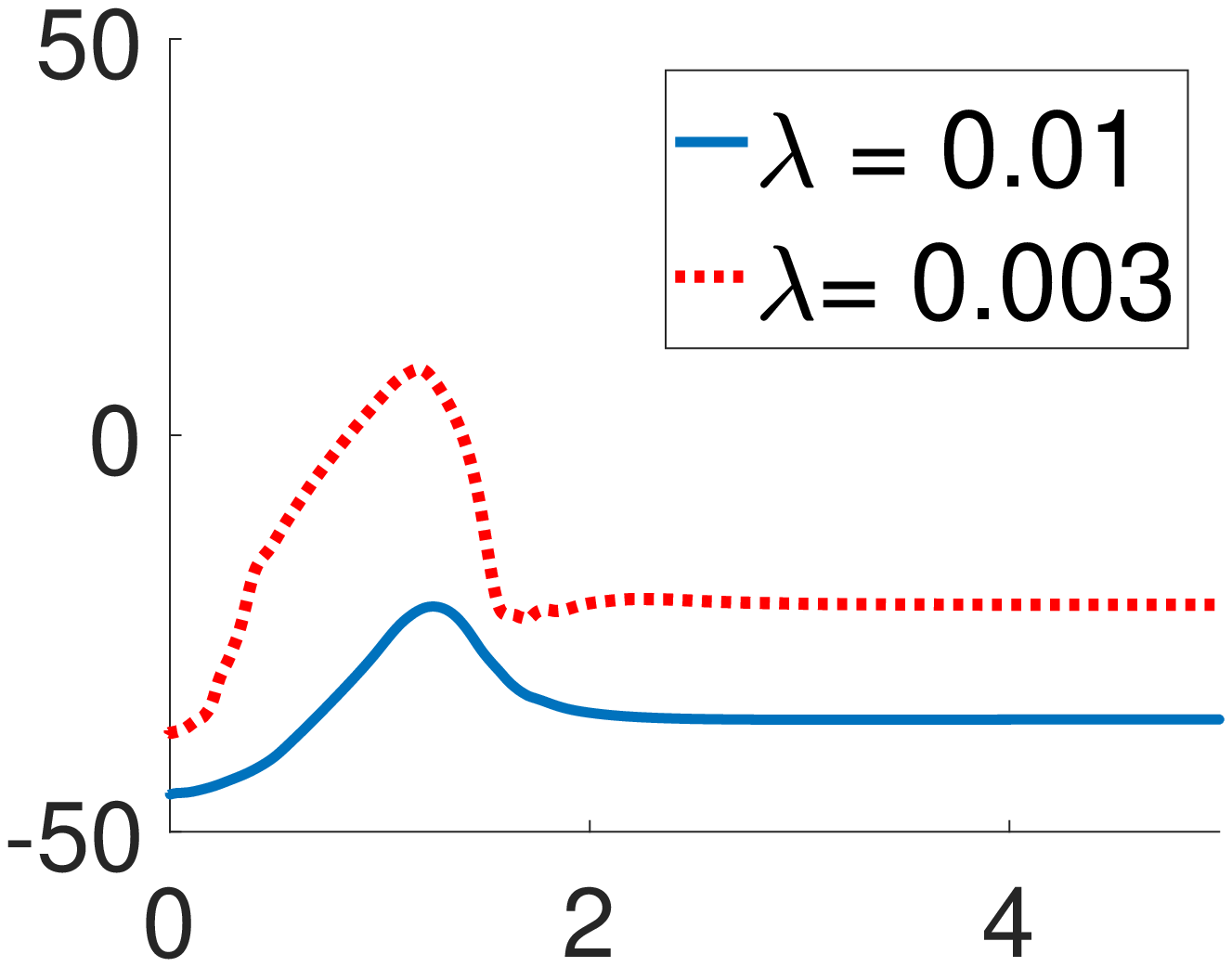}\\
 \hspace*{0.7cm} {\small time $t$ [h]}  \hspace*{1.5cm} {\small time $t$ [h]}  \hspace*{1.3cm} {\small time $t$ [h]} 
\end{minipage}
}
\end{tabular}
\caption{Gas pipeline network. Temporal evolution of solution resulting from input profile $u$, for two choices of friction factor~$\lambda$. Illustrated for ends of pipes $\omega_2$, $\omega_4$ and $\omega_7$ \textit{(from left to right)}. \textit{Top:} Density $\rho$. \textit{Bottom:} Mass flow~$A  m$.\label{fig:gl38-time}}
\end{figure}

In the second example, we consider the isothermal Euler equations in a friction-dominated regime relevant in the context of gas transport networks, cf. \cite{art:gaslib-2017,art:model-gasdistrib-2009,art:domschke-adj-based2015}. Each edge represents a cross-sectionally averaged pipe, where the reference values for density and mass flow are taken as $\rho_\star=1~[\mathrm{kg\,m^{-3}}]$ and $(Am)_\star= 1~[\mathrm{kg\,s^{-1}}]$. We use an isothermal law for the pressure $p$ $[\mathrm{Pa}]$  
\begin{align*}
p(\rho) &= RT \frac{\rho}{1-RT\alpha \rho} 
\end{align*}
with $T= 283~[\mathrm{K}]$,  $R = 518~[\mathrm{J} (\mathrm{kg \, K})^{-1}]$ and $\alpha = -3 \cdot 10^{-8}~[\mathrm{Pa}^{-1}]$, which implies the pressure potential $P(\rho) = RT \rho \log{\left( {(1-RT\alpha \rho)\rho_\star}/{\rho} \right)}$. 
The friction factor $\lambda$ is set constant over the whole network. With $\lambda = 0.01$ we choose a small, but typical value for gas networks. For illustration purposes we additionally consider $\lambda =0.003$ being much smaller than observed in gas network modeling.
The used pipeline network is visualized in Fig.~\ref{fig:net-topol-ex}. It consists of 38 pipes with diameters $D^\Onepipe$ between 0.4 and 1~$[\mathrm{m}]$ and lengths $l^\Onepipe$ between 5 to 74~$[\mathrm{km}]$. The total pipe length of the network is 1008~$[\mathrm{km}]$. The topology is a slight modification of \cite[GasLib-40]{art:gaslib-2017}, whereby the compressors are replaced by pipes. Consequently, there is no compensation of the energy loss related to friction inside of the network. For the realization of network simulations with active elements, a coupling of the sub-components over their boundary efforts and flows can be used. We refer to \cite{inproc:blsEcmi2016} for details and to \cite{code:bls22-phapprox}, where a simulation script for a network with a compressor is provided. The circular cross-sectional pipe areas  $A^\Onepipe$ act as edge weights for $\Onepipe \in \mathcal{E}$ (cf.\ Appendix~\ref{app:weighted-edges}).
At the six boundary nodes $\nu_i$, $i=1,...,6$ we prescribe the following boundary conditions 
\begin{align*}
	\rho(t,\nu_1) &= (65 + u(t))\,\rho_\star, \qquad \rho(t,\nu_2) = (50 + u(t))\,\rho_\star, \qquad   \rho(t,\nu_4) = (60 - u(t))\,\rho_\star, \\
	\rho(t,\nu_5) &= 60\,\rho_\star, \hspace*{2.cm} \rho(t,\nu_6) = 45\, \rho_\star, \hspace*{1.65cm}  A  m (t,\nu_3) = - 100\, (Am)_\star,
\end{align*}
with the time-varying input $u$,
 \begin{align*}
 	u(t) = \begin{cases} 10 \,{t}/{t_\star}, & \hspace*{0.5cm} 0 \le t<t_\star \\ 10\, (2-{t}/t_\star), & \hspace*{0.4cm} t_\star \le t< 1.5 t_\star  \\ 5, & 1.5 t_\star \le t, \end{cases} 
 \end{align*}
with reference time $t_\star=1~[\mathrm{h}]$.
As initial condition the stationary solution belonging to the boundary conditions at $t=0$ is taken. The end time is $T=5\,t_\star$. 

\begin{figure}[tb]
\centering
\begin{tabular}{rllrl}
\begin{minipage}{0.022\textwidth}
{\hspace{0.3cm}
{\small \rotatebox{90}{density $\rho$}\\ \vspace{0.0cm}
}}
\end{minipage}
&
{\hspace{-0.6cm}
\begin{minipage}{0.4\textwidth}
\center
\hspace*{0.00cm} \includegraphics[height = 0.68\textwidth, width = 0.85\textwidth]{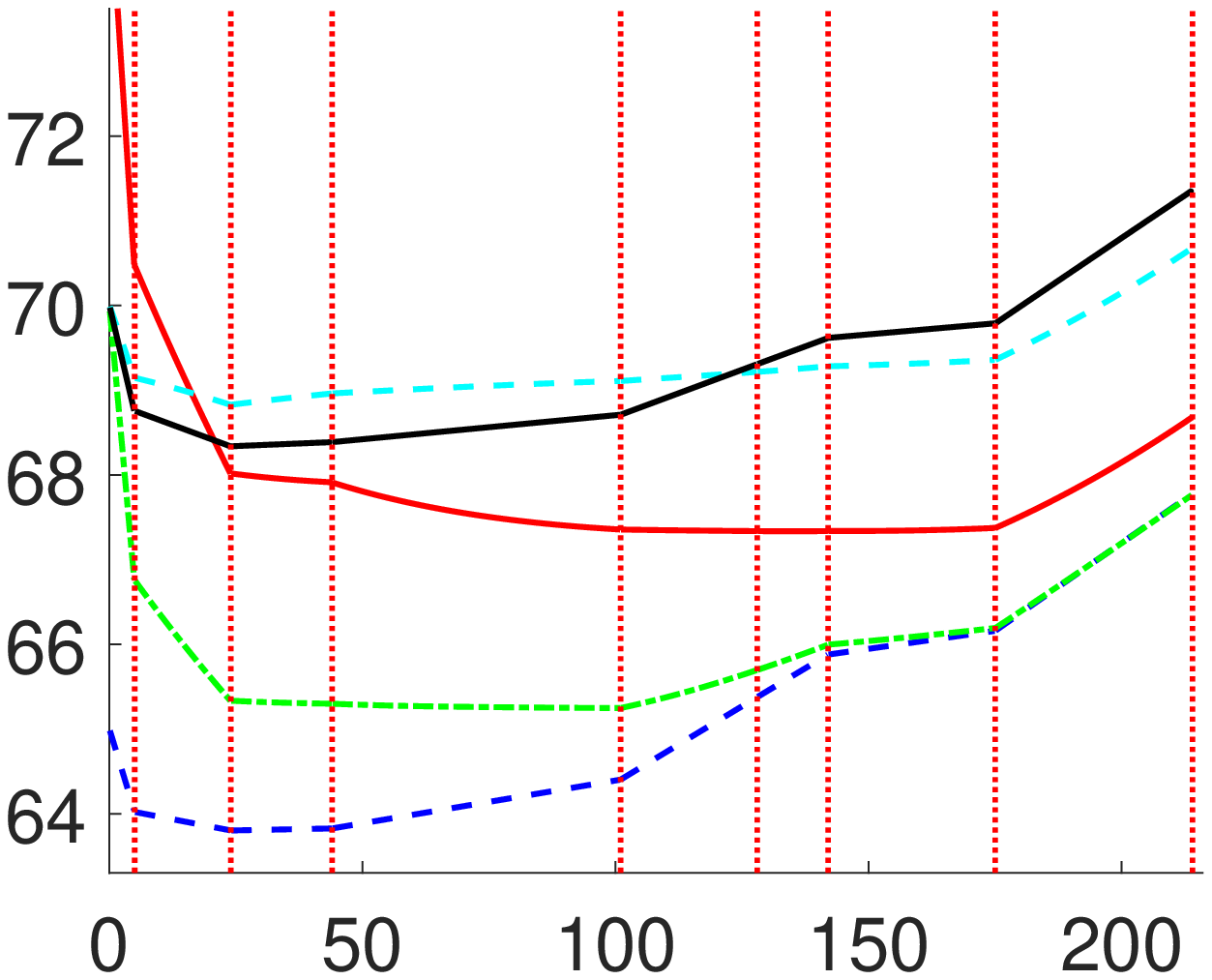}\\
 \hspace{0.5cm} {\small space $x$ [km]}
\end{minipage}
}
&
\begin{minipage}{0.00\textwidth}
\end{minipage}
&
\begin{minipage}{0.022\textwidth}
{\hspace{0.3cm}
{\small \vspace{0.0cm}\\
\rotatebox{90}{mass flow $A  m$\\ \vspace{0.cm}} 
}}
\end{minipage}
&
{\hspace{-0.6cm}
\begin{minipage}{0.42\textwidth}
\center
\includegraphics[height = 0.68\textwidth, width = 0.88\textwidth]{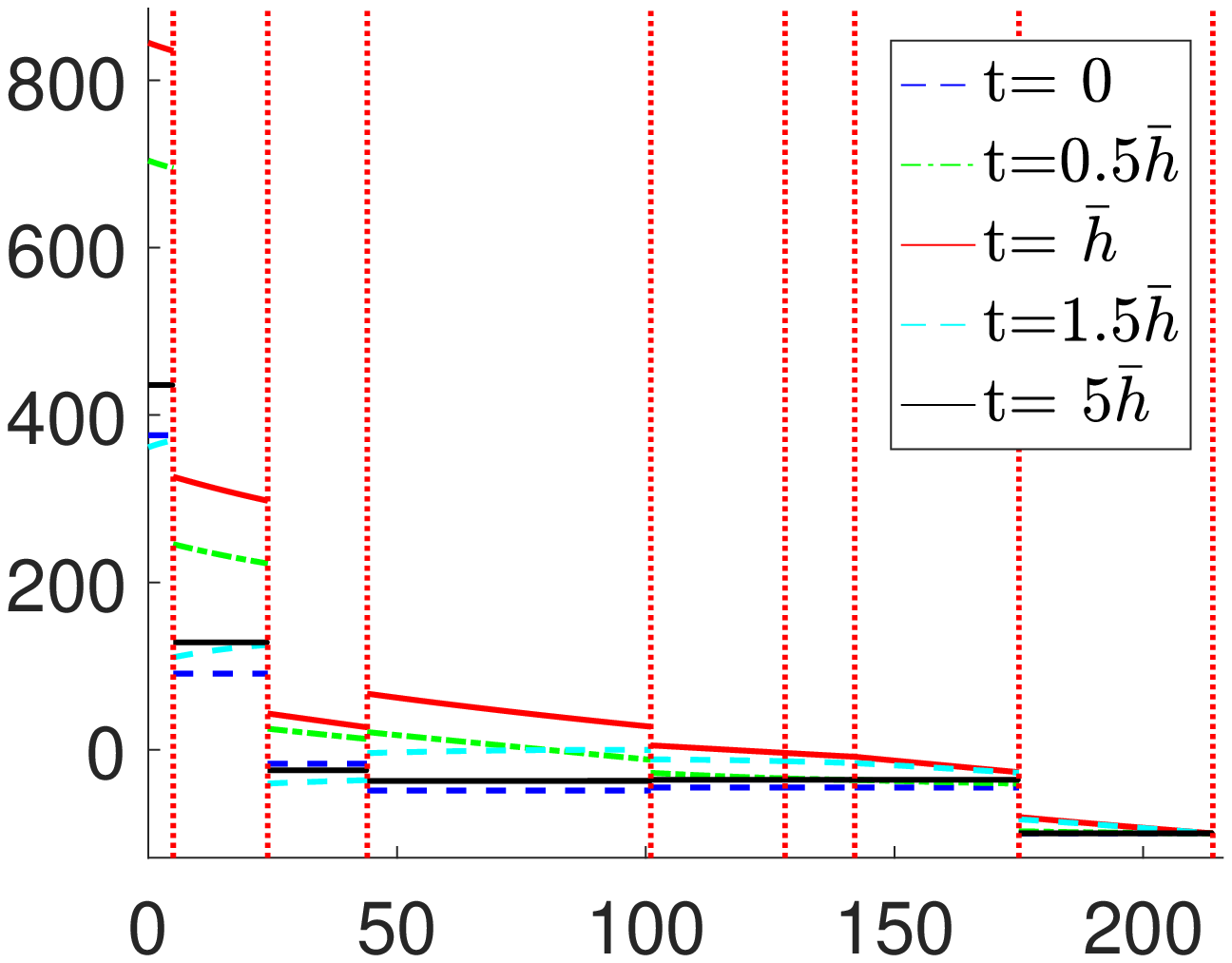}\\
 \hspace{0.5cm} {\small space $x$ [km]}
\end{minipage}
}
\end{tabular}
\caption{Gas pipeline network, friction factor $\lambda = 0.01$. Solution along the pipes $\omega_j$, $j=1,...,8$, path from $\nu_1$ to $\nu_3$ (cf.\ Fig.~\ref{fig:net-topol-ex}) for various times. Pipe junctions are indicated by vertical red lines.\label{fig:gl38-space}}
\end{figure}

The space discretization is performed by means of the lowest-order finite elements with a uniform mesh on each pipe. The maximal length of a finite element satisfies $\Delta_x \leq 200~[\mathrm{m}]$, which yields $10\,156$ degrees of freedom in space. The time step is chosen as  $\Delta_t = 1~[s]$, i.e., in total $18\,000$ steps for the simulation time of 5 hours. For this setup, the runtime on our machine is about 35 minutes per simulation. The computational cost is dominated by the nonlinear solves needed in each time step.
The resulting density and mass flow are presented along the pipes $\omega_j$, $j=1,...,8$, that build a path from the supplier node $\nu_1$ to the consumer node $\nu_3$ (see Fig.~\ref{fig:net-topol-ex}). To illustrate the impact of the damping on the temporal evolution, the solution at fixed space points along the path is plotted in Fig.~\ref{fig:gl38-time}. Stronger damping effects and lower peaks are observed for the larger, realistic friction factor ($\lambda=0.01$). Moreover, the sharp input profile $u$ is more clearly transferred to the mass flow, whereas the dynamics of the density is smoothed out stronger. 
The simulation time the moving profile needs to travel from the end of pipe $\omega_2$ to the end of pipe $\omega_7$ is about $0.5 t_\star$ for density and mass flow. For $t \geq 3t_\star$ the state seems to be almost at rest in the considered spatial domain for both choices of $\lambda$, whereby the equilibrium is reached slightly faster for the model with the application-relevant friction factor ($\lambda = 0.01$). For this model, we also visualize the solution in space in Fig.~\ref{fig:gl38-space}, using
time snapshots along the path. The pipe junctions are indicated by vertical lines. Note that the mass flow is discontinuous at junctions with more than two pipes, which is in accordance to the coupling conditions \eqref{bls-eq:abstr-coup}. As the spatial representation illustrates, no shocks are to be expected in the considered friction-dominated regime, and the solution has a rather simple structure in space. The latter motivates the use of model  reduction techniques to speed up gas network simulations, cf., \cite{phd:liljegren,art:mor-gas-grundel-jansen,inbook:mathEnergy21}, which in contrast to finite element methods rely on ansatz functions with global support over the full network. 

The results of this paper lay the theoretical foundation for generalizing our model order reduction approach \cite{art:lilsailer-dwemor} for the linear damped wave equation to a model order- and complexity-reduction approach for a general nonlinear flow problem class. In our follow-up paper \cite{art:bls21-snapshotbased}, the algorithmic aspects related to snapshot-based model reduction in our structure-preserving framework are investigated.

\section{Conclusion and outlook}

For many practically relevant dynamical systems, an underlying Hamiltonian structure can be identified. In this paper we proposed a structure-preserving approximation approach relying on a port-Hamiltonian formulation for a class of nonlinear flows on networks. Its basis is a variational principle that inherits the Hamiltonian structure. Further, the parametrization of the solution takes a prominent role. We take a mix of energy and co-energy variables and thus avoid to fully change to the co-energy variables (which are equal to the entropy variables for the Euler equations). To treat the parametrization and the nonlinearities in a systematic manner, we employ the theory on Legendre transformations. Apart from the port-Hamiltonian structure, our approximations ensure local mass conservation and an energy bound under mild assumptions. In particular, quite general Galerkin projections and complexity reduction of nonlinearities by a quadrature-type ansatz are covered by our analysis. We showcased the applicability and good stability properties of our approximation at the example of the barotropic Euler equations. However, our approach can also be applied to any p-system and  other symmetrizable hyperbolic systems in two variables, e.g., in the context of electromagnetic waves.

The results of this paper can also be used as a theoretical basis for structure-preserving model reduction methods, e.g., for gas network systems. Note 
that the realization of compatible reduced models involves some non-trivial algorithmic issues. In parts, they are similar to the ones occurring in symplectic model order reduction \cite{art:symplHamMor,art:morParamHamHesthaven,art:moramHamDissHesthaven}, some are covered in our former works \cite{phd:liljegren,inproc:bls-hyp2018,art:lilsailer-dwemor}. A thorough discussion of the algorithmic aspects related to snapshot-based model order- and complexity-reduction in our framework can be found in our follow-up paper \cite{art:bls21-snapshotbased}. 
 Another interesting direction for future research is the extension and adaption of our approach to other classes of model problems and more involved applications, e.g., to multi-dimensional systems with a similar underlying Hamiltonian structure  \cite{art:gieselmann-rel-energy,art:antonelli-quantumHydrodyn,inbook:ruggeri05-dissiphyp}, or thermal Euler equations with additional dissipation terms \cite{inproc:BadMBM18,inproc:pH-discrtictheat}.

\subsection*{\textbf{Data Availability}} The \texttt{MATLAB} code used to generate the presented numerical results and some supplementary tests, including a simulation of a network with a compressor, can be found under the DOI \texttt{10.5281/zenodo.6372667}, see \cite{code:bls22-phapprox}.

\begin{appendix}

\section{Inclusion of dissipation}  \label{subsec:dissip}
Our model problem \eqref{bls-eq:abstr} is of hyperbolic type. This changes when dissipation effects in form of a second-order derivative in space are included. We assume them to be described by a non-negative term $d: \mathbb{R}^2 \rightarrow \mathbb{R}^+$, similar as in \cite{art:egger-mfem-compressEuler}.
In generalization to \eqref{bls-eq:abstr}, let the state \mbox{$\ubar{\bv{z}}: [0,T] \times \blsPipe \rightarrow \mathbb{R}^2 $} be governed by 
\begin{align*} 
	\partial_t	\ubar{\bv{z}}(t,x) =
	\begin{bmatrix}
		 & -\partial_x \\
		 -\partial_x & -r(\ubar{\bv{z}}(t,x)	)
	\end{bmatrix}
	 \derst{\bv{z}} h(\ubar{\bv{z}}(t,x)) + 
	\begin{matrixbj}
		0 \\
		\partial_x \left( d(\ubar{\bv{z}}(t,x)) \partial_x \derst{2} h(\ubar{\bv{z}}(t,x)) \right)
	\end{matrixbj}
\end{align*}
with coupling conditions at $\blsVertex \in {\mathcal{N}}_0$ given by
\begin{align*} 
	\sum_{\Onepipe \in \mathcal{E}(\blsVertex) } n^{\Onepipe}[\blsVertex]  \nabla_{z_2} h(\ubar{\bv{z}}_{|\Onepipe}(t,\blsVertex)) = 0  , \qquad
	s(\ubar{\bv{z}}_{|\Onepipe}(t,\blsVertex)) = s(\ubar{\bv{z}}_{|\tilde{\Onepipe}}(t,\blsVertex)) \quad \text{for } \Onepipe,\tilde{\Onepipe} \in \mathcal{E}(\blsVertex),	
\end{align*}
where $s({\ubar{\bv{z}}}) = \nabla_{1} h({\ubar{\bv{z}}}) - d({\ubar{\bv{z}}}) \partial_x \derst{2} h({\ubar{\bv{z}}})$ for ${\ubar{\bv{z}}}=[ {z}_1;{z}_2] \in \funSpace{C}_{pw}^1(\mathcal{E})$.
To close the system, initial conditions and one boundary condition per boundary node $\blsVertex\in \mathcal{N}_\partial$ have to be prescribed, similarly to the case without dissipation.
Sufficiently smooth solutions can be shown to fulfill the energy dissipation equality
\begin{align*}
	\frac{d}{dt} \tilde{\HamPDE}(\ubar{\bv{z}}) &=  \hspace*{-0.25cm} \sum_{\blsVertex \in \mathcal{N}_\partial, \, \Onepipe \in \mathcal{E}(\blsVertex) } \hspace*{-0.5cm} n^{\Onepipe}[\blsVertex] s(\ubar{\bv{z}}_{|\Onepipe}[\blsVertex]) \derst{2} h(\ubar{\bv{z}}_{|\Onepipe}[\blsVertex]) -\Lscal r(\ubar{\bv{z}}), (\derst{2} h(\ubar{\bv{z}}))^2 \Rscal  - \Lscal  d(\ubar{\bv{z}}), (\partial_x \derst{2} h(\ubar{\bv{z}}))^2 \Rscal
	 \\ & \leq \sum_{\blsVertex \in \mathcal{N}_\partial, \, \Onepipe \in \mathcal{E}(\blsVertex) } n^{\Onepipe}[\blsVertex] s(\ubar{\bv{z}}_{|\Onepipe}[\blsVertex]) \derst{2} h(\ubar{\bv{z}}_{|\Onepipe}[\blsVertex]) .
\end{align*}
Our spatial approximation approach is transferable to this problem. Respective discrete energy bounds and port-Hamiltonian structure can be shown with very minor adjustments to our derivations from the main part.

\section{Edge weights}\label{app:weighted-edges}
Pipelines in gas transport networks are typically modeled with cross-sectionally averaged dynamics. Hence, we include the cross-sec{\-}tional pipe area $A^{\Onepipe}$ for $\Onepipe \in \mathcal{E}$ as edge weight in our approach. Note that the edge weighting affects all integral-related expressions and definitions from  Section~\ref{subsec-netdescrip}. The inner product $\Lscal \cdot, \cdot \Rscal$ becomes $\Lscal b, \tilde{b} \Rscal = \sum_{\Onepipe \in \mathcal{E}} A^{\Onepipe} \int_{\Onepipe} b[x] \tilde{b}[x] dx$ and the incidence mapping 
\begin{align*}
	 n^\Onepipe[\blsVertex] =
	 \begin{cases}
	 	\,\,\,\, A^{\Onepipe} & \text{ for } \Onepipe = ( \blsVertex, \bar{\blsVertex}) \text{ for some }  \bar{\blsVertex} \in \mathcal{N} \\
	 	 -A^{\Onepipe}  & \text{ for } \Onepipe = ( \bar{\blsVertex}, \blsVertex ) \text{ for some }  \bar{\blsVertex} \in \mathcal{N},
	 \end{cases}
\end{align*}
which modifies boundary and coupling conditions, and thus the boundary operator $\TraceOp^{ }: \funSpace{H}_{pw}^1(\mathcal{E}) \rightarrow \mathbb{R}^{p}$ and the function space $\funSpace{H}_{div}(\mathcal{E})$.
Moreover, the Hamiltonian is altered to $\HamPDE(\ubar{\bv{z}}) = \Lscal h(\ubar{\bv{z}}), 1 \Rscal = \sum_{\Onepipe \in \mathcal{E}} A^{\Onepipe} \int_{\Onepipe} h(\ubar{\bv{z}}) dx$.
\end{appendix}

\bibliographystyle{abbrv}
\bibliography{NLFlow}

\begin{thebibliography}{10}

\bibitem{book:matrixManifoldsAbsil08}
P.-A. Absil, R.~Mahony, and R.~Sepulchre.
\newblock {\em Optimization Algorithms on Matrix Manifolds}.
\newblock Princeton University Press, 2008.

\bibitem{art:morParamHamHesthaven}
B.~Afkham and J.~Hesthaven.
\newblock Structure preserving model reduction of parametric {H}amiltonian
  systems.
\newblock {\em SIAM J. Sci. Comput.}, 39(6):A2616--A2644, 2017.

\bibitem{art:moramHamDissHesthaven}
B.~Afkham and J.~Hesthaven.
\newblock Structure-preserving model-reduction of dissipative {H}amiltonian
  systems.
\newblock {\em J. Sci. Comput.}, 81(1):3--21, 2019.

\bibitem{art:efficient-integration-cubature}
S.~S. An, T.~Kim, and D.~L. James.
\newblock Optimizing cubature for efficient integration of subspace
  deformations.
\newblock {\em ACM Trans. Graph.}, 27(5):1--10, 2008.

\bibitem{art:antonelli-quantumHydrodyn}
P.~Antonelli and P.~Marcati.
\newblock The quantum hydrodynamics system in two space dimensions.
\newblock {\em Arch. Ration. Mech. Anal.}, 203(2):499--527, 2012.

\bibitem{book:arnold2006-comp}
D.~N. Arnold, P.~B. Bochev, R.~B. Lehoucq, R.~A. Nicolaides, and M.~Shashkov,
  editors.
\newblock {\em Compatible Spatial Discretizations}.
\newblock The IMA Volumes in Mathematics and its Applications. Springer, 1
  edition, 2006.

\bibitem{art:badlyan2018open}
A.~M. Badlyan, B.~Maschke, C.~A. Beattie, and V.~Mehrmann.
\newblock Open physical systems: from {GENERIC} to port-{H}amiltonian systems.
\newblock arXiv e-prints 1804.04064, 2018.

\bibitem{art:morBauBF14}
U.~Baur, P.~Benner, and L.~Feng.
\newblock Model order reduction for linear and nonlinear systems: A
  system-theoretic perspective.
\newblock {\em Arch. Comput. Methods Eng.}, 21(4):331--358, 2014.

\bibitem{art:beattie2017porthamiltonian}
C.~A. Beattie, V.~Mehrmann, H.~Xu, and H.~Zwart.
\newblock Linear port-{H}amiltonian descriptor systems.
\newblock {\em Math. Control Signals Systems}, 30(17), 2018.

\bibitem{book:dimred2003}
P.~Benner, V.~Mehrmann, and D.~C. Sorensen, editors.
\newblock {\em Dimension Reduction of Large-Scale Systems}.
\newblock Lecture Notes in Computational Science and Engineering. Springer, 1
  edition, 2005.

\bibitem{book:braess07}
D.~Braess.
\newblock {\em Finite {E}lements. Theory, {F}ast {S}olvers and {A}pplications
  in {E}lasticity {T}heory}.
\newblock Springer, 4 edition, 2007.

\bibitem{book:BrezzisFuncAna}
H.~Brezis.
\newblock {\em Functional Analysis, Sobolev Spaces and Partial Differential
  Equations}.
\newblock Springer, 2011.

\bibitem{art:hyperreduction-preserving-lagrangian-structure}
K.~Carlberg, R.~Tuminaro, and P.~Boggs.
\newblock Preserving {L}agrangian structure in nonlinear model reduction with
  application to structural dynamics.
\newblock {\em {SIAM} J. Sci. Comput.}, 37(2):B153--B184, 2015.

\bibitem{art:avfPDE}
E.~Celledoni, V.~Grimm, R.~I. McLachlan, D.~I. McLaren, D.~O'Neale, B.~Owren,
  and G.~R.~W. Quispel.
\newblock Preserving energy resp. dissipation in numerical {PDE}s using the
  {Average Vector Field} method.
\newblock {\em J. Comput. Phys.}, 231(20):6770--6789, 2012.

\bibitem{art:CharlotHughes}
F.~Chalot, T.~J. Hughes, and F.~Shakib.
\newblock Symmetrization of conservation laws with entropy for high-temperature
  hypersonic computations.
\newblock {\em Comput. Syst. Eng.}, 1(2-4):495--521, 1990.

\bibitem{art:chan-EntropROM}
J.~Chan.
\newblock Entropy stable reduced order modeling of nonlinear conservation laws.
\newblock {\em J. Comput. Phys.}, 423:109789, 2020.

\bibitem{art-mor-structure-preserving-nonlinear-pH}
S.~Chaturantabut, C.~Beattie, and S.~Gugercin.
\newblock Structure-preserving model reduction for nonlinear port-{H}amiltonian
  systems.
\newblock {\em {SIAM} J. Sci. Comput.}, 38(5):B837--B865, 2016.

\bibitem{art:topics-in-strpresdisc}
S.~H. Christiansen, H.~Z. Munthe-Kaas, and B.~Owren.
\newblock Topics in structure-preserving discretization.
\newblock {\em Acta Numerica}, 20:1--119, 2011.

\bibitem{inbook:mathEnergy21}
T.~Clees, A.~Baldin, P.~Benner, S.~Grundel, C.~Himpe, B.~Klaassen,
  F.~K{\"u}sters, N.~Marheineke, L.~Nikitina, I.~Nikitin, J.~Pade, N.~Stahl,
  C.~Strohm, C.~Tischendorf, and A.~Wirsen.
\newblock Math{E}nergy -- {M}athematical key technologies for evolving energy
  grids.
\newblock In S.~G{\"o}ttlich, M.~Herty, and A.~Milde, editors, {\em
  Mathematical Modeling, Simulation and Optimization for Power Engineering and
  Management}, pages 233--262. Springer, 2021.

\bibitem{art:domschke-adj-based2015}
P.~Domschke, O.~Kolb, and J.~Lang.
\newblock Adjoint-based error control for the simulation and optimization of
  gas and water supply networks.
\newblock {\em Appl. Math. Comput.}, 259:1003--1018, 2015.

\bibitem{art:egger-mfem-compressEuler}
H.~Egger.
\newblock A robust conservative mixed finite element method for compressible
  flow on pipe networks.
\newblock {\em {SIAM} J. Sci. Comput.}, 40(1):A108--A129, 2018.

\bibitem{art:discret-dissip-Egger19}
H.~Egger.
\newblock Structure preserving approximation of dissipative evolution problems.
\newblock {\em Numer. Math.}, 143(1):85--106, 2019.

\bibitem{art:egger2020stability}
H.~Egger and J.~Giesselmann.
\newblock Stability and asymptotic analysis for instationary gas transport via
  relative energy estimates.
\newblock arXiv e-prints 2012.14135, 2020.

\bibitem{art:kugler-dwe-net}
H.~Egger and T.~Kugler.
\newblock Damped wave systems on networks: Exponential stability and uniform
  approximations.
\newblock {\em Numer. Math.}, 138(4):839--867, 2018.

\bibitem{inproc:bls-hyp2018}
H.~Egger, T.~Kugler, and B.~Liljegren-Sailer.
\newblock Stability preserving approximations of a semilinear hyperbolic gas
  transport model.
\newblock In {\em Hyperbolic {P}roblems: {T}heory, {N}umerics, {A}pplications},
  volume~10, pages 427--433. AIMS Series on Appl. Math, 2020.

\bibitem{art:lilsailer-dwemor}
H.~Egger, T.~Kugler, B.~Liljegren-Sailer, N.~Marheineke, and V.~Mehrmann.
\newblock On structure-preserving model reduction for damped wave propagation
  in transport networks.
\newblock {\em SIAM J. Sci. Comput.}, 40(1):A331--A365, 2018.

\bibitem{art:comred-ecsw}
C.~Farhat, P.~Avery, T.~Chapman, and J.~Cortial.
\newblock Dimensional reduction of nonlinear finite element dynamic models with
  finite rotations and energy-based mesh sampling and weighting for
  computational efficiency.
\newblock {\em Int. J. Numer. Meth. Eng.}, 98(9):625--662, 2014.

\bibitem{inproc:farle-pH-FE}
O.~{Farle}, D.~{Klis}, M.~{Jochum}, O.~{Floch}, and R.~{Dyczij-Edlinger}.
\newblock A port-{H}amiltonian {F}inite-{E}lement formulation for the {M}axwell
  equations.
\newblock In {\em 2013 International Conference on Electromagnetics in Advanced
  Applications (ICEAA)}, pages 324--327, 2013.

\bibitem{art:FISHER13-entrop}
T.~C. Fisher and M.~H. Carpenter.
\newblock High-order entropy stable finite difference schemes for nonlinear
  conservation laws: Finite domains.
\newblock {\em J. Comput. Phys.}, 252:518--557, 2013.

\bibitem{art:gieselmann-rel-energy}
J.~Giesselmann, C.~Lattanzio, and A.~E. Tzavaras.
\newblock Relative energy for the {K}orteweg theory and related {H}amiltonian
  flows in gas dynamics.
\newblock {\em Arch. Ration. Mech. Anal.}, 223(3):1427--1484, 2017.

\bibitem{art:GoloTSM04}
G.~Golo, V.~Talasila, A.~van~der Schaft, and B.~Maschke.
\newblock {H}amiltonian discretization of boundary control systems.
\newblock {\em Autom.}, 40(5):757--771, 2004.

\bibitem{art:mor-gas-grundel-jansen}
S.~Grundel, L.~Jansen, N.~Hornung, T.~Clees, C.~Tischendorf, and P.~Benner.
\newblock Model order reduction of differential algebraic equations arising
  from the simulation of gas transport networks.
\newblock In S.~Sch{\"o}ps, A.~Bartel, M.~G{\"u}nther, W.~E.~J. ter Maten, and
  C.~P. M{\"u}ller, editors, {\em Progress in Differential-Algebraic Equations:
  Deskriptor 2013}, pages 183--205. Springer, 2014.

\bibitem{art:gugat17noblowup}
M.~Gugat and S.~Ulbrich.
\newblock The isothermal euler equations for ideal gas with source term:
  Product solutions, flow reversal and no blow up.
\newblock {\em J. Math. Anal. Appl.}, 454(1):439--452, 2017.

\bibitem{inbook:interpolation-based-port-Hamiltonian-Systems}
S.~Gugercin, R.~V. Polyuga, C.~A. Beattie, and A.~van~der Schaft.
\newblock Interpolation-based {$H_2$} model reduction for port-{H}amiltonian
  systems.
\newblock In {\em Proceedings of the 48th IEEE Conference on Decision and
  Control, and the 28th Chinese Control Conference, Shanghai}, pages
  5362--5369, 2009.

\bibitem{book:hairerGeomInt}
E.~Hairer, C.~Lubich, and G.~Wanner.
\newblock {\em Geometric Numerical Integration}.
\newblock Springer, 2 edition, 2006.

\bibitem{art:harten1983}
A.~Harten.
\newblock On the symmetric form of systems of conservation laws with entropy.
\newblock {\em J. Comput. Phys.}, 49(1):151--164, 1983.

\bibitem{inproc:pH-discrtictheat}
S.-A. Hauschild, N.~Marheineke, V.~Mehrmann, J.~Mohring, A.~M. Badlyan,
  M.~Rein, and M.~Schmidt.
\newblock Port-{H}amiltonian modeling of district heating networks.
\newblock In S.~Grundel, T.~Reis, and S.~Sch{\"o}ps, editors, {\em Progress in
  Differential-Algebraic Equations II}, pages 333--355. Springer, 2020.

\bibitem{art:model-gasdistrib-2009}
A.~Herran-Gonzalez, J.~M. D.~L. Cruz, B.~D. Andres-Toro, and J.~L.
  Risco-Martin.
\newblock Modeling and simulation of a gas distribution pipeline network.
\newblock {\em Appl. Math. Model.}, 33(3):1584--1600, 2009.

\bibitem{art:Jameson2008}
A.~Jameson.
\newblock The construction of discretely conservative finite volume schemes
  that also globally conserve energy or entropy.
\newblock {\em J. Sci. Comput.}, 34(2):152--187, 2008.

\bibitem{art:weak-pH-koty}
P.~Kotyczka, B.~Maschke, and L.~Lefevre.
\newblock Weak form of {S}tokes--{D}irac structures and geometric
  discretization of port-{H}amiltonian systems.
\newblock {\em J. Comput. Phys.}, 361:442--476, 2018.

\bibitem{art:lee-mixededmimFE}
D.~Lee and A.~Palha.
\newblock A mixed mimetic spectral element model of the rotating shallow water
  equations on the cubed sphere.
\newblock {\em J. Comput. Phys.}, 375:240--262, 2018.

\bibitem{book:leveque-finite-volume-methods}
R.~J. LeVeque.
\newblock {\em Finite Volume Methods for Hyperbolic Problems}.
\newblock Cambridge University Press, 2002.

\bibitem{phd:liljegren}
B.~Liljegren-Sailer.
\newblock {\em On {P}ort-{H}amiltonian {M}odeling and {S}tructure-{P}reserving
  {M}odel {R}eduction}.
\newblock PhD thesis, Universit{\"a}t Trier, 2020.

\bibitem{code:bls22-phapprox}
B.~Liljegren-Sailer.
\newblock Code for the paper {'On port-Hamiltonian approximation of a nonlinear
  flow problem on networks'}.
\newblock https://doi.org/10.5281/zenodo.6372667, 2022.

\bibitem{inproc:blsEcmi2016}
B.~Liljegren-Sailer and N.~Marheineke.
\newblock A structure-preserving model order reduction approach for
  space-discrete gas networks with active elements.
\newblock In {\em Progress in Industrial Mathematics at ECMI 2016}, pages
  439--446. Springer, 2017.

\bibitem{art:bls21-snapshotbased}
B.~Liljegren-Sailer and N.~Marheineke.
\newblock On snapshot-based model reduction under compatibility conditions for
  a nonlinear flow problem on networks.
\newblock arXiv e-prints 2110.04777, 2021.

\bibitem{art:pcH-maschke92}
B.~M. Maschke and A.~van~der Schaft.
\newblock Port-controlled {H}amiltonian systems: Modelling origins and
  systemtheoretic properties.
\newblock {\em IFAC Proceedings Volumes}, 25(13):359--365, 1992.

\bibitem{inbook:Maschke2001}
B.~M. Maschke and A.~van~der Schaft.
\newblock {H}amiltonian representation of distributed parameter systems with
  boundary energy flow.
\newblock In {\em Nonlinear Control in the Year 2000}, volume~2, pages
  137--142. Springer, 2001.

\bibitem{book:diffgeo-mcInerney}
A.~McInerney.
\newblock {\em First Steps in Differential Geometry: Riemannian, Contact,
  Symplectic}.
\newblock Johns Hopkins series in information sciences and systems. Springer, 1
  edition, 2013.

\bibitem{inprod-MehrmannM19}
V.~Mehrmann and R.~Morandin.
\newblock Structure-preserving discretization for port-{H}amiltonian descriptor
  systems.
\newblock In {\em 58th {IEEE} Conference on Decision and Control, {CDC} 2019,
  Nice, France, December 11-13, 2019}, pages 6863--6868. {IEEE}, 2019.

\bibitem{art:mei-psystems09}
M.~Mei.
\newblock Nonlinear diffusion waves for hyperbolic p-system with nonlinear
  damping.
\newblock {\em J. Differ. Eqs.}, 247(4):1275--1296, 2009.

\bibitem{art:hierarchGasMindt2019}
P.~Mindt, J.~Lang, and P.~Domschke.
\newblock Entropy-preserving coupling of hierarchical gas models.
\newblock {\em SIAM J. Math. Anal.}, 51(6):4754--4775, 2019.

\bibitem{art:mock1980}
M.~S. Mock.
\newblock Systems of conservation laws of mixed type.
\newblock {\em J. Differ. Eqs.}, 37(1):70--88, 1980.

\bibitem{inproc:BadMBM18}
A.~Moses~Badlyan, B.~Maschke, C.~Beattie, and V.~Mehrmann.
\newblock Open physical systems: {F}rom {GENERIC} to port-{H}amiltonian
  systems.
\newblock In {\em Proceedings of the 23rd International Symposium on
  Mathematical Theory of Systems and Networks}, pages 204--211, 2018.

\bibitem{book:NovotnyCompressible}
A.~Novotny and I.~Straksraba.
\newblock {\em Introduction to the Mathematical Theory of Compressible Flow.},
  volume~27 of {\em Johns Hopkins series in information sciences and systems}.
\newblock Oxford Lecture Series in Mathematics and its Applications, 1 edition,
  2004.

\bibitem{art:PasumarthyAS12}
R.~Pasumarthy, V.~Ambati, and A.~van~der Schaft.
\newblock Port-{H}amiltonian discretization for open channel flows.
\newblock {\em Syst. Control. Lett.}, 61(9):950--958, 2012.

\bibitem{art:symplHamMor}
L.~Peng and K.~Mohseni.
\newblock Symplectic model reduction of {H}amiltonian systems.
\newblock {\em SIAM J. Sci. Comput.}, 38(1):A1--A27, 2016.

\bibitem{art:euler-reigstad}
G.~A. Reigstad.
\newblock Existence and uniqueness of solutions to the generalized {R}iemann
  problem for isentropic flow.
\newblock {\em SIAM J. Appl. Math.}, 75(2):679--702, 2015.

\bibitem{book:RockWets98}
R.~Rockafellar and R.~Wets.
\newblock {\em Variational Analysis}.
\newblock Springer, 1998.

\bibitem{rockafellar-1970a}
R.~T. Rockafellar.
\newblock {\em Convex {A}nalysis}.
\newblock Princeton Mathematical Series. Princeton University Press, 1970.

\bibitem{inbook:ruggeri05-dissiphyp}
T.~Ruggeri.
\newblock Global existence of smooth solutions and stability of the constant
  state for dissipative hyperbolic systems with applications to extended
  thermodynamics.
\newblock In S.~Rionero and G.~Romano, editors, {\em Trends and Applications of
  Mathematics to Mechanics}, pages 215--224. Springer, 2005.

\bibitem{art:gaslib-2017}
M.~Schmidt, D.~A{\ss}mann, R.~Burlacu, J.~Humpola, I.~Joormann, N.~Kanelakis,
  T.~Koch, D.~Oucherif, M.~Pfetsch, L.~Schewe, R.~Schwarz, and M.~Sirvent.
\newblock {G}as{L}ib -- {A} {L}ibrary of {G}as {N}etwork {I}nstances.
\newblock {\em Data}, 2(4), 2017.

\bibitem{art:Nordstrom14-SumByParts}
M.~Sv{\"a}rd and J.~Nordstr{\"o}m.
\newblock Review of summation-by-parts schemes for initial--boundary-value
  problems.
\newblock {\em J. Comput. Phys.}, 268:17--38, 2014.

\bibitem{art:SchJ14}
A.~van~der Schaft and D.~Jeltsema.
\newblock Port-{H}amiltonian systems theory: An introductory overview.
\newblock {\em Foundations and Trends in Systems and Control}, 1:173--378,
  2014.

\bibitem{art:SchM13}
A.~van~der Schaft and B.~M. Maschke.
\newblock Port-{H}amiltonian systems on graphs.
\newblock {\em SIAM J. Control Optim.}, 51:906--937, 2013.

\bibitem{art:winters-euler-entropy}
A.~R. Winters, C.~Czernik, M.~B. Schily, and G.~J. Gassner.
\newblock Entropy stable numerical approximations for the isothermal and
  polytropic {E}uler equations.
\newblock {\em {BIT}}, 60:791--824, 2020.

\bibitem{art:WOLF2010401}
T.~Wolf, B.~Lohmann, R.~Eid, and P.~Kotyczka.
\newblock Passivity and structure preserving order reduction of linear
  port-{H}amiltonian systems using {K}rylov subspaces.
\newblock {\em Eur. J. Control}, 16(4):401--406, 2010.

\bibitem{book:Zeidler1984NonlinearFA}
E.~Zeidler.
\newblock {\em Nonlinear Functional Analysis and its Applications}, volume~3.
\newblock Springer, 1985.

\end{thebibliography}

\end{document}